\documentclass[a4paper,10pt,reqno]{amsart}

\usepackage[top=2.5cm, bottom=2.5cm, left=2.5cm, right=2.5cm]{geometry}
\usepackage{amssymb,amsmath,amsthm,ascmac,array,bm,multirow}
\usepackage[dvipdfmx]{graphicx}
\usepackage{xcolor}

\usepackage[dvipdfmx,colorlinks=true,bookmarks=true,%citecolor=blue,
bookmarksnumbered=true,bookmarkstype=toc,linktocpage=true
]{}

\allowdisplaybreaks
\numberwithin{equation}{section}

\newtheorem{thm}{Theorem}[section]

\newtheorem{lem}[thm]{Lemma}
\newtheorem{prop}[thm]{Proposition}
\newtheorem{cor}[thm]{Corollary}

\newtheorem{rem}[thm]{Remark}
\newtheorem{ass}[thm]{Assumption}

\newcommand{\mcc}{\mathcal{C}}
\newcommand{\mcd}{\mathcal{D}}
\newcommand{\mce}{\mathcal{E}}
\newcommand{\mcf}{\mathcal{F}}

\newcommand{\mcl}{\mathcal{L}}
\newcommand{\mcm}{\mathcal{M}}
\newcommand{\mcq}{\mathcal{Q}}

\newcommand{\mcu}{\mathcal{U}}
\newcommand{\mfb}{\mathfrak{b}}

\newcommand{\mfF}{\mathfrak{F}}

\newcommand{\mbbh}{\mathbb{H}}
\newcommand{\mbbn}{\mathbb{N}}
\newcommand{\mbbr}{\mathbb{R}}
\newcommand{\mbbrp}{\mathbb{R}_{+}}
\newcommand{\mbbs}{\mathbb{S}}

\newcommand{\mbby}{\mathbb{Y}}
\newcommand{\mbbz}{\mathbb{Z}}
\newcommand{\mbbzp}{\mathbb{Z}_{+}}

\newcommand{\al}{\alpha}
\newcommand{\del}{\delta}

\newcommand{\ep}{\epsilon}
\newcommand{\D}{\Delta}
\newcommand{\Sig}{\Sigma}
\newcommand{\lam}{\lambda}

\newcommand{\gam}{\gamma}
\newcommand{\Gam}{\Gamma}
\newcommand{\p}{\partial}

\newcommand{\cil}{\xrightarrow{\mcl}} % <- Convergence in law
\newcommand{\cip}{\xrightarrow{p}} % <- Convergence in probability
\newcommand{\argmin}{\mathop{\rm argmin}} 
\newcommand{\argmax}{\mathop{\rm argmax}}
\newcommand{\diag}{\mathop{\rm diag}} 
\newcommand{\trace}{\mathop{\rm trace}} 

\def\nn{\nonumber}

\def\wp{Wiener process}
\def\lp{L\'{e}vy process}
\def\lm{L\'{e}vy measure}

\def\sumj{\sum_{j=1}^{n}}

\def\pr{P}
\def\E{E}

\def\dim{\mathrm{dim}}

\def\tz{\theta_{0}}

\def\tes{\hat{\theta}_{n}}

\def\aes{\hat{\alpha}_{n}}

\def\ges{\hat{\gamma}_{n}}

\def\gqaic{\mathrm{GQAIC}}
\def\gqbic{\mathrm{GQBIC}}

\newcommand{\X}{\bm{X}}

\title[Gaussian quasi-information criteria for ergodic L\'{e}vy driven SDE]
{Gaussian quasi-information criteria for ergodic L\'{e}vy driven SDE}

\author[S. Eguchi]{Shoichi Eguchi}
\address[S. Eguchi]{Faculty of Information Science and Technology, Osaka Institute of Technology, 1-79-1 Kitayama, Hirakata City, Osaka, 573-0196, Japan.}
\email{shoichi.eguchi@oit.ac.jp}

\author[H. Masuda]{Hiroki Masuda}
\address[H. Masuda]{Faculty of Mathematics, Kyushu University, 744 Motooka Nishi-ku Fukuoka 819-0395, Japan
\and 
Graduate School of Mathematical Sciences, The University of Tokyo, 3-8-1 Komaba Meguro-ku Tokyo 153-8914, Japan.
}
\email{hmasuda@ms.u-tokyo.ac.jp}

\date{\today}

\begin{document}
\setlength{\baselineskip}{4.5mm}

\maketitle

\begin{abstract}
We consider relative model comparison for the parametric coefficients of an ergodic L\'{e}vy driven model observed at high-frequency. Our asymptotics is based on the fully explicit two-stage Gaussian quasi-likelihood function (GQLF) of the Euler-approximation type. For selections of the scale and drift coefficients, we propose explicit Gaussian quasi-AIC (GQAIC) and Gaussian quasi-BIC (GQBIC) statistics through the stepwise inference procedure, and prove their asymptotic properties.
In particular, we show that the mixed-rates structure of the joint GQLF, which does not emerge in the case of diffusions, gives rise to the non-standard forms of the regularization terms in the selection of the scale coefficient, quantitatively clarifying the relation between estimation precision and sampling frequency.
Also shown is that the stepwise strategies are essential for both the tractable forms of the regularization terms and the derivation of the asymptotic properties of the Gaussian quasi-information criteria.
Numerical experiments are given to illustrate our theoretical findings.
\end{abstract}

%%%%%
%%%%%
\section{Introduction}

Suppose that we observe an equally spaced high-frequency sample $\X_{n}=(X_{t_{j}^{n}})_{j=0}^{n}$ for $t_{j}^{n}=t_j=jh$, where $X=(X_{t})_{t\in\mbbr_{+}}$ is a solution to the stochastic differential equation (SDE)
\begin{align}
dX_{t} = A(X_{t})dt + C(X_{t-})dZ_{t},
\label{se:model1}
\end{align}
where $A:\,\mbbr^d\to\mbbr^d$ and $C:\,\mbbr^d\to\mbbr^{d}\otimes\mbbr^{r}$, and $Z=(Z_t)_{t\in\mbbrp}$ is an $r$-dimensional L\'{e}vy process independent of the initial value $X_{0}$. We will suppose that $Z$ is standardized in the sense that $Z_1$ is zero-mean and has the identity covariance matrix.
The sampling stepsize $h=h_n>0$ is a known real such that
\begin{align*}
T_{n}:=nh\to\infty, \quad nh^{2}\to 0
\end{align*}
as $n\to\infty$. We want to infer the coefficients $A$ and $C$ based on a sample $\X_n$, without specifying $\mcl(Z)$, the distribution of the process $Z$.
Also, suppose that we are given the following candidates
\begin{align}
& c_{1}(x,\gamma_{1}),\ldots,c_{M_{1}}(x,\gamma_{M_{1}}),
%\label{se:ms.c} 
\nn\\
& a_{1}(x,\alpha_{1}),\dots,a_{M_{2}}(x,\alpha_{M_{2}}), \nn
%\label{se:ms.a} 
\end{align}
for the scale and drift coefficients, respectively.
Then, each candidate SDE model $\mathcal{M}_{m_{1},m_{2}}$ is described by
\begin{equation}
dX_t=c_{m_{1}}(X_{t-},\gam_{m_{1}})dZ_t + a_{m_{2}}(X_t,\al_{m_{2}})dt.
\label{hm:SDE.cand.model}
\end{equation}
%Then, each candidate model $\mathcal{M}_{m_{1},m_{2}}$ is the L\'{e}vy driven SDE
%\begin{equation*}
%dX_{t} = a_{m_{2}}(X_{t},\alpha_{m_{2}})dt + c_{m_{1}}(X_{t-},\gamma_{m_{1}})dZ_{t},
%\end{equation*}
The distribution of $\X_n$ is seldom explicitly given, hence we need to resort to some approximation.
In this paper, we will consider the Gaussian approximation
%stochastic differential process model defined as
\begin{equation}
\mcl(X_{t_{j}}|X_{t_{j-1}}=x)\sim N\left(x+ha_{m_{2}}(x,\alpha_{m_{2}}),hc_{m_{1}}(x,\gamma_{m_{1}})c_{m_{1}}(x,\gamma_{m_{1}})^\top\right)
\label{hm:SDE.stat.model}
\end{equation}
as our statistical model corresponding to $\mathcal{M}_{m_{1},m_{2}}$,
where $\gamma_{m_{1}}\in\Theta_{\gamma_{m_{1}}}\subset\mathbb{R}^{p_{\gamma_{m_{1}}}}$ ($m_{1}=1,\dots,M_{1}$) and $\alpha_{m_{2}}\in\Theta_{\alpha_{m_{2}}}\subset\mbbr^{p_{\alpha_{m_{2}}}}$ ($m_{2}=1,\dots,M_{2}$) are finite-dimensional unknown parameters. 
The parameter spaces $\Theta_{\gamma_{m_{1}}}$ and $\Theta_{\alpha_{m_{2}}}$ are assumed to be bounded convex domains.
 main objective of this paper is to develop a model selection procedure for selecting the best model $\mcm_{\hat{m}_{1,n},\hat{m}_{2,n}}$ among the candidate models.
For selecting an appropriate model, we will develop
the Akaike information criterion (AIC, \cite{Aka73}) and Bayesian information criterion (BIC, \cite{Sch78}) type model comparison for semiparametric L\'{e}vy driven SDE \eqref{se:model1}.
Although we are interested in the SDE model \eqref{hm:SDE.cand.model}, we consider the Gaussian (logarithmic) quasi-likelihood function (GQLF) based on \eqref{hm:SDE.stat.model} instead of the true likelihood for the inference. 
In this sense, our statistical models are all misspecified.
%Some backgrounds are described below.
%, with the parameter spaces $\Theta_{\gamma_{m_{1}}}$ and $\Theta_{\alpha_{m_{2}}}$ being bounded convex domains.
%\tcr{
%To derive the criterion for selections of the scale and drift coefficients, we assume that candidate coefficients $c_{1},\ldots,c_{M_{1}}$ and $a_{1},\ldots,a_{M_{2}}$ are correctly specified in the sense that there exist $\gamma_{m_{1},0}\in\Theta_{\gamma_{m_{1}}}$ and $\alpha_{m_{2},0}\in\Theta_{\alpha_{m_{2}}}$ such that $c_{m_{1}}(\cdot,\gamma_{m_{1},0})=C(\cdot)$ and $a_{m_{2}}(\cdot,\alpha_{m_{2},0})=A(\cdot)$.
%}

The information criteria are one of the most convenient and powerful tools for model selection, and the AIC and BIC are derived from two different classical principles:
the AIC and the GIC (generalized information criterion \cite{KonKit96}), which is an extension of AIC, are predictive model selection criteria minimizing the Kullback-Leibler divergence which measures the deviation from the true model to the prediction model; the BIC is given by the Bayesian principle and used for finding better model descriptions.
The AIC is not intended to select the true model consistently even if the true model is included in the set of candidate models, while the BIC puts importance on both underfitting and overfitting.
Based on the same classical principles as AIC and BIC, several studies have been conducted on model selection for SDEs: the contrast-based information criterion for ergodic diffusion processes \cite{Uch10}, the BIC type information criterion for locally asymptotically quadratic models \cite{EguMas18a}, and the BIC type information criterion for possibly misspecified ergodic SDEs \cite{EguUeh21}.

%In this paper, we consider the L\'{e}vy driven SDEs as underlying stochastic processes, for which we have the same bottleneck as in the diffusion case: the genuine log-likelihood is rarely available in an explicit and tractable form.
Asymptotic inference based on the GQLF for L\'{e}vy driven SDE has been developed by several previous works, of course at the expense of asymptotic estimation efficiency: see \cite{Mas13} and \cite{MasUeh17-2}, %and \cite{Ueh19}, 
as well as the references therein.
The \textit{Gaussian quasi-AIC(BIC)}, which we will introduce and term \textit{GQAIC (GQBIC)} for short, is based on the fully explicit two-stage GQLF of the Euler-approximation type \eqref{hm:SDE.stat.model}.
Our study develops the GQAIC (GQBIC) under the condition that the scale and drift coefficients are correctly specified and clarifies that taking the two steps will be inevitable for the simple form of the GQAIC (GQBIC) to be in force in the sense that they provide us with the specific asymptotic selection probabilities for GQAIC and the selection consistency for GQBIC:
the details will be given in  Sections \ref{hm:sec_GQAIC}, \ref{hm:sec_GQBIC}, and \ref{se:sec_modsele}.
%This result is different from the case of diffusions which can obtain the information criteria that work for both joint and stepwise GQLF.
The GQAIC and GQBIC which are derived through the GQLF of the first stage will be called $\gqaic_{1}$ and $\gqbic_{1}$, respectively.
Also, the GQAIC and GQBIC based on the GQLF of the second stage will be called $\gqaic_{2}$ and $\gqbic_{2}$, respectively.

The two-stage procedure proposed in this paper is summarized as follows:
for the AIC type, first, we select a scale-coefficient model as a minimizer of $\gqaic_{1}$ over the candidates $c_{1},\ldots,c_{M_{1}}$, and then select a drift-coefficient model as a minimizer of $\gqaic_{2}$ over the candidates $a_{1}\dots,a_{M_{2}}$.
Our model comparisons will be presented in Section \ref{se:sec_modsele}.
There, we will consider the cases where the candidate coefficients $c_{1},\ldots,c_{M_{1}}$ and $a_{1},\ldots,a_{M_{2}}$ contain both correctly specified coefficients and misspecified coefficients. 
Also, we formally use the $\gqaic_{1,n}$ and $\gqaic_{2,n}$ even for the possibly misspecified coefficients, although the assumption that the candidate scale and drift coefficients are correctly specified is necessary for the derivation of GQAIC.
As for the BIC-type, we follow the same way, replacing $\gqaic_{1}$ and $\gqaic_{2}$ by $\gqbic_{1}$ and $\gqbic_{2}$, respectively.
In particular, concerned with both AIC- and BIC-type model comparisons of the scale coefficient, it turned out that we should employ some non-standard forms of the regularization term.
Especially for the BIC type methodology, 
it turned out that
%, to execute \textit{consistent} model selection, 
the conventional stochastic expansion of the marginal (quasi-)likelihood is not appropriate for \textit{consistent} model selection: we needed to ``heated up'' free energy (Section \ref{hm:sec_GQBIC.scale}).

The very different features compared with the ergodic diffusions will be presented in this paper. They are
essentially due to the mixed-rates structure \cite{Rad08} of the joint GQLF given by \eqref{hm:def_joint.GQLF} below, which does not emerge for the case of diffusions where $Z$ is a standard {\wp}; 
%Remark \ref{hm:rem_mixed.rates.M}
%See Remark \ref{hm:rem_diffusion+1}
see Remark \ref{hm:rem_diffusion}(1).
Informally speaking, the use of the GQLF against non-Gaussian {\lp es} causes the non-standard phenomena in inference for the scale coefficient, quantitatively clarifying the relation between estimation precision and sampling frequency.
Remarkably, we could still obtain explicit results of practical value.

In the rest of this paper, we give some prerequisites in Section \ref{hm:sec_pre}.
Then, in Sections \ref{hm:sec_GQAIC} and \ref{hm:sec_GQBIC} we present how the classical AIC- and BIC-type arguments can work in our model setup, respectively.
In Section \ref{se:sec_modsele}, we introduce the stepwise model comparison procedure and discuss the asymptotic probability of relative model selection.
Section \ref{hm:sec_sim} presents illustrative numerical results supporting our findings.
The proofs are gathered in Section \ref{hm:sec_proofs}.

%%%%%
%%%%%
\section{Preliminaries}
\label{hm:sec_pre}

%%%%%
\subsection{Basic notation}

The following basic notation will be used throughout this paper.
We denote by $|A|$ the determinant of a square matrix $A$, and by $\|A\|$ the Frobenius norm of a matrix $A$.
Write $A^{\otimes 2}=AA^\top$ for any matrix $A$, with $\top$ denoting transposition.
For a $K$th-order multilinear form $M=\{M^{(i_1\dots i_K)}:i_k=1,\dots, d_k;k=1,\dots,K\}\in\mbbr^{d_1}\otimes\dots\otimes\mbbr^{d_K}$ and $d_k$-dimensional vectors $u_k=\{u_k^{(j)}\}$, we let 
$M[u_1,\dots,u_K]:=\sum_{i_1=1}^{d_1}\dots\sum_{i_K=1}^{d_K} M^{(i_1,\dots,i_K)}u_1^{(i_1)}\dots u_K^{(i_K)}$; in particular, $A[B]:=\trace(AB^{\top})$ in case of $K=2$ for matrices $A$ and $B$ of the same sizes.
The symbol $\p_{a}^{k}$ stands for $k$-times partial differentiation with respect to variable $a$, and $I_{r}$ denotes the $r\times r$-identity matrix.
We write $C>0$ for a universal positive constant which may vary at each appearance, and 
$a_{n} \lesssim b_{n}$ for possibly random nonnegative sequences $(a_n)$ and $(b_n)$ if $a_{n}\le C b_{n}$ a.s. holds for every $n$ large enough.
The density of the Gaussian distribution $N_d(\mu,\Sig)$ will be denoted by $\phi_d(x;\mu,\Sig)$.

The basic setting is as follows.
We denote by $(\Omega,\mcf,\pr)$ the underlying probability space and by $\E$ the associated expectation operator.
%; sometimes, we will denote by $\pr_\theta$ the distribution of a solution process to \eqref{hm:SDE.single.model} with the parameter value $\theta$.
For notational convenience, instead of \eqref{hm:SDE.cand.model} we look at a single model
\begin{equation}
dX_t = c(X_{t-},\gam)dZ_t + a(X_t,\al)dt,
\label{hm:SDE.single.model}
\end{equation}
where $\gam=(\gam_k)\in\Theta_\gam\subset\mbbr^{p_\gam}$ and $\al=(\al_l)\in\Theta_\al\subset\mbbr^{p_\al}$, both parameter spaces being bounded convex domains. Let $p:=p_\al+p_\gam$.
%write $a^{(i)}$ and $c^{(i,j)}$ for the $i$th and $(i,j)$th entry of $a$ and $c$, respectively, and also 
Let $\D_{j}Y:=Y_{t_{j}}-Y_{t_{j-1}}$ for a process $Y$, and $f_{j-1}(\theta):=f(X_{t_{j-1}},\theta)$ for any measurable function on $f:\mbbr^{d}\times\Theta$.
The symbols $\cip$ and $\cil$ denote the convergence in probability and distribution, respectively.
In Sections 2, 3 and 4, we assume that coefficients $c$ and $a$ are correctly specified in the sense that there exist $\gamma_{0}\in\Theta_{\gamma}$ and $\alpha_{0}\in\Theta_{\alpha}$ such that $c(\cdot,\gamma_{0})=C(\cdot)$ and $a(\cdot,\alpha_{0})=A(\cdot)$.
%Let $\tz=(\al_0,\gam_0) \in \Theta:=\Theta_\al\times\Theta_\gam$ denote the true value of $\theta=(\al,\gam)$.

%%%%%
\subsection{Two-stage Gaussian quasi-likelihood estimation}
\label{hm:sec_stepwise.GQMLE}

Write $S(x,\gam)=c(x,\gam)^{\otimes 2}$ for the scale matrix, which will play the role of diffusion matrix in the diffusion context.
Let $\nu(dz)$ denote the {\lm} of $Z$, and then for $i_{1},\dots,i_{m}\in\{1,\dots,r\}$ with $m\ge 3$ we write $\nu(m)$ the tensor consisting of all the $m$th-mixed moments of $\nu$:
\begin{equation}
\nu(m)=\{\nu_{i_{1}\dots i_{m}}(m)\}_{i_{1},\dots,i_{m}} := \left\{\int z_{i_{1}}\dots z_{i_{m}} \nu(dz)\right\}_{i_{1},\dots,i_{m}}.
\nonumber
\end{equation}
% [Note that in particular, $\nu(2)=I_r$ as a matrix.]

Denote by $\lambda_{\min}\{S(x,\gamma)\}$ the minimum eigenvalue of $S(x,\gamma)$.
The symbol $\mcc^{k,l}_{\sharp}$ for nonnegative integers $k$ and $l$ denotes the function space consisting of all measurable $f:\,\mbbr^d\times\overline{\Theta} \to \mbbr$ such that:
\begin{itemize}
\item $f(\cdot,\theta)$ is globally Lipschitz uniformly in $\theta\in\overline{\Theta}$;
\item $f(x,\theta)$ is $k$-times (resp. $l$-times) continuously differentiable in $x$ (resp. in $\theta$), respectively, all the partial derivatives are continuous over $\overline{\Theta}$ for each $x$, and the estimates
\begin{equation}
\max_{i\le k} \max_{j\le l} \sup_{\theta\in\overline{\Theta}} |\p_\theta^j \p_x^i f(x,\theta)| \lesssim 1+ |x|^{C_{k,l}}
\nonumber
\end{equation}
holds for some constant $C_{k,l}\ge 0$.
\end{itemize}

We need some regularity conditions on the process $(X, Z)$ to ensure the asymptotic normality and the uniform tail-probability estimate.

The following conditions are standard in the literature, essentially borrowed from \cite{Mas13} and \cite{MasUeh17-2}.

\begin{ass}[Moments]
\label{ass1}
$\E[Z_{1}]=0$, $\E[Z_{1}^{\otimes2}]=I_{r}$, and $\E[|Z_{1}|^{q}]<\infty$ for all $q>0$. 
\end{ass}

\begin{ass}[Smoothness and non-degeneracy]
\label{ass2}
The components of $a$ and $c$ belong to the class $\mcc^{2,4}_{\sharp}$, and
%$a^{(i)},\, c^{(i,j)} \in \mcc^{2,4}_{\sharp}$ for each $i$ and $j$, and 
\begin{equation}
\sup_{\gamma\in\overline{\Theta}_{\gamma}}
\lam_{\min}\{S(x,\gamma)\}^{-1} \lesssim 1+|x|^{C_{0}}
%\inf_{\gamma\in\overline{\Theta}_{\gamma}}\lam_{\min}\{S(x,\gamma)\} \gtrsim \frac{1}{(1+\|x\|)^{C_{0}}}.
\nonumber
\end{equation}
for some constant $C_0\ge 0$.
%\begin{itemize}
%\item[(i)] For $x_{1},x_{2}\in\mbbr^{d}$,
%\begin{align*}
%\sup_{\alpha\in\overline{\Theta}_{\alpha}}\|a(x_{1},\alpha)-a(x_{2},\alpha)\|+\sup_{\gamma\in\overline{\Theta}_{\gamma}}\|c(x_{1},\gamma)-c(x_{2},\gamma)\|
%\lesssim \|x_{1}-x_{2}\|.
%\end{align*}
%%
%\item[(ii)] $a,c\in\mcc^{2,5}(\mbbr^{d}\times \Theta)$ and there exists a constant $C_{0}\ge 0$ such that for $x\in\mbbr^{d}$ and $i\in \{0,1,2,3\}$ and $j \in \{0,1,2\}$,
%\begin{align}
%& \sup_{\alpha\in\overline{\Theta}_{\alpha}}\|\p_{x}^{j}\p_{\alpha}^{i}a(x,\alpha)\| + \sup_{\gamma\in\overline{\Theta}_{\gamma}}\|\p_{x}^{j}\p_{\gamma}^{i}c(x,\gamma)\|  
%\lesssim (1+\|x\|)^{C_{0}}, \nn\\
%& \inf_{\gamma\in\overline{\Theta}_{\gamma}}\lam_{\min}\{S(x,\gamma)\} \gtrsim (1+\|x\|)^{-C_{0}}.
%\nonumber
%\end{align}
%\end{itemize}
\end{ass}

%The Gershgorin circle theorem says that
%\begin{equation}
%\inf_{\al}\lam_{\min}\{S(x,\al)\} \ge \inf_{\al} \min_{1\le i\le d}\bigg( S_{ii}(x,\al) - \sum_{j\ne i}|S_{ij}(x,\al)| \bigg),
%\nonumber
%\end{equation}
%so that an easy sufficient condition for the last inequality in (ii) is that this lower bound is bounded below by $C(1+|x|)^{-C_{0}}$.

%Let $(P_{t})$ denote the transition probability of $X_{t}$.
%Given a function $\rho:\mbbr^{d}\to\mbbr_{+}$ and a signed measure $m$ on $d$-dimensional Borel space, we define
%\begin{align*}
%\|m\|_{\rho}&=\sup\left\{\{|m(f)|:f \text{ is } \mbbr \text{-valued and measurable, such that } |f|\leq\rho\right\}.
%\end{align*}

\begin{ass}[Stability]
\label{ass3}
There exists a probability measure $\pi=\pi_{\tz}$ such that for every $q>0$ we can find positive constant $a$ for which
\begin{align*}
\sup_{t\in\mbbr_{+}}e^{at}
\sup_{f:\,|f|\le g}\left| \int f(y) P_{t}(x,dy) - \int f(y) \pi(dy)\right| \lesssim g(x), \qquad x\in\mbbr^d,
\end{align*}
where $g(x):=1+\|x\|^{q}$ and $P_t(x,dy):=\pr(X_t\in dy|X_0=x)$.
Further, for every $q>0$,
\begin{align}
\sup_{t\in\mbbr_{+}}\E[|X_{t}|^{q}]<\infty.
\label{hm:bdd.moments}
\end{align}
%\begin{itemize}
%\item[(i)] There exists a probability measure $\pi=\pi_{\tz}$ such that for every $q>0$ we can find positive constant $a$ for which
%\begin{align*}
%\sup_{t\in\mbbr_{+}}e^{at}
%\sup_{f:\,|f|\le g}\left| \int f(y) P_{t}(x,dy) - \int f(y)\pi_{0}(dy)\right| \lesssim g(x), \qquad x\in\mbbr^d,
%\end{align*}
%where $g(x):=1+\|x\|^{q}$.
%%
%\item[(ii)] For every $q>0$, 
%\begin{align*}
%\sup_{t\in\mbbr_{+}}\E[|X_{t}|^{q}]<\infty,
%\end{align*}
%\end{itemize}
%where $\E$ denotes the expectation operator with respect to $P_{0}$.
%%There exists a probability measure $\pi=\pi_{0}$ such that
%%\begin{align*}
%%\frac{1}{T}\int_{0}^{T}g(X_{t})dt\cip\int_{\mbbr^{d}}g(x)\pi(dx),\qquad T\to\infty,
%%\end{align*}
%%for any measurable function $g \in L^{1}(\pi)$.
%%In addition, $\sup_{t\in\mbbrp}\E(|X_{t}|^{q}) <\infty$ for all $q>0$ in case where the constant $C_{0}$ in Assumption \ref{ass2}(ii) is positive.
\end{ass}

It follows from Assumptions \ref{ass2} and \ref{ass3} that
\begin{equation}
\frac{1}{n}\sumj g(X_{t_{j-1}},\theta) \cip \int g(x,\theta)\pi(dx),\qquad n\to\infty,
\nn%\label{hm:disc.LLN}
\end{equation}
uniformly in $\theta$ for sufficiently smooth function $g(x,\theta)$ whose partial derivative with respect to $x$ are of at most polynomial growth in $x$ uniformly in $\theta$. This can be seen in the standard moment estimates and the tightness argument: see \cite[p.1598 and Section 4.1.1]{Mas13}.
Also to be noted is that the seemingly stringent moment condition \eqref{hm:bdd.moments} could be removed in compensation for the boundedness of the coefficients and the uniform non-degeneracy of $S$: see \cite[Theorem 2.9]{Mas13}.

\medskip

The Euler approximation for \eqref{hm:SDE.single.model} under $\pr_\theta$ is given by
\begin{align}
X_{t_{j}}\approx X_{t_{j-1}}+a_{j-1}(\alpha)h+c_{j-1}(\gamma)\D_{j}Z.
\label{Euler_app}
\end{align}
Taking the small-time Gaussian approximation
\begin{align}
\mathcal{L}(X_{t_{j}}|X_{t_{j-1}}=x) \approx N_{d}\left(x+a(x,\alpha)h, hS(x,\gamma)\right)
\label{local_gauss}
\end{align}
into account, we are led to the joint GQLF $\mbbh_n(\theta)=\mbbh_n(\X_{n},\theta)$:
\begin{align}
\mbbh_n(\theta)
&:= \sumj \log \phi_d\left(X_{t_j};\, X_{t_{j-1}} + a_{j-1}(\al)h,\, h S_{j-1}(\gam)\right) \nn\\
&=-\frac{1}{2}\sumj\left( \log\left|2\pi hS_{j-1}(\gamma)\right|+\frac{1}{h}S_{j-1}^{-1}(\gamma)\left[\left(\D_{j}X-ha_{j-1}(\alpha)\right)^{\otimes2}\right]\right).
\label{hm:def_joint.GQLF}
\end{align}
%\tcr{We regard \eqref{local_gauss} as our statistical model, which is obviously misspecified; still, it is possible to estimate the true coefficients when they are correctly specified (see Theorem \ref{hm:thm_GQMLE} below).}
In the present study, although the model of interest is described by \eqref{hm:SDE.single.model}, we will not consider the associated exact likelihood. Instead, we will regard \eqref{local_gauss} as our statistical model and deal with the explicit GQLF \eqref{hm:def_joint.GQLF} based on the theoretically incorrect \eqref{local_gauss} for inference purposes; in this sense, our statistical model is misspecified.
Still, it is possible to estimate the true coefficients when they are correctly specified (see Theorem \ref{hm:thm_GQMLE} below).

We can write $\mbbh_n(\theta)=\mbbh_{1,n}(\gam) + \mbbh_{2,n}(\theta)$ where
\begin{align}
\mbbh_{1,n}(\gam;\X_{n})= \; \mbbh_{1,n}(\gam) &:= \sumj \log \phi_d\left(X_{t_j};\, X_{t_{j-1}},\, h S_{j-1}(\gam)\right),
\nn\\
%\label{hm:def_stepwise.GQLF_1}\\
\mbbh_{2,n}(\theta;\X_{n})= \; \mbbh_{2,n}(\theta) &:= 
\sumj \left(
S_{j-1}^{-1}(\gam)\left[ \D_j X,\,a_{j-1}(\al)\right] -\frac{h}{2}S_{j-1}^{-1}(\gam)\left[a_{j-1}^{\otimes 2}(\al)\right]
\right).
\label{hm:def_stepwise.GQLF_2}
\end{align}
The joint GQLF $\mbbh_n(\theta)$ has two different ``resolutions'' (see \eqref{hm:2-resolutions} below), which comes from the fact that the last term $c_{j-1}(\gamma)\D_{j}Z$ in the right-hand side of \eqref{Euler_app} is stochastically dominant compared with the second one $a_{j-1}(\alpha)h$.
It was seen in \cite{MasUeh17-2} that under suitable conditions including the ergodicity of $X$ that both $n^{-1}\mbbh_{1,n}(\gam)$ and $T_n^{-1}\mbbh_{2,n}(\al,\gam)$ have non-trivial limits (of the ergodic theorem) for each $\theta$, in particular the former limits depending on $\gam$ only.
Building on this observation, the following two-stage estimation strategy is suggested: \textit{first}, we estimate $\gam$ by $\ges\in\argmax_{\gam}\mbbh_{1,n}(\gam)$; \textit{and then}, estimate $\al$ by $\aes\in\argmax_{\al}\mbbh_{2,n}(\al)$ where, with a slight abuse of notation,
\begin{equation}
\mbbh_{2,n}(\al):=\mbbh_{2,n}(\al,\ges).
\nonumber
\end{equation}
Note that maximizing $\al\mapsto\mbbh_{2,n}(\al,\gam)$ \textit{given a value of $\gam$} amounts to maximizing the discrete-time approximation of the log-likelihood function corresponding to the continuous-time observation (see \cite[Chapter 7]{LipShi01-1} for details), and also to maximizing
\begin{equation}
%\al\mapsto
\mbbh_{2,n}^\ast(\theta) := \sumj \log \phi_d\left(X_{t_j};\, X_{t_{j-1}} + a_{j-1}(\al)h,\, h S_{j-1}(\gam)\right).
\label{hm:def_stepwise.GQLF_2ast}
\end{equation}
Therefore, in either case, the second stage itself may be recognized as a GQLF.
We denote this two-stage \textit{Gaussian quasi-maximum likelihood estimator (GQMLE)} by $\tes=(\aes,\ges)$, which is essentially the same as in the one considered in \cite{MasUeh17-2}.

%Let
%\begin{align}
%\mbby_{1,n}(\gam) &:= -\frac12 \int\left\{\trace\left(S(x,\gam)^{-1}S(x,\gam_0) - I_r\right) - \log\frac{|S(x,\gam)|}{|S(x,\gam_0)|}\right\}\pi(dx), \nn\\
%\mbby_{2,n}(\gam) &:= -\frac12 \int S^{-1}(x,\gam_0) \left[\left( a(x,\al)-a(x,\al_0) \right)^{\otimes 2}\right] \pi(dx).
%\nonumber
%\end{align}

\begin{ass}[Identifiability]
\label{hm:ass_iden}
There exist positive constants $\chi_\gam$ and $\chi_\al$ such that
\begin{align}
-\frac12 \int\left\{\trace\left(S(x,\gam)^{-1}S(x,\gam_0) - I_d\right) + \log\frac{|S(x,\gam)|}{|S(x,\gam_0)|}\right\}\pi(dx)
&\le -\chi_\gam |\gam-\gam_0|^2, \nn\\
-\frac12 \int S^{-1}(x,\gam_0) \left[\left( a(x,\al)-a(x,\al_0) \right)^{\otimes 2}\right] \pi(dx) &\le -\chi_\al |\al-\al_0|^2,
\nonumber
\end{align}
for every $\gam$ and $\al$.
\end{ass}
The two integrals in the left-hand sides in Assumption \ref{hm:ass_iden} correspond to the Kullback-Leibler divergences associated with $\mbbh_{1,n}$ and $\mbbh_{2,n}$, respectively.

To state the result we need to introduce the matrix $V(\tz)=\Gam(\tz)^{-1}\Sig(\tz)\Gam(\tz)^{-1}$ where
\begin{align}
\Sig(\tz)=(\Sig^{(kl)}(\tz))_{k,l} &:=\left(
\begin{array}{cc}
\Gam_\al(\tz) & W_{\al,\gam}(\tz) \\
W_{\al,\gam}(\tz)^{\top} & W_\gam(\gam_0)
\end{array}
\right), \nn\\
\Gam(\tz)=(\Gam^{(kl)}(\tz))_{k,l} &:=\diag\{\Gam_\al(\tz),\Gam_\gam(\gam_0)\},
\nonumber
\end{align}
with, letting $\Psi =(\Psi^{(kl)})_{k,l}:=\p_\gam(S^{-1})= -S^{-1}(\p S)S^{-1}$ where $\Psi^{(kl)} \in \mbbr^{p_\gam}$,
%and $\Xi =(\Xi^{(k)})_{k}:=\p_\al a$ where $\Xi^{(k)} \in \mbbr^{p_\al}$,
\begin{align}
\Gam_\al^{(kl)}(\tz)&:= \int S^{-1}(x,\beta_{0})\left[\p_{\al_k}a(x,\al_{0}), \p_{\al_l}a(x,\al_{0})\right]\pi(dx),
\nn\\
\Gam_\gam^{(kl)}(\gam_0) &:= \frac12 \int\mathrm{trace}\left[\left\{(S^{-1}\p_{\gam_k}S)(S^{-1}\p_{\gam_l}S)\right\}(x,\gam_{0})\right]\pi(dx),
\nn\\
W_{\al,\gam}^{(qr)}(\tz) &:= -\frac12 \sum_{k,l,k',l',s,t,t'}\nu_{stt'}(3)\int
\left((\Psi^{(kl)})^{(q)}(\p_{\al_{r}}a^{(k')})c^{(ks)}c^{(lt)}c^{(l't')}(S^{-1})^{(k'l')}\right)(x,\tz) \pi(dx),
\nn\\
W_{\gam}^{(qr)}(\gam_0) &:= \frac14 \sum_{k,l,k',l',s,t,s',t'}\nu_{sts't'}(4)
\int\left( (\Psi^{(kl)})^{(q)}(\Psi^{(k'l')})^{(r)} c^{(ks)}c^{(lt)}c^{(k's')}c^{(l't')} \right)(x,\gam_0) \pi(dx).
\nonumber
\end{align}
%with, for each $v'_{1}, v'_{2}\in\mbbr^{p_{\al}}$ and $v''_{1}, v''_{2}\in\mbbr^{p_{\gam}}$ and with $c^{(\cdot k)}(x,\gam)\in\mbbr^{d}$ denoting the $k$th column of $c(x,\gam)$,
%\begin{align}
%\Gam_\al(\tz)[v''_{1}, v''_{2}] &:= \int S^{-1}(x,\beta_{0})\left[\p_{\al}a(x,\al_{0})[v'_{1}], \p_{\al}a(x,\al_{0})[v'_{2}]\right]\pi(dx),
%\nn\\
%\Gam_\gam(\gam_0)[v'_{1}, v'_{2}] &:= \frac12 \int\mathrm{trace}\left[\left\{(S^{-1}\p_{\gam}S)\otimes (S^{-1}\p_{\gam}S)\right\}(x,\gam_{0})
%[v''_{1},v''_{2}]\right]\pi(dx), \nn\\
%W_{\al,\gam}(\tz)[v'_1,v''_{1}] &:= -\frac12 \int\sum_{k',l',s'}\nu_{k'l's'}(3)
%S^{-1}(x,\gam_{0})\left[\p_{\al}a(x,\al_{0})[v_1'],c^{(\cdot s')}(x,\gam_{0})\right] 
%\nonumber \\
%&{}\qquad\cdot\left\{\p_{\gam}(S^{-1})(x,\gam_{0})\right\}
%[v''_{1},c^{(\cdot k')}(x,\gam_{0}),c^{(\cdot l')}(x,\gam_{0})]\pi(dx), \nonumber \\
%W_{\gam}(\tz)[v''_{1},v''_{2}] &:= \frac14 \int\sum_{s,t,s',t'}\nu_{sts't'}(4)
%\left\{\p_{\gam}S^{-1}(x,\gam_{0})[v''_{1},c^{(\cdot s)}(x,\gam_{0}),c^{(\cdot t)}(x,\gam_{0})]\right\}
%\nonumber \\
%&{}\qquad\cdot
%\left\{\p_{\gam}S^{-1}(x,\gam_{0})[v''_{2},c^{(\cdot s')}(x,\gam_{0}),c^{(\cdot t')}(x,\gam_{0})]\right\}\pi(dx);
%\nonumber
%\end{align}
%we note that $\p_{\gam}S^{-1}$ denotes $\p_{\gam}(S^{-1})=-S(\p S) S$.
%
Under Assumption \ref{hm:ass_iden}, the matrix $\Gam(\tz)$ is positive definite.
We introduce the following empirical counterparts: let $\hat{V}_n := \hat{\Gam}_n^{-1}\hat{\Sig}_n \hat{\Gam}_n^{-1}$, where
\begin{align}
\hat{\Sig}_n := \left(
\begin{array}{cc}
\hat{\Gam}_{\al,n} & \hat{W}_{\al,\gam,n} \\
\hat{W}_{\al,\gam,n}^{\top} & \hat{W}_{\gam,n}
\end{array}
\right), \qquad 
\hat{\Gam}_n :=\diag\{\hat{\Gam}_{\al,n},\,\hat{\Gam}_{\gam,n}\},
\nonumber
\end{align}
with, using the shorthand $\p_\theta\hat{f}_{j-1}$ for $(\p_\theta f_{j-1})(\tes)$ and so on, the entries of $\hat{\Gam}_{\al,n}$ and $\hat{\Gam}_{\gam,n}$ being given by
\begin{align}
\hat{\Gam}_{\al,n}^{(kl)} &:= \frac1n \sumj \hat{S}^{-1}_{j-1} \left[\p_{\al_k}\hat{a}_{j-1}, \p_{\al_l}\hat{a}_{j-1}\right],
\nn\\
\hat{\Gam}_{\gam,n}^{(kl)} &:= \frac{1}{2n} \sumj 
\trace\left[(\hat{S}^{-1}_{j-1}\p_{\gam_k}\hat{S}_{j-1})(\hat{S}^{-1}_{j-1}\p_{\gam_l}\hat{S}_{j-1})\right],
\nonumber
\end{align}
and also those of $\hat{W}_{\al,\gam,n}$ and $\hat{W}_{\gam,n}$ by, writing $\hat{\chi}_j = \D_j X - h\hat{a}_{j-1}$,
%\tcm{
\begin{align}
\hat{W}_{\al,\gam,n}^{(qr)} &:= 
\frac{1}{2T_n}\sumj \left\{\left(\hat{S}^{-1}_{j-1}(\p_{\gam_q}\hat{S}_{j-1})\hat{S}^{-1}_{j-1}\right)[\hat{\chi}_j^{\otimes 2}]\right\} \left\{(\hat{S}^{-1}_{j-1})[\hat{\chi}_j, \p_{\al_r}\hat{a}_{j-1}]\right\},
\nn\\
\hat{W}_{\gam,n}^{(qr)} &:= \frac{1}{4T_n} \sumj 
\left\{\left(\hat{S}^{-1}_{j-1}(\p_{\gam_q}\hat{S}_{j-1})\hat{S}^{-1}_{j-1}\right)[\hat{\chi}_j^{\otimes 2}]\right\}
\left\{\left(\hat{S}^{-1}_{j-1}(\p_{\gam_r}\hat{S}_{j-1})\hat{S}^{-1}_{j-1}\right)[\hat{\chi}_j^{\otimes 2}]\right\}.
\label{hm:W4-hat}
\end{align}
%} %% HM

\begin{thm}
\label{hm:thm_GQMLE}
Suppose that Assumptions \ref{ass1}, \ref{ass2}, \ref{ass3}, and \ref{hm:ass_iden} hold true. %, with $\Gam(\tz)$ being positive definite.
Then, for any continuous function $f:\,\mbbr^p\to\mbbr$ of at most polynomial growth, we have the convergence of moments
\begin{equation}
\E\left[f\left(\sqrt{T_n}(\tes-\tz)\right)\right] \to \int f(u) \phi(u;0,V(\tz))du;
\nn%\label{hm:mighty.conv}
\end{equation}
in particular, we have
\begin{equation}
\hat{u}_n:=\sqrt{T_n}(\tes-\tz) \cil N_p\left(0,V(\tz)\right),
\nonumber
\end{equation}
$\E(\hat{u}_n)\to 0$, $\E(\hat{u}_n^{\otimes 2})\to V(\tz)$, and also $\sup_n\E(|\hat{u}_n|^q)<\infty$ for every $q>0$.
Further, we have $\hat{\Sig}_n\cip\Sig(\tz)$ and $\hat{\Gam}_n\cip\Gam(\tz)$, followed by $\hat{V}_n \cip V(\tz)$, so that
\begin{equation}
\hat{V}_n^{-1/2}\sqrt{T_n}(\tes-\tz) \cil N_p(0,I_p),
\label{hm:GQMLE.asn}
\end{equation}
as soon as $\Sig(\tz)$ is positive definite (hence so is $V(\tz)$).
\end{thm}

\medskip

Comparing \cite[Theorem 2.7]{Mas13} and Theorem \ref{hm:thm_GQMLE} shows that the two GQMLEs have the same asymptotic distribution, so that there is no loss of asymptotic efficiency when using the two-stage procedure instead of the joint one.
%The previous study \cite{Mas13} was concerned with the moment convergence of the joint GQMLE $\tes'\in\argmax\mbbh_n$ through $Z$-estimation through the joint GQLF $\mbbh_n(\theta)$. We could also consider the $M$-estimation framework where, typically, regularity conditions become slightly less restrictive.
For convenience, we give a sketch of the proof in Section \ref{hm:sec_proof_stepwise.GQMLE}.
The positive definiteness of $\Sig_0$ seems not straightforward to be verified. Yet, if $\mcl(Z_1)$ is assumed to be symmetric from the beginning so that $\nu(3)=0$, then $\Sig(\tz)=\diag\{\Gam_\al(\tz),W_\gam(\gam_0)\}$ and the assumption reduce to the positive definiteness of $W_\gam(\gam_0)$.

\begin{rem}[Asymptotic normality of the joint GQMLE]\normalfont
\label{hm:rem_mixed.rates.M}
Although \cite{Mas13} used a $Z$-estimation framework, we may follow an $M$-estimation one; typically, regularity conditions in the latter case are less restrictive.
Recalling that $\mbbh_n(\theta)=\mbbh_{1,n}(\gam) + \mbbh_{2,n}(\theta)$, we have the mixed-rates expression
\begin{align}
\mbbh_n(\tz+T_n^{-1/2}u) - \mbbh_n(\tz) &= \frac1h\left\{ h\left(\mbbh_{1,n}(\gam_0+T_n^{-1/2}u_\gam) - \mbbh_{1,n}(\gam_0)\right)\right\} \nn\\
&{}\qquad + \left(\mbbh_{2,n}(\tz+T_n^{-1/2}u) - \mbbh_{2,n}(\tz)\right) \nn\\
&=: \frac1h \log \mbbz_{1,n}(u_\gam) + \log \mbbz_{2,n}(u),
\label{hm:2-resolutions}
\end{align}
where $u=(u_\al,u_\gam)\in\mbbr^{p_\al}\times\mbbr^{p_\gam}$.
Importantly, it can be shown that both of the random functions $u_\gam\mapsto \log \mbbz_{1,n}(u_\gam)$ and $u\mapsto \log \mbbz_{2,n}(u)$ are locally asymptotically quadratic with
\begin{align}
\log \mbbz_{1,n}(u_\gam) &\cil \D_{1}(\gam_0)[u_\gam] - \frac12 \Gam_\gam(\gam_0)[u_\gam^{\otimes 2}] =: f_0(u_\gam), \nn\\
\log \mbbz_{2,n}(u) &\cil 
\D_{21}(\tz)[u_\al] %+ \D_{22}(\tz)[u_\gam] 
- \frac12 \Gam_{\al}(\tz)[u_\al^{\otimes 2}] - \frac12 \Gam_{2,\gam}(\tz)[u_\gam^{\otimes 2}] =: g_0(u),
\nonumber
\end{align}
where $\D_{1}(\gam_0)$ and $\D_{21}(\tz)$ are weak limits of $T_n^{-1/2} h\p_\gam\mbbh_{1,n}(\gam_0)$ and $T_n^{-1/2} \p_\al\mbbh_{2,n}(\tz)$, respectively, where $\Gam_{2,\gam}(\tz)$ is the limit in probability of 
$-(1/2)n^{-1}\sumj \p_\gam^2(S^{-1})_{j-1}(\gam_0)[a_{j-1}(\al_0)^{\otimes 2}]$, and where $\argmax_{u_{\gam}} f_0(u_{\gam}) = \{\Gam_\gam(\gam_0)^{-1}\D_1(\gam_0)\} =:\{\hat{u}_{\gam,0}\}$ and $\argmax_{u_{\al}} g_0(u_{\al}, \hat{u}_{\gam,0}) = \{\Gam_\al(\al_0)^{-1}\D_{21}(\tz)\}=:\{\hat{u}_{\al,0}\}$ a.s.
Note that in the limit we have no ``cross term'' involving both $u_\al$ and $u_\gam$, entailing that $\hat{u}_0:=(\hat{u}_{\al,0},\hat{u}_{\gam,0})\sim N_p(0,V(\tz))$.
Now, by means of \cite[Theorem 1]{Rad08}, it is possible to deduce the asymptotic normality $\sqrt{T_n}(\tes'-\tz) \cil N_p\left(0,V(\tz)\right)$ of the joint GQMLE $\tes'\in\argmax_{\theta}\mbbh_n(\theta)$.
\end{rem}

\begin{rem}[Univariate case]\normalfont
\label{hm:rem_1dim}
If $d=r=1$, then the asymptotic covariance matrix $V(\tz)$ takes the following much simpler forms:
\begin{align}
& \Gam_\al(\tz)= \int \left(\frac{\p_\al a(x,\al_0)}{c(x,\gam_0)}\right)^{\otimes 2}\pi(dx),
\qquad 
\Gam_\gam(\gam_0)= \frac12 \int \left(\frac{\p_\gam S(x,\gam_0)}{S(x,\gam_0)}\right)^{\otimes 2}\pi(dx),
\nn\\
& W_{\al,\gam}(\tz)=\frac{\nu(3)}{2} \int \frac{(\p_\gam S) \otimes (\p_\al a)}{c^3}(x,\tz) \pi(dx),
\nn\\
& W_{\gam}(\gam_0)=\frac{\nu(4)}{4} \int \left(\frac{\p_\gam S(x,\gam_0)}{S(x,\gam_0)}\right)^{\otimes 2}\pi(dx)
= \frac{\nu(4)}{2}\Gam_\gam(\gam_0).
\nonumber
\end{align}
\end{rem}

\begin{rem}\normalfont
\label{hm:rem_diffusion+1}
We are primarily interested in cases of a driving {\lp} $Z$ having a non-null jump part.
The studentization (asymptotic standard normality) \eqref{hm:GQMLE.asn} is, however, valid as it is even for the case of diffusion where $Z=w$, an $r$-dimensional standard {\wp}. This implies that the factor $\hat{V}_n$ automatically distinguishes whether or not the driving noise is Gaussian or not.
To see this, we recall the well-known fact $D_n(\tes-\tz) \cil N_p(0,\Gam(\tz)^{-1})$ with the different (partly faster) rate of convergence $D_n:=\diag(\sqrt{T_n}I_{p_\al}, \sqrt{n}I_{p_{\gam}})$, see e.g. \cite{UchYos12}. By writing
\begin{align}
\hat{V}_n^{-1/2}\sqrt{T_n}(\tes-\tz) = \left( \sqrt{T_n}\hat{V}_n D_n^{-1}\right) D_n(\tes-\tz)
\nonumber
\end{align}
and then taking into account the difference of the orders of the conditional moments of $\D_j X$ given $\mcf_{t_{j-1}}$ (see e.g. \cite[Lemma 4.5]{Mas13}), it is possible to deduce that $\sqrt{T_n}\hat{V}_n D_n^{-1} \cip \Gam(\tz)^{1/2}$, hence \eqref{hm:GQMLE.asn}; the proof is standard and omitted.
Concerned with the rate of convergence, if $Z$ has a diffusion component and compound-Poisson jumps, 
%the asymptotic normality is established with 
the same rate as in the diffusion case can be achieved: $\sqrt{T_n}$ for the drift (and jump) component and $\sqrt{n}$ for the diffusion one (see \cite{ShiYos06} and \cite{OgiYos11}).
If $Z$ is a locally $\beta$-stable pure-jump L\'{e}vy process, the drift parameter has the rate of convergence $\sqrt{n}h^{1-1/\beta}$ (see \cite{CleGlo20}, \cite{Mas10ejs}, \cite{Mas19spa}, and \cite{masuda2023optimal}).
It should be noted that the rates of convergence affect the regularization term in the quasi-BIC type statistics.
See Remark \ref{hm:rem_other.QLF}.
\end{rem}

%%%%%
%%%%%
\section{Gaussian quasi-AIC}
\label{hm:sec_GQAIC}

Building on Theorem \ref{hm:thm_GQMLE}, we turn to AIC-type model selection.
Before proceeding, let us briefly describe the classical Akaike paradigm in a general setting:
for the unknown true distribution $g(x)\mu(x)$ of a sample $\X_n$, we are given a statistical model, say $\{f(\cdot;\theta):\,\theta\in\Theta\}$, and a divergence $\mcd(f;g)$ measuring deviation from $g$ to $f$.
We will follow the standard route by taking the Kullback-Leibler divergence:
we estimate $g$ by $\hat{f}_n(\cdot)=f(\cdot;\tes)$ for some estimator $\tes=\tes(\X_n)$, and look at the random quantity
$\mcd(\hat{f}_n;g)$, where
\begin{equation}
%\mcd(\hat{f}_n;g) = \left.\left(\int g \log\left(\frac{g}{f_n}\right) d\mu \right)\right|_{f_n=\hat{f}_n},
%\mcd(\hat{f}_n;g) = \int \log\left(\frac{g(z)}{\hat{f}_n(z;\tes)}\right) g(z)\mu(dz),
\mcd(f;g) = \int \log\left(\frac{g}{f}\right) g d\mu,
\nonumber
\end{equation}
which we want to minimize over given candidate models.
It amounts to minimizing the relative entropy 
$\mce(\hat{f}_n;g)$ where $\mce(f;g) := -\int (\log f) g d\mu$.
%$\mce(\hat{f}_n;g) := -\left.\left(\int g \log f_n d\mu \right)\right|_{f_n=\hat{f}_n}$.
As $g$ is unknown, we substitute the empirical counterpart 
%$\mce(\hat{f}_n;\del_{\X_n}):=-\sumj \log f(X_j;\tes)$ for $\mce(\hat{f}_n;g)$.
$\mce(\hat{f}_n;\del_{\X_n})$ for $\mce(\hat{f}_n;g)$.
Then, by removing the randomness by integrating out $\X_n$ with respect to $g$, it is desired to derive a computable corrector $\hat{\mfb}_n$ such that %for some $c_n>0$ \tcr{such that $\liminf$},
\begin{equation}
%c_n \times 
\E\left[ \mce(\hat{f}_n;g) - \left(\mce(\hat{f}_n;\del_{\X_n}) + \hat{\mfb}_n\right) \right] = o(1), \qquad n\to\infty.
\label{hm:AIC_intro1}
\end{equation}
The routine way is to first compute the ``leading'' term(s) of
\begin{equation}
\mfb_n := \E\left[ \mce(\hat{f}_n;g)-\mce(\hat{f}_n;\del_{\X_n})\right],
\label{hm:AIC_intro2}
\end{equation}
and then construct an asymptotically unbiased estimator $\hat{\mfb}_n=\hat{\mfb}_n(\X_n)$, namely $\E[\hat{\mfb}_n - \mfb_n]=o(1)$.
The stochastic order of $\hat{\mfb}_n$ should be strictly smaller than that of $\mce(\hat{f}_n;\del_{\X_n})$;
%and typically we can take $c_n\equiv 1$, as 
%in many cases including the i.i.d. models, we have $\mce(\hat{f}_n;\del_{\X_n}) = -\sumj \log f(X_j;\tes) = O_p(n)$ and $\hat{\mfb}_n$ can be taken to be the number of the unknown parameters in the model \cite{Aka73}.
For regular models, very often $\hat{\mfb}_n$ equals the number of the unknown parameters involved in the model, although it may not be the case, depending on each situation; see \cite{Aka73} and \cite{KonKit96} with the references therein for details.

The above strategy remains the same and makes sense, whatever the probabilistic structure of $\mcl(\X_n)$ is and when $\mcd(f;g)$ is replaced by other divergences.
Moreover, since we are considering a relative comparison among the candidate models, the above claim is satisfied as long as the suitable asymptotic properties of estimators are given even if the model $\{f(\cdot;\theta):\,\theta\in\Theta\}$ is misspecified. In other words, what is essential is that we can find a suitable corrector $\hat{\mfb}_n$ satisfying the property \eqref{hm:AIC_intro1}; again, note that our statistical models are all misspecified since the intractable true log-likelihood function was replaced by the fake Gaussian ones.
Still, the candidate coefficients $c_{m_1}(x,\gamma_{m_1})$ and $a_{m_2}(x,\alpha_{m_2})$ 
%are assumed to be correctly specified and they 
are estimable as specified in Theorem \ref{hm:thm_GQMLE}.

%%%%%
\subsection{Expansion of moments: joint case}
\label{hm:sec_joint.GQAIC}

Turning to our setup, we keep considering the SDE model \eqref{hm:SDE.single.model}.
We first observe what will occur if we look at the AIC associated not with the two-stage GQLF $(\mbbh_{1,n},\mbbh_{2,n})$, but with the \textit{joint} GQLF $\mbbh_n(\theta)=\mbbh_n(\theta;\X_n)=\mbbh_{1,n}(\gam)+\mbbh_{2,n}(\theta)$ of \eqref{hm:def_joint.GQLF}; most of the computations will be of direct use in the stepwise case of our primary interest as well (Section \ref{hm:sec_stepwise.GQAIC}).

In what follows, we will omit a large portion of the technical details, for they are essentially based on and quite analogous to the basic computations in \cite[Section 4]{Mas13}.

Denote by $\tilde{\X}_n$ the independent copy of $\X_n=(X_{t_j})_{j=0}^n$, and by $\tilde{E}$ the expectation operator with respect to $\mcl(\tilde{\X}_n)$. Let $\hat{u}_{\gam,n} := \sqrt{T_n}(\ges-\gam_0)$ and $\hat{u}_{\al,n} := \sqrt{T_n}(\aes-\al_0)$, so that $\hat{u}_n=(\hat{u}_{\al,n},\hat{u}_{\gam,n})$; by Theorem \ref{hm:thm_GQMLE}, both quantities are $L^q(\pr)$-bounded for every $q>0$. %in $\bigcap_{q>0}L^{q}(\pr)$.
As was mentioned in \eqref{hm:AIC_intro2}, we want to compute (approximate) the quantity
\begin{align}
\mfb_n &:= \E\left[ \mbbh_n(\tes(\X_n);\X_n) - \tilde{\E}\left[\mbbh_n(\tes(\X_n);\tilde{\X}_n)\right]\right] \nn\\
&= \E\tilde{\E}\left[ \left(\mbbh_n(\tes) - \mbbh_n(\tz)\right) - \left(\tilde{\mbbh}_n(\tes) - \tilde{\mbbh}_n(\tz)\right)\right].
%\nn\\
%&= \E\tilde{\E}\left[ (\mbbh -\tilde{\mbbh})_n(\tes) - (\mbbh -\tilde{\mbbh})_n(\tz)\right] \nn\\
\nn
\end{align}
Concerned with the first term $\mbbh_n(\tes) - \mbbh_n(\tz)$ inside the sign $E\tilde{E}$, we note that
\begin{equation}
\forall k\ge 1\,\forall l\ge 0,\quad \p_\al^k\p_\gam^l\mbbh_{n}(\theta)=\p_\al^k\p_\gam^l\mbbh_{2,n}(\theta).
\nonumber
\end{equation}
%$\p_\al^k\p_\gam^l\mbbh_{n}(\theta)=\p_\al^k\p_\gam^l\mbbh_{2,n}(\theta)$ for $k\ge 1$ and $l\ge 0$.
Let $\D_{1,n}(\gam_0):=h\,T_n^{-1/2}\p_\gam \mbbh_{1,n}(\gam_0)$ and $\D_{2,n}(\tz):=T_n^{-1/2}\p_\al \mbbh_{2,n}(\tz)$.
%Let further $\Gam_{1,n}(\gam_0):=-n^{-1} \p_\gam^2 \mbbh_{1,n}(\gam_0)$ and $\Gam_{2,n}(\tz):=- T_n^{-1} \p_\theta^2 \mbbh_{2,n}(\tz)$.
We have
\begin{align}
\Gam_{\gam,n}(\tz) &:= -\frac1n \p_\gam^2 \mbbh_{n}(\tz) \cip \Gam_\gam(\gam_0),
\label{hm:Gam1->}\\
\Gam_{\al,n}(\tz) &:= -\frac{1}{T_n} \p_\alpha^2 \mbbh_{2,n}(\tz) \cip \Gam_\al(\tz).
\nonumber
\end{align}
Then, using the third-order Taylor expansion, we can derive the following key asymptotically quadratic structure having two different resolutions for $\al$ and $\gam$: 
for some $\check{\theta}_n=\check{\theta}_n(s)=s\tz+(1-s)\tes$ with random $s\in[0,1]$,
\begin{align}
\mbbh_n(\tes) - \mbbh_n(\tz) 
&= 
\frac{1}{\sqrt{T_n}}\p_\theta\mbbh_n(\tz) [\hat{u}_n] - \frac12 \left(-\frac{1}{T_n}\p_\theta^2 \mbbh_n(\tz)\right) [\hat{u}_n^{\otimes 2}]
\nn\\
&{}\qquad + \frac{1}{\sqrt{T_n}} \left(\frac{1}{6T_n}\p_\theta^3 \mbbh_n(\check{\theta}_n)\right) [\hat{u}_n^{\otimes 3}]
\nn\\
&=: \frac1h \mcq_{1,n}(\gam_0)  %\sqrt{T_n}\, \mcr_n(\tz) 
+ \mcq_{2,n}(\tz),
\label{hm:joint.se_1}
\end{align}
where
\begin{align}
\mcq_{1,n}(\gam_0) 
&:= \D_{1,n}(\gam_0) [\hat{u}_{\gam,n}] - \frac12 \Gam_\gam(\gam_0) [\hat{u}_{\gam,n}^{\otimes 2}]
+ \frac{1}{\sqrt{T_n}} R_{1,n}(\tes;\tz),
\label{hm:Q1_se}\\
%\mcr_{n}(\tz) 
%&:= \frac{1}{T_n}\p_\gam\mbbh_{2,n}(\tz) [\hat{u}_{\gam,n}],
%\nn\\
\mcq_{2,n}(\tz) 
&:= \D_{2,n}(\tz) [\hat{u}_{\al,n}] - \frac12 \Gam_{\al}(\tz) [\hat{u}_{\al,n}^{\otimes 2}] + \frac{1}{\sqrt{T_n}} R_{2,n}(\tes;\tz).
\label{hm:Q2_se}
\end{align}
Here the ``remainder'' terms are given as follows:
\begin{align}
R_{1,n}(\tes;\tz) &= nh^2\left(\frac{1}{T_n}\p_\gam\mbbh_{2,n}(\tz)\right) [\hat{u}_{\gam,n}]
-\frac12\left(\sqrt{T_n} \left( \Gam_{n,\gam}(\gam_0) - \Gam_\gam(\gam_0) \right)\right) [\hat{u}_{\gam,n}^{\otimes 2}]
\nn\\
&{}\qquad + \left(\frac{h}{6T_n}\p_\gam^3 \mbbh_n(\check{\theta}_n)\right) [\hat{u}_{\gam,n}^{\otimes 3}],
\nn\\
R_{2,n}(\tes;\tz) &= 
\left(\frac{1}{\sqrt{T_n}}\p_\al \p_\gam \mbbh_{2,n}(\tz)\right)[\hat{u}_{\al,n},\hat{u}_{\gam,n}]
-\frac12\left(\sqrt{T_n} \left( \Gam_{n,\al}(\tz) - \Gam_\al(\tz) \right)\right) [\hat{u}_{\al,n}^{\otimes 2}]
\nn\\
&{}\qquad + \sum_{r\in\mbbzp^p;\,|r|=3,\, |r_\al|\ge 1} \frac{1}{6T_n}\p_\theta^r \mbbh_{2,n}(\check{\theta}_n) \hat{u}_n^r,
\nonumber
\end{align}
where the standard multi-index notation is used for the summation sign in the latter ($\mbbzp:=\{0,1,2,\dots\}$).
Following the proofs of the basic lemmas in \cite[Sections 4.1.2 and 4.1.3]{Mas13}, we can deduce that all of the random sequences $\{\D_{1,n}(\gam_0)\}_n$, $\{\D_{2,n}(\tz)\}_n$, $\{R_{1,n}(\tes;\tz)\}_n$, and $\{R_{2,n}(\tes;\tz)\}_n$ are $L^q(\pr)$-bounded for every $q>0$.
%we note that $\{T_n^{-1}\p_\gam\mbbh_{2,n}(\tz)\}_n$ is $L^q(\pr)$-bounded for any $q>0$.
Moreover, we have (see \cite[Lemma 4.6]{Mas13})
\begin{equation}
\left(\D_{2,n}(\tz),\,\D_{1,n}(\gam_0) \right) \cil N_p(0,\Sig(\tz)).
\nonumber
\end{equation}
Completely analogously to \eqref{hm:joint.se_1}, with replacing $\mbbh_n$ by $\tilde{\mbbh}_n$ and so on we obtain the expression
\begin{equation}
\tilde{\mbbh}_n(\tes) - \tilde{\mbbh}_n(\tz) = \frac1h \tilde{\mcq}_{1,n}(\gam_0) + \tilde{\mcq}_{2,n}(\tz)
\nonumber
\end{equation}
with similar structures to \eqref{hm:Q1_se} and \eqref{hm:Q2_se}. Therefore
\begin{equation}
\mfb_n = \frac1h E\tilde{E}\left[\mcq_{1,n}(\gam_0)-\tilde{\mcq}_{1,n}(\gam_0)\right] 
+ E\tilde{E}\left[\mcq_{2,n}(\tz) - \tilde{\mcq}_{2,n}(\tz)\right].
\label{hm:joint.se_2}
\end{equation}

Concerning $\mcq_{1,n}(\gam_0) - \tilde{\mcq}_{1,n}(\gam_0)$, we apply the standard Taylor-expansion argument to obtain
\begin{equation}
\hat{u}_{\gam,n} = \Gam_{\gam}(\gam_0)^{-1}\D_{1,n}(\gam_0) + \frac{1}{\sqrt{T_n}}\del_{\gam,n},
%\cil N_{p_\gam}\left( 0,\Gam_\gam(\gam_0)^{-1}W_\gam(\gam_0)\Gam_\gam(\gam_0)^{-1} \right),
\label{hm:u.gam_se}
\end{equation}
with $\{\del_{\gam,n}\}_n$ being $L^q(\pr)$-bounded for every $q>0$.
From the stochastic expansions \eqref{hm:Q1_se} and \eqref{hm:u.gam_se}, it follows that
\begin{align}
\E\tilde{\E}\left[ \mcq_{1,n}(\tz) \right] 
&= \trace \left\{ \Gam_\gam(\gam_0)^{-1}\E[\D_{1,n}(\gam_0)^{\otimes 2}] + o(1) \right\} - \frac12 \Gam_{\gam}(\gam_0) \left[\E[\hat{u}_{\gam,n}^{\otimes 2}]\right] + O(T_n^{-1/2}) \nn\\
&= \trace \left\{ \Gam_\gam(\gam_0)^{-1}W_\gam(\gam_0)\right\} 
- \frac12 \Gam_{\gam}(\gam_0) \left[ \Gam_\gam(\gam_0)^{-1} W_\gam(\gam_0) \Gam_\gam(\gam_0)^{-1}\right] + O(T_n^{-1/2}).
\label{hm:u.EtE1_se+1}
\end{align}
Likewise, we can show that $\E\tilde{\E}[ \hat{\mcq}_{1,n}(\tz)]$ admits the same expansion as in \eqref{hm:u.EtE1_se+1} except that the first term on the right-hand side is replaced by $o(h\sqrt{T_n})=o(1)$: indeed, using the independence between $\X_n$ and $\tilde{\X}_n$, Burkholder's inequality, $\E[\hat{u}_{\gam,n}]=o(1)$, and also the obvious notation with tilde, we have
\begin{equation}
\E\tilde{\E}[\tilde{\D}_{1,n}(\gam_0) [\hat{u}_{\gam,n}]]=\tilde{\E}[\tilde{\D}_{1,n}(\gam_0)] \big[\E[\hat{u}_{\gam,n}]\big] = O(h\sqrt{T_n})\cdot o(1) = o(nh^3)=o(1).
\nonumber
\end{equation}
By subtraction and recalling that $nh^2\to 0$, we conclude that
\begin{equation}
E\tilde{E}\left[\mcq_{1,n}(\gam_0)-\tilde{\mcq}_{1,n}(\gam_0)\right] = 
\trace \left\{ \Gam_\gam(\gam_0)^{-1}W_\gam(\gam_0)\right\} + O(T_n^{-1/2}).
\label{hm:u.EtE1_se}
\end{equation}

\begin{rem}\normalfont
\label{hm:rem_1dim.GQAIC}
In case of $d=r=1$ we have $W_\gam(\gam_0)=\nu(4) \Gam_\gam(\gam_0)/2$ (Remark \ref{hm:rem_1dim})
the first term on the right-hand side of \eqref{hm:u.EtE1_se} becomes $p_\gam \nu(4)/2$.
\end{rem}

We can handle $\mcq_{2,n}(\tz)-\tilde{\mcq}_{2,n}(\tz)$ similarly. As in \eqref{hm:u.gam_se} we can derive the stochastic expansion
\begin{equation}
\hat{u}_{\al,n} = \Gam_{\al}(\tz)^{-1}\D_{2,n}(\tz) + \frac{1}{\sqrt{T_n}}\del_{\al,n},
\label{hm:u.al_se}
\end{equation}
with $\{\del_{\al,n}\}_n$ being $L^q(\pr)$-bounded for every $q>0$.
Substituting \eqref{hm:u.al_se} in \eqref{hm:Q2_se} and proceeding as in the case of $\mcq_{1,n}(\gam_0)-\tilde{\mcq}_{1,n}(\gam_0)$, we obtain
\begin{equation}
E\tilde{E}\left[\mcq_{2,n}(\tz)-\tilde{\mcq}_{2,n}(\tz)\right] = \trace \left\{ \Gam_\al(\tz)^{-1}\Gam_\al(\tz)\right\} + O(T_n^{-1/2}) = p_\al + O(T_n^{-1/2}).
\label{hm:u.EtE2_se}
\end{equation}

We have derived the bias expressions for the $\gam$ and $\al$ parts separately.
Combining \eqref{hm:joint.se_2} with \eqref{hm:u.EtE1_se} and \eqref{hm:u.EtE2_se}, we obtain the following proposition.
%Combining \eqref{hm:joint.se_2} with \eqref{hm:u.EtE1_se} and \eqref{hm:u.EtE2_se}, we arrive at the expression
\begin{prop}
Suppose that Assumptions \ref{ass1}, \ref{ass2}, \ref{ass3}, and \ref{hm:ass_iden} hold.
Then,
\begin{align}
\mfb_n %&= \E\left[ \mbbh_n(\tes(\X_n);\X_n) - \tilde{\E}\left[\mbbh_n(\tes(\X_n);\tilde{\X}_n)\right]\right]  \nn\\
&= \frac{1}{h} \left(\trace \left\{ \Gam_\gam(\gam_0)^{-1}W_\gam(\gam_0)\right\} + O(T_n^{-1/2})\right)+p_\al + O(T_n^{-1/2}).
\label{hm:bias.joint}
\end{align}
\label{se:prop_bias.joint}
\end{prop}
Note that $h^{-1}T_n^{-1/2}\to\infty$.
This means that the above form of $\mfb_n$ is inconvenient, for the residual term in the $\gam$ part becomes stochastically larger than the leading term of the $\al$ part.
To make the expansion fully explicit in decreasing order, we thus have to further expand $\mbbh_{1,n}(\gam)$. 
Rather roughly, recalling \eqref{hm:joint.se_1} and the subsequent paragraphs, we may formally write 
%(recall $\log \mbbz_{1,n}(u_\gam)$ in Remark \ref{hm:rem_mixed.rates.M})
\begin{align}
\mbbh_{n}(\tes) - \mbbh_{n}(\tz) 
%&= \sum_{k\in\mbbzp} \frac{1}{k!}\left\{\left(\frac{1}{n}\p_\theta^k\mbbh_n(\tz)\right) [\hat{u}_n^{\otimes k}]\right\}
%n\, T_n^{-k/2}
&= \sum_{k\in\mbbn} \frac{1}{k!}\left(\frac1n \p_\theta^k\mbbh_n(\tz)\right) [\hat{u}_n^{\otimes k}]\, n\, T_n^{-k/2}
%\nn\\
%&=: \sum_{k\in\mbbn} \, \frac{1}{k!} \zeta_{1,n}^{(k)}(\gam_0) \, n\, T_n^{-k/2}
\nonumber
\end{align}
with %$\{\zeta_{1,n}^{(k)}(\gam_0)\}_n$ being $L^q(\pr)$-bounded for any $q>0$.
each $n^{-1}\p_\theta^k\mbbh_n(\tz)$ being asymptotically non-null (as in \eqref{hm:Gam1->}).
This implies that we need to pick up the terms up to the order $k_0$ which is the minimal integer such that $n\, T_n^{-k_0/2} \to 0$;
since we are assuming that $nh^2\to 0$, the number $k_0$ is necessarily greater than or equal to $5$.
%\cite{KonKit03}

In sum, because of the mixed-rates structure, the direct evaluation of $\mfb_n$ based on the \textit{joint} GQLF $\mbbh_n(\theta)$ necessitates the higher-order derivatives of $\mbbh_n$, resulting in rather complicated expressions.
In the next section, we are going to take a different route through a \textit{stepwise} manner to bypass this annoying point.

%\begin{rem}\normalfont
%\label{hm:rem_diffusion+2}
%It would be worth mentioning that we have no mixed-rates structure in the case of diffusions.
%...
%\end{rem}

%%%%%
\subsection{Stepwise bias corrections}
\label{hm:sec_stepwise.GQAIC}

Building on the observations in the previous section, we can expect that the stepwise AIC procedure will work, making a simple formula of the bias correction for the scale coefficient. 
Recall the stepwise GQMLE $\ges\in\argmax_{\gam}\mbbh_{1,n}(\gam)$ and $\hat{\al}_{n}\in\argmax_{\al}\mbbh_{2,n}(\al)$ introduced in Section \ref{hm:sec_stepwise.GQMLE}.
%Specifically, we proceed as follows.

First, we focus on the relative comparison of the scale coefficient, looking at the quasi-likelihood $\mbbh_{1,n}(\gam)$.
By inspecting the derivation of \eqref{hm:bias.joint} in Section \ref{hm:sec_joint.GQAIC}, we see that the bias
\begin{equation}
\mfb_{\gam,n} := h \E\left[ \mbbh_{1,n}(\ges(\X_n);\X_n) - \tilde{\E}\left[\mbbh_{1,n}(\ges(\X_n);\tilde{\X}_n)\right]\right]
\nonumber
\end{equation}
admits the expression
\begin{equation*}
\mfb_{\gam,n} =\trace \left\{ \Gam_\gam(\gam_0)^{-1}W_\gam(\gam_0)\right\} + O(T_n^{-1/2}).
%\label{hm:stepwise.bias-1}
\end{equation*}
Hence it is natural to define (dividing by $h$)
\begin{equation}
\gqaic_{1,n} := -2\,\mbbh_{1,n}(\ges) + \frac2h \trace \left( \hat{\Gam}_{\gam,n}^{-1}\hat{W}_{\gam,n}\right)
\label{hm:GQAIC_gam}
\end{equation}
as the first-stage GQAIC; recall that $\hat{\Sig}_n\cip\Sig(\tz)$ and $\hat{\Gam}_n\cip\Gam(\tz)$ (Theorem \ref{hm:thm_GQMLE}).
If in particular $d=r=1$, then we may define (Remark \ref{hm:rem_1dim.GQAIC})
\begin{equation}
\gqaic_{1,n} = -2\,\mbbh_{1,n}(\ges) + \frac{p_\gam}{h} \hat{\nu}_n(4),
\label{hm:GQAIC_gam-2}
\end{equation}
where $\hat{\nu}_n(4)$ is a suitable consistent estimator of $\nu(4)$; 
it can be conveniently estimated by
%since $\nu(4)$ is here related to the fourth-order moment of $Z_h$, we may make use of
\begin{equation}
\hat{\nu}_n(4) := \frac{1}{T_n}\sumj \left(\frac{\D_j X}{c_{j-1}(\ges)}\right)^4 \cip \nu(4).
\nonumber
\end{equation}
Note that this convergence does hold in $L^1(\pr)$.
The penalty term in \eqref{hm:GQAIC_gam}, hence in \eqref{hm:GQAIC_gam-2} as well, is stochastically divergent at the non-standard order $1/h$, which is in sharp contrast to the classical AIC and also to the CIC of \cite{Uch10}.
Still, it can be seen that the first term $-2\,\mbbh_{1,n}(\ges)$ is the leading one:
$-2 hn^{-1}\,\mbbh_{1,n}(\ges)$ has a non-trivial constant limit in probability, while $2 (nh)^{-1} \trace( \hat{\Gam}_{\gam,n}^{-1}\hat{W}_{\gam,n}) = O_p(T_n^{-1})=o_p(1)$.
%The situation is therefore the same as in the classical AIC, say $-2n^{-1}\ell_n(\tes) + 2n^{-1}\text{dim}(\theta)$ with the normalized maximized log-likelihood $-2n^{-1}\ell_n(\tes)$ having a non-null constant limit of law of large numbers.

\medskip

Having the estimate $\ges$ in hand, we proceed to the bias evaluation concerning the drift coefficient.
Again inspecting the derivation of \eqref{hm:bias.joint} in Section \ref{hm:sec_joint.GQAIC}, we see that
\begin{align}
\mfb_{\al,n} &:= \E\left[ \mbbh_{2,n}(\aes,\ges) - \tilde{\E}\left[\tilde{\mbbh}_{2,n}(\aes,\ges)\right]\right]
\nn\\
&= \E\left[ \mbbh_{2,n}(\aes,\ges) - \mbbh_{2,n}(\al_0,\ges)\right]
-\E\tilde{\E}\left[\tilde{\mbbh}_{2,n}(\aes,\ges) - \tilde{\mbbh}_{2,n}(\al_0,\ges)\right]
\nn\\
&{}\qquad + \E\tilde{\E}\left[ \mbbh_{2,n}(\al_0,\ges) - \tilde{\mbbh}_{2,n}(\al_0,\ges)\right]
\nn\\
&=: \mfb_{\al,n}^{(1)} - \mfb_{\al,n}^{(2)} + \mfb_{A,n}.
\end{align}
Using the same devices as in Section \ref{hm:sec_joint.GQAIC} together with the results in \cite[Sections 4.1.2 and 4.1.3]{Mas13}, we can deduce that
\begin{align}
\mfb_{\al,n}^{(1)} &= 
-\E\tilde{\E}\left[ -\frac{1}{\sqrt{T_n}}\p_\al \mbbh_{2,n}(\aes,\ges) [\hat{u}_{\al,n}]
+\frac12 \left( -\frac{1}{T_n}  \p_\al^2 \mbbh_{2,n}(\check{\al}_n,\ges)[\hat{u}_{\al,n}^{\otimes 2}]\right)\right]
\nn\\
&= o(1) + \frac12 \trace\left(\Gam_{\al}(\tz) \E[\hat{u}_{\al,n}^{\otimes 2}]\right)
= o(1) + \frac12 p_\al, 
\label{se:stepwise.bias-2-1}
\end{align}
and that for some $\check{\al}^{\prime}_{n}=\check{\al}^{\prime}_{n}(s^{\prime})=s^{\prime}\al_{0}+(1-s^{\prime})\hat{\al}_{n}$ with random $s^{\prime}\in[0,1]$,
\begin{align}
\mfb_{\al,n}^{(2)} &=
-\E\tilde{\E}\left[ \frac{1}{\sqrt{T_n}}\p_\al \tilde{\mbbh}_{2,n}(\aes,\ges) [\hat{u}_{\al,n}]
-\frac12 \left( -\frac{1}{T_n}  \p_\al^2 \tilde{\mbbh}_{2,n}(\check{\al}^{\prime}_n,\ges)[\hat{u}_{\al,n}^{\otimes 2}]\right)\right]
\nn\\
&= o(1) - \frac12 \trace\left(\Gam_{\al}(\tz) \E[\hat{u}_{\al,n}^{\otimes 2}]\right)
= o(1) - \frac12 p_\al;
\label{se:stepwise.bias-2-2}
\end{align}
in part, we used the facts that $\pr[\p_\al \tilde{\mbbh}_{2,n}(\aes,\ges) = 0] \to 1$ and that $\{T_n^{-1/2}\p_\al \tilde{\mbbh}_{2,n}(\aes,\ges) [\hat{u}_{\al,n}]\}_n$ is $L^q$ bounded for any $q>0$.
Thus we have obtained
\begin{equation*}
\mfb_{\al,n}=p_\al + \mfb_{A,n} + o(1).
%\label{hm:stepwise.bias-2}
\end{equation*}
Recall that we are assuming that scale and drift coefficients are correctly specified.
Since $\mfb_{A,n}=\E\tilde{\E}[ \mbbh_{2,n}(\al_0,\ges) - \tilde{\mbbh}_{2,n}(\al_0,\ges)]$ is independent of the drift estimator $\aes$ and, given a $\ges$, is common to all the candidates, we may and do ignore $\mfb_{A,n}$ in relative model comparison.
Therefore, we define the second-stage GQAIC as the usual form
\begin{equation}
\gqaic_{2,n} := -2\,\mbbh_{2,n}(\aes) + 2 p_\al.
\label{se:stepwise.gqaic_dri}
\end{equation}

Summarizing the above observations yields the following result.

\begin{thm}
\label{hm:thm_GQAIC.stepwise.se}
Suppose that Assumptions \ref{ass1}, \ref{ass2}, \ref{ass3}, and \ref{hm:ass_iden} hold true, %with $\Gam(\tz)$ being positive definite, 
and that
\begin{equation}
\E\left[\trace \left( \hat{\Gam}_{\gam,n}^{-1}\hat{W}_{\gam,n}\right)\right] \to 
\trace \left\{ \Gam_\gam(\gam_0)^{-1}W_\gam(\gam_0)\right\}.
\label{hm:gqaic_thm-1}
\end{equation}
Then, we have $\mfb_{n}=h^{-1}\mfb_{\gam,n}+\mfb_{\al,n}$, where %\eqref{hm:stepwise.bias-1} and \eqref{hm:stepwise.bias-2} hold true.
\begin{align}
\mfb_{\gam,n}&=\trace \left\{ \Gam_\gam(\gam_0)^{-1}W_\gam(\gam_0)\right\} + O(T_n^{-1/2}), \label{hm:stepwise.bias-1}\\
\mfb_{\al,n}&=p_\al + \mfb_{A,n} + o(1). \label{hm:stepwise.bias-2}
\end{align}
\end{thm}

The equations \eqref{hm:u.EtE1_se}, \eqref{hm:u.EtE2_se}, and \eqref{hm:stepwise.bias-1} are derived similarly, while \eqref{hm:stepwise.bias-2} is derived by dividing the bias into three parts and expanding them.
The point here is that, by considering $\mfb_{\gam,n}$ and $\mfb_{\al,n}$ separately, we can bypass the problem of the residual term in the $\gamma$ part being stochastically larger than the leading term in the $\alpha$ part; recall the expression \eqref{hm:bias.joint}.
%Furthermore, considering $\mfb_{\gam,n}$ and $\mfb_{\al,n}$ separately avoids the problem of the residual term in the $\gamma$ part being stochastically larger than the leading term in the $\alpha$ part; recall the expression \eqref{hm:bias.joint}.

\medskip

Here is a simple sufficient condition for \eqref{hm:gqaic_thm-1}.

\begin{prop}
\label{hm:prop+1}
Under the assumptions in Theorem \ref{hm:thm_GQAIC.stepwise.se}, \eqref{hm:gqaic_thm-1} is implied by
\begin{equation}
\exists\del>0,\,\exists N\in\mbbn, \quad 
\sup_{n\ge N} \E\left[\lam_{\min}^{-(1+\del)}(\hat{\Gam}_{\gam,n})\right] < \infty.
\label{hm:v0.6_add1}
\end{equation}
\end{prop}

Unfortunately, verification of \eqref{hm:v0.6_add1} may not be technically trivial.
Naively, if the off-diagonal elements of $\hat{\Gam}_{\gam,n}$ are small enough in magnitude, a simple sufficient condition for \eqref{hm:v0.6_add1} can be given through the Gerschgorin circle theorem, which ensures that $|\lam_i| \ge |a_{ii}| - \sum_{j\ne i}|a_{ij}|$ for any eigenvalue $\lam_i$ of a square matrix $A=(a_{ij})$:
writing $M(x,\gam)=[M^{(k)}(x,\gam)]_{k=1}^{p_\gam}$ with $M^{(k)}(x,\gam):=(S^{-1}\p_{\gam_k}S)(x,\gam)$, we have
\begin{align}
\lam_{\min}(\hat{\Gam}_{\gam,n}) 
&\ge \min_{1\le k\le d} \bigg( \hat{\Gam}_{\gam,n}^{(kk)} - \sum_{l\ne k}|\hat{\Gam}_{\gam,n}^{(kl)}| \bigg) \nn\\
&= \min_{1\le k\le d} 
\frac{1}{2n} \sumj \left(
\trace\left[ (\hat{M}_{j-1}^{(k)})^2
\right]
-\sum_{l\ne k}\left|\trace\left[ \hat{M}_{j-1}^{(k)} \hat{M}_{j-1}^{(l)} \right]\right|
\right).
\nonumber
\end{align}
Then, \eqref{hm:v0.6_add1} holds if
\begin{equation}
\inf_{\gam}\min_{1\le k\le d} \left(
\trace\left[M^{(k)}(x,\gam)^2\right]
-\sum_{l\ne k}\left|\trace\left[
M^{(k)}(x,\gam)\,M^{(l)}(x,\gam)
\right]\right|
\right)
\gtrsim (1+|x|^{C})^{-1}.
\nonumber
\end{equation}
Things become simpler if, for example, we assume the spectral representation $S(x,\gam)=\sum_k \lam_k(x,\gam)\Pi_k(x)$ for some positive functions $\lam_k(x,\gam)$ and some projection-valued ones $\Pi_k(x)$.
%%\begin{rem}[Practical consideration]\normalfont
%%The bias-correction term in \eqref{hm:GQAIC_gam} can be unstable if $\hat{\Gam}_{\gam,n}^{-1}$ is very small, even if $\lam_{\min}(\Gam_\gam(\gam_0)>0$.
%%In practice, one may use a conveniently modified version of $\hat{\Gam}_{\gam,n}$ such as
%%\begin{align}
%%\hat{\Gam}_{\gam,n}^{-1}(\ep) 
%%&:= \hat{\Gam}_{\gam,n}^{-1}I(|\hat{\Gam}_{\gam,n}| > \ep) + \Gam_\ast^{-1}I_p I(|\hat{\Gam}_{\gam,n}| \le \ep) \nn\\
%%&:= \hat{\Gam}_{\gam,n}^{-1} + (\Gam_\ast^{-1} - \hat{\Gam}_{\gam,n}^{-1}) I_p I(|\hat{\Gam}_{\gam,n}| \le \ep)
%%\nonumber
%%\end{align}
%%for any invertible $\Gam_\ast$ and a sufficiently small constant $\ep>0$.
%%\end{rem}
We will get a little bit further into the issue of bounding inverse moments \eqref{hm:v0.6_add1} in Section \ref{hm:sec_GQAIC.verification}.

For the convergence of moments \eqref{hm:gqaic_thm-1}, we could bypass the non-trivial bound \eqref{hm:v0.6_add1} by suitably truncating the minimum eigenvalue $\lam_{\min}(\hat{\Gam}_{\gam,n})$. Specifically, define
\begin{equation}
\hat{\Gam}_{\gam,n}^{-1}(b_n) := \hat{\Gam}_{\gam,n}^{-1}I\left( \lam_{\min}(\hat{\Gam}_{\gam,n}) \ge b_n \right)
\nonumber
\end{equation}
for some positive sequence $(b_n)$ such that $b_n\to 0$ and that 
\begin{equation}
\exists \kappa >0,\quad b_n \gtrsim T_n^{-(1-\kappa)/2}.
\label{hm:add.v2-1}
\end{equation}
Let $\lam_n := \lam_{\min}(\hat{\Gam}_{\gam,n})$ and $\lam_0 := \lam_{\min}(\Gam_{\gam}(\gam_0))>0$ for brevity.
Then, since $\lam_n \cip \lam_0$, we have
\begin{equation}
\hat{\Gam}_{\gam,n}^{-1}(b_n) \cip \Gam_{\gam}(\gam_0)^{-1}.
\nonumber
\end{equation}
Observe that $\lam_n = \inf_{|u|=1}u^\top \hat{\Gam}_{\gam,n}u \ge \lam_0 - |\hat{\Gam}_{\gam,n}-\Gam_{\gam}(\gam_0)|$.
Also, $\sup_n \E[|\sqrt{T_n}(\hat{\Gam}_{\gam,n}-\Gam_{\gam}(\gam_0))|^q]<\infty$ for any $q>0$ under the present assumptions; see \cite{Mas13} for details.
%$\hat{G}_n :=\sqrt{T_n}(\hat{\Gam}_{\gam,n}-\Gam_{\gam}(\gam_0))$
%
Building on the above observations, for $K>1$ and $q>0$ and for $n$ large enough,
\begin{align}
\E\left[\big|\hat{\Gam}_{\gam,n}^{-1}(b_n)\big|^{K}\right]
&= \E\left[\big|\hat{\Gam}_{\gam,n}^{-1}\big|^{K};\, \lam_n \ge b_n\right]
\lesssim \E\left[ \lam_n^{-K};\, \lam_n \ge b_n\right] \nn\\
&= \int_0^\infty \pr\left[ \lam_n^{-K} I(\lam_n\ge b_n) \ge x\right]dx \nn\\
&= \int_0^\infty \pr\left[ b_n \le \lam_n \le x^{-1/K}\right]dx \nn\\
&\le 1 + \int_1^{b_n^{-K}}\pr\big[ \lam_n \le x^{-1/K}\big]dx \nn\\
&\le 1 + \int_1^{b_n^{-K}}\pr\left[ \lam_0 \le T_n^{-1/2}|\sqrt{T_n}(\hat{\Gam}_{\gam,n}-\Gam_{\gam}(\gam_0))| + x^{-1/K}\right]dx \nn\\
&
\lesssim 1+ \int_1^{b_n^{-K}}\left(\pr\left[ \bigl|\sqrt{T_n}(\hat{\Gam}_{\gam,n}-\Gam_{\gam}(\gam_0))\bigr| \ge \sqrt{T_n}\,x^{-1/K}\right] + \pr\big[\lam_0 \le 2 x^{-1/K}\big]\right)dx
\nn\\
&\lesssim 1+ \int_1^{b_n^{-K}}\pr\left[ \bigl|\sqrt{T_n}(\hat{\Gam}_{\gam,n}-\Gam_{\gam}(\gam_0))\bigr| \ge \sqrt{T_n}\,x^{-1/K}\right]dx \nn\\
&\lesssim 1+  \left(\sup_n \E\left[\bigl|\sqrt{T_n}(\hat{\Gam}_{\gam,n}-\Gam_{\gam}(\gam_0))\bigr|^q\right]\right) T_n^{-q/2} \, \int_1^{b_n^{-K}} x^{q/K}ds \nn\\
&\lesssim 1+  T_n^{-q/2} b_n^{-q-K}.
\nonumber
\end{align}
%\hmrev{
%Here, we used $\sup_n \E[|\sqrt{T_n}(\hat{\Gam}_{\gam,n}-\Gam_{\gam}(\gam_0))|^q]<\infty$ in the last step.
%}
The rightmost side is bounded in $n$ if $b_n \gtrsim T_n^{-(1-K/(q+K))/2}$, which holds under \eqref{hm:add.v2-1}.
%%by letting $q$ large enough (given any $K>1$).
%
In sum, under the additional ad-hoc tuning \eqref{hm:add.v2-1}, we could remove the requirement \eqref{hm:gqaic_thm-1} by adopting
\begin{equation}
\gqaic_{1,n}^{\flat} := -2\,\mbbh_{1,n}(\ges) + \frac2h \trace \left( \hat{\Gam}_{\gam,n}^{-1}(b_n)\hat{W}_{\gam,n}\right)
\label{hm:add.v2-2}
\end{equation}
instead of \eqref{hm:GQAIC_gam}; a similar modification can be applied to \eqref{hm:GQAIC_gam_m} below as well.

\begin{rem}\normalfont
Inspecting the derivation it is trivial that the same bias corrections as in Theorem \ref{hm:thm_GQAIC.stepwise.se} hold even if we replace $\mbbh_{2,n}(\theta)$ of \eqref{hm:def_stepwise.GQLF_2} by $\mbbh_{2,n}^\ast(\theta)$ of \eqref{hm:def_stepwise.GQLF_2ast} all through bias evaluation at the second step.
\end{rem}
 
\begin{rem}[GQAIC for diffusion]\normalfont
\label{hm:rem_diffusion}
%Although our primary interest in this paper is an SDE driven by a non-Gaussian {\lp},
In relation to Remark \ref{hm:rem_diffusion+1}, it is worth mentioning the case of a diffusion process where $Z$ is an $r$-dimensional standard {\wp} $w$.
\begin{enumerate}
\item Then, the first and second GQAICs become $-2\mbbh_{1,n}(\ges)+2p_\gam$ and $-2\mbbh_{2,n}(\aes)+2p_\al$, respectively.
Moreover, in the case of the joint estimation through the GQLF $\mbbh_n(\theta)$ of \eqref{hm:def_joint.GQLF}, the GQAIC is given by $-2\mbbh_{n}(\tes)+2(p_\al+p_\gam)$, showing that the GQAIC takes the same form as in the CIC of \cite{Uch10}, the contrast information criterion.
We omit the technical details of these observations, for they can be derived in an analogous way to the L\'{e}vy SDE case, with an essential difference that, although in this case the rates of convergence are different for the diffusion and drift parameters, we can simultaneously normalize the associated random field by a single matrix to conclude the locally asymptotically quadratic structure (see \cite[Section 6]{Yos11} for details):
\begin{align}
u=(u_\al,u_\gam) &\mapsto \mbbh_n(\al_0+T_n^{-1/2}u_\al, \gam_0+n^{-1/2}u_\gam) - \mbbh_n(\tz) \nn\\
&= \D_n(\tz)[u] - \frac12 \Gam(\tz)[u,u] + o_p(1),\qquad u\in\mbbr^p.
\nonumber
\end{align}
This implies that the associated statistical random field does not have the mixed-rates structure of such as \eqref{hm:2-resolutions}.
The difference between the Gaussian and the non-Gaussian cases comes from the fact that the random sequence $(h^{-1/2}w_h)_{h>0}$ is $L^K(\pr)$-bounded for any $K>0$ whereas it is not the case for $(h^{-1/2}Z_h)_{h>0}$ with non-Gaussian {\lp} $Z$;
more specifically, it is known that in the one-dimensional case,
\begin{equation}
\lim_{h\to 0} \frac1h \E[|Z_h|^K] = \int |z|^K\nu(dz)
\nonumber
\end{equation}
for $K>2$ if $\E[Z_1]=0$ and $\E[|Z_1|^K]<\infty$ (see \cite[Lemma 3.1]{AsmRos01}).
This is also why the matrix rate of convergence $D_n$ mentioned in Remark \ref{hm:rem_diffusion+1} emerges in the diffusion case.
We refer to \cite[Section 4.1.1]{Mas13} for related remarks.

\item 
Write $o^\ast_p(1)$ for a random sequence $(\zeta_n)_n$ such that $\E(|\zeta_n|^q)\to 0$ for any $q>0$.
It can be shown that $\hat{\Gam}_{\gam,n}^{(kl)} \cip \Gam_\gam^{(kl)}(\gam_0)$ and that (recall \eqref{hm:W4-hat}), through repeated compensations,
\begin{align}
%\hat{\Gam}_{\gam,n}^{(kl)} &= \frac{1}{2n} \sumj 
%\trace\left[(\hat{S}^{-1}_{j-1}\p_{\gam_k}\hat{S}_{j-1})(\hat{S}^{-1}_{j-1}\p_{\gam_l}\hat{S}_{j-1})\right] \cip \Gam_\gam^{(kl)}(\gam_0), \nn\\
\frac2h \hat{W}_{\gam,n}^{(qr)} &= \frac{1}{2n} \sumj 
\bigg\{\bigg(\hat{S}^{-1}_{j-1}(\p_{\gam_q}\hat{S}_{j-1})\hat{S}^{-1}_{j-1}\bigg)\bigg[\bigg(\frac{\hat{\chi}_j}{\sqrt{h}}\bigg)^{\otimes 2}\bigg]\bigg\}
\nn\\
&{}\qquad \times
\bigg\{\bigg(\hat{S}^{-1}_{j-1}(\p_{\gam_r}\hat{S}_{j-1})\hat{S}^{-1}_{j-1}\bigg)\bigg[\bigg(\frac{\hat{\chi}_j}{\sqrt{h}}\bigg)^{\otimes 2}\bigg]\bigg\}
\nn\\
&= \frac{1}{2n} \sumj 
\bigg\{\bigg(S^{-1}_{j-1}(\p_{\gam_q}S_{j-1})\bigg)\bigg[\bigg(\frac{\D_j w}{\sqrt{h}}\bigg)^{\otimes 2}\bigg]\bigg\}
\nn\\
&{}\qquad \times
\bigg\{\bigg(S^{-1}_{j-1}(\p_{\gam_r}S_{j-1})\bigg)\bigg[\bigg(\frac{\D_j w}{\sqrt{h}}\bigg)^{\otimes 2}\bigg]\bigg\} + o^\ast_p(1)
\nn\\
&= \frac{1}{2n} \sumj 
\E^{j-1}\Bigg[
\bigg\{\bigg(S^{-1}_{j-1}(\p_{\gam_q}S_{j-1})\bigg)\bigg[\bigg(\frac{\D_j w}{\sqrt{h}}\bigg)^{\otimes 2}\bigg]\bigg\}
\nn\\
&{}\qquad \times
\bigg\{\bigg(S^{-1}_{j-1}(\p_{\gam_r}S_{j-1})\bigg)\bigg[\bigg(\frac{\D_j w}{\sqrt{h}}\bigg)^{\otimes 2}\bigg]\bigg\}
\Bigg] + o^\ast_p(1)
\nn\\
&= \frac{1}{2n} \sumj 
\bigg\{ 2\trace\left(S^{-1}_{j-1}(\p_{\gam_q}S_{j-1}) S^{-1}_{j-1}(\p_{\gam_r}S_{j-1})\right) \nn\\
&{}\qquad 
+\trace\Big(\left(S^{-1}_{j-1}(\p_{\gam_q}S_{j-1})\right)\otimes \left(S^{-1}_{j-1}(\p_{\gam_r}S_{j-1})\right)\Big)
\bigg\} + o^\ast_p(1)
\nn\\
&= 2\Gam_\gam^{(qr)}(\gam_0) + 
\frac{1}{2n} \sumj 
\trace\left(\left(\hat{S}^{-1}_{j-1}(\p_{\gam_q}\hat{S}_{j-1})\right)\otimes \left(\hat{S}^{-1}_{j-1}(\p_{\gam_r}\hat{S}_{j-1})\right)\right)
+ o^\ast_p(1),
\label{hm:rem_diffusion-1}
\end{align}
where we used the identity \cite[Theorem 4.2(i)]{MagNeu79} for the fourth equality.
Write $\hat{A}^{(qr)}_{\gam,n}$ for the second term in \eqref{hm:rem_diffusion-1}, and let $\hat{A}_{\gam,n}:=(\hat{A}^{(qr)}_{\gam,n})_{q,r}$; obviously, we have $\hat{A}_{\gam,n}=O^\ast_p(1)$. The identity \eqref{hm:rem_diffusion-1} suggests us use the following modified version
\begin{equation}
\gqaic_{1,n} = -2\,\mbbh_{1,n}(\ges) + \trace \left\{ \hat{\Gam}_{\gam,n}^{-1} \left(\frac2h\hat{W}_{\gam,n}
-\hat{A}_{\gam,n}\right)\right\}
\label{hm:GQAIC_gam_m}
\end{equation}
instead of \eqref{hm:GQAIC_gam} as an alternative that can be used for both diffusion and L\'{e}vy driven SDE in common, in exchange for a slight additional computational cost.
In particular for $d=1$, instead of \eqref{hm:GQAIC_gam-2} we could use
\begin{equation}
\gqaic_{1,n} = -2\,\mbbh_{1,n}(\ges) + p_\gam\left(\frac{1}{h} \hat{\nu}_n(4) - \hat{\nu}_n(2)^2\right),
\nn%\label{hm:GQAIC_gam-2_m}
\end{equation}
since $h^{-1}\hat{\nu}_n(4) \cip 3$ and $\hat{\nu}_n(2) := T_n^{-1}\sumj c_{j-1}(\ges)^{-2}(\D_j X)^2 \cip 1$, with the latter holding for both diffusion and L\'{e}vy driven SDE while the former only for diffusion; or, more simply we could use $\gqaic_{1,n} = -2\,\mbbh_{1,n}(\ges) + p_\gam\left(h^{-1}\hat{\nu}_n(4) - 1\right)$ in common.
\end{enumerate}
\end{rem}

%%%%%
\subsection{Inverse-moment bound}
\label{hm:sec_GQAIC.verification}

In this section, we revisit \eqref{hm:v0.6_add1}. 
For notational simplicity, we write
\begin{equation}
\zeta(x,\gam) = (\p_\gam \log|S|)(x,\gam)
\nonumber
\end{equation}
and $\mbbs:=\{u\in\mbbr^{p_\gam }:\,|u|=1\}$.
We will prove the following criterion.

\begin{lem}
\label{hm:lem_imb}
Suppose that the assumptions given in Section \ref{hm:sec_stepwise.GQMLE} hold.
Moreover, suppose that there exist positive constants $\rho$ and $C'$ such that for each $\ep\in(0,1]$,
\begin{equation}
\sup_n \sup_{2\le j\le n}\sup_{u\in\mbbs}\sup_{\gam}
\left.\pr\left[ \big|\zeta(X_{t_{j-1}},\gam)[u]\big| < \ep \right| \mcf_{t_{j-2}}\right]
\le C'\ep^\rho\quad \text{a.s.}
\label{hm:lem_imb_A*}
\end{equation}
%%
%%Moreover, suppose that there exist positive constants $\Delta$, $\rho$ and $C'$ such that for each $\ep\in(0,1]$,
%%\begin{equation}
%%\sup_{t\ge 0}\sup_{u\in\mbbs}\sup_{\gam}
%%\left.\pr\left[ \big|\zeta(X_{t+\Delta},\gam)[u]\big| < \ep \right| \mcf_{t}\right]
%%\le C'\ep^\rho \quad \text{a.s.}
%%\label{hm:lem_imb_A}
%%\end{equation}
%%
Then \eqref{hm:v0.6_add1} holds.% for any $N$ large enough.
\end{lem}

Note that the bound \eqref{hm:lem_imb_A*} is implied by
\begin{equation}
\sup_n \sup_{2\le j\le n}\sup_{\gam}
\left.\pr\left[ 
\lam_{\min}\left(\zeta(X_{t_{j-1}},\gam)^{\otimes 2}\right)<\ep^2 \right| \mcf_{t_{j-2}}\right]
\le C'\ep^\rho\quad \text{a.s.}
\nn%\label{hm:lem_imb_A*+}
\end{equation}

The proof of Lemma \ref{hm:lem_imb} utilizes the technique which dates back to \cite{BhaPap91}, later improved and generalized by \cite{FinWei02} and \cite{ChaIng11}:
%\cite{SriBos88}.
For the sake of reference, we will give the proof in an almost self-contained form.
Write $q=1+\del$ % and $\mbbs:=\{u\in\mbbr^{p_\gam }:\,|u|=1\}$ 
in what follows. 
We recall the following lemma for later reference.

\begin{lem}[Lemma A.3 in \cite{FinWei02}]
\label{hm:lem_FinWei.A3}
For $\del\in(0,1)$, there exists a finite subset $\mbbs(\del)\subset\mbbs$ such that:
\begin{enumerate}
%\item $\limsup_{\del\downarrow 0}\del^{p_\gam -1} |\mbbs(\del)|<\infty$;
\item $\mbbs(\del)$ has at most $\lfloor C_{p_\gam} \del^{-(p_\gam -1)}\rfloor$ elements for some constant $C_{p_\gam}$ only depending on $p_\gam$;
\item For each $u\in\mbbs$ there exists an element $v\in\mbbs(\del)$ for which $|u-v|<\del$.
\end{enumerate}
\end{lem}

Here is a direct corollary to Lemma \ref{hm:lem_imb}.

\begin{cor}
\label{hm:lem_imb-cor1}
Suppose that the assumptions given in Section \ref{hm:sec_stepwise.GQMLE} hold.
Further, suppose that there exists a nonnegative measurable function $\underline{\lam}(x)$ for which
\begin{equation}
\inf_{\gam}\lam_{\min}\left(
\zeta(x,\gam)^{\otimes 2}\right) \ge \underline{\lam}(x),
\nn
\end{equation}
and for every $\Delta>0$ small enough,
\begin{equation}
\sup_{t\ge 0}\pr\bigg( \underline{\lam}(X_{t+\Delta}) \le \ep \,\bigg|\,\mcf_t  \bigg) \lesssim \ep^\rho\quad\text{a.s.}
\nonumber
\end{equation}
for $\ep\in(0,1]$.
Then \eqref{hm:v0.6_add1} holds for any $N$ large enough.
\end{cor}

We remark that the ``tuning'' parameter $k$ in the proof of Lemma \ref{hm:lem_imb} plays a role to relieve possible high concentration probability of $\lam_{\min}(\hat{\Gam}_{\gam,n})$ around the origin; obviously, it is redundant under the present assumptions if the stronger non-degeneracy condition of $\zeta(x,\gam)$ holds:

\begin{cor}
\label{hm:lem_imb-cor2}
Suppose that the assumptions given in Section \ref{hm:sec_stepwise.GQMLE} hold.
Further, suppose that
\begin{equation}
%\exists C\ge 0~
\inf_{\gam}\lam_{\min}\left(
\zeta(x,\gam)^{\otimes 2}\right) \gtrsim (1+|x|)^{-C}.
\label{hm:imb-6}
\end{equation}
Then \eqref{hm:v0.6_add1} holds for any $N$ large enough.
\end{cor}

%%%%%
%%%%%
\section{Gaussian quasi-BIC}
\label{hm:sec_GQBIC}

In this section, we consider a two-stage Schwarz's type Bayesian information criterion, termed Gaussian quasi-BIC ($\gqbic$), through the GQLF.
We keep using the notation introduced in Section \ref{hm:sec_pre}.
Suppose that Assumptions \ref{ass1}, \ref{ass2}, \ref{ass3}, and \ref{hm:ass_iden} hold true.
%Suppose that Assumptions \ref{ass1} to \ref{hm:ass_iden} hold true.
%Recall that we have two GQLFs $\mbbh_{1,n}(\gam)$ and $\mbbh_{2,n}(\al)=\mbbh_{2,n}(\al,\ges)$.
In addition, we consider the prior densities $\pi_1(\gam)$ and $\pi_2(\al)$ for $\al$ and $\gam$, respectively.
We assume that both $\pi_1$ and $\pi_2$ are continuous and bounded in $\overline{\Theta}_\gam$ and $\overline{\Theta}_\al$ respectively, and moreover that $\pi_1(\gam_0)>0$ and $\pi_2(\al_0)>0$.
Moreover, for a technical reason, throughout this section we assume that there exists a constant $c_1\in(0,1)$ for which
\begin{equation}
T_n \gtrsim n^{c_1}.
\label{hm:h-rate-add}
\end{equation}
This is a real restriction in addition to $nh^2\to 0$ as is seen by the example $h=n^{-1}\log n$.

%%%%%
\subsection{Scale}
\label{hm:sec_GQBIC.scale}

We introduce the stochastic expansion of the \textit{free energy at the inverse temperature $\mfb>0$}, which is defined using the negative normalized logarithmic partition function (we refer to \cite{Wat13} for relevant backgrounds):
\begin{equation}
\mfF_{1,n}(\mfb) := -\frac{1}{n \mfb}\log\left(\int_{\Theta_\gam}\exp\{\mfb\, \mbbh_{1,n}(\gam)\} \pi_1(\gam)d\gam\right).
\nonumber
\end{equation}
Here, the terminology ``normalized'' means that $\mfF_{1,n}(\mfb)$ has non-trivial limit (in probability) for each $\mfb>0$.
The normalized marginal quasi-log likelihood corresponds to $\mfF_{1,n}(1)$, and the classical BIC methodology is based on a stochastic expansion of $\mfF_{1,n}(1)$. See \cite{EguMas18a} and the references therein.

We will prove the following expansions in Section \ref{hm:sec_thm_gqbic1_proof}.

\begin{thm}
\label{hm:thm_gqbic1}
%\label{hm:thm_qbic_sqrt(n)}
%\label{hm:thm_qbic_it.h&sqrt(T)}
We have the following stochastic expansions:
\begin{align}
\mfF_{1,n}(1) &= - \frac1n \mbbh_{1,n}(\ges) + \frac{p_\gam}{2 n} \log n + O_p\left(\frac1n\right), \label{hm:gqbic_expa1_EU} \\
\mfF_{1,n}(h) &= - \frac1n \mbbh_{1,n}(\ges) + \frac{p_\gam}{2 T_n} \log T_n + O_p\left(\frac{1}{T_n}\right). \label{hm:gqbic_expa1}
\end{align}
\end{thm}

The first one \eqref{hm:gqbic_expa1_EU} was previously given in \cite{EguUeh21}, based on which the authors introduced the ${\gqbic}$ for the scale by
\begin{equation}
\gqbic_{1,n}^{\sharp} := - 2\mbbh_{1,n}(\ges) + p_\gam \log n.
\nonumber
\end{equation}
Theorem \ref{se:thm_sele.prob2} below revises the incorrect part of \cite[Theorem 3.2]{EguUeh21}, showing that $\gqbic_{1,n}^{\sharp}$ does \textit{not} bring about the model-selection consistency.

Instead, building on \eqref{hm:gqbic_expa1}, we propose to use
\begin{equation}
\gqbic_{1,n} := -2\,\mbbh_{1,n}(\ges) + \frac{p_\gam}{h} \log T_n.
\label{hm:GQBIC_gam}
\end{equation}
In Theorem \ref{se:thm_sele.prob3} below, it will show that this form has the model-selection consistency.
This implies that in the present L\'{e}vy driven SDE setting, we need to ``heat up'' the quasi-likelihood $\mbbh_{1,n}$ by multiplying $h^{-1}$.

%%%%%
\subsection{Drift}

Different from the previous scale case, the second stage QBIC corresponding to $\mbbh_{2,n}(\al)$ is standard: we do not need to heat up the second-stage quasi-likelihood $\mbbh_{1,n}(\al)$.
We can directly look at the normalized marginal quasi-log likelihood
\begin{equation}
\mfF_{2,n}=\mfF_{2,n}(1) := -\frac{1}{T_n}\log\left(\int_{\Theta_\al}\exp\{\mbbh_{2,n}(\al)\} \pi_2(\al)d\al\right).
\nonumber
\end{equation}

We have the following stochastic expansion:
\begin{thm}
\label{hm:thm_gqbic2}
\begin{equation}
\mfF_{2,n}(1) = - \frac{1}{T_n}\mbbh_{2,n}(\aes) + \frac{p_\al}{2 T_n} \log T_n + O_p\left(\frac{1}{T_n}\right).
\label{hm:gqbic_expa2}
\end{equation}
\end{thm}

Theorem \ref{hm:thm_gqbic2} can be proved similarly to the proof of \eqref{hm:gqbic_expa1} in Theorem \ref{hm:thm_gqbic1}. The proof of the stochastic expansion \eqref{hm:gqbic_expa2} is much simpler, and we omit the proof.

Ignoring the vanishing term $O_p(T_n^{-1})$ of $\mfF_{2,n}$ (just as in \cite{EguMas18a}), we introduce the $\gqbic$ for the drift in the same form as in \cite{EguUeh21}:
\begin{equation}
\gqbic_{2,n} = -2\,\mbbh_{2,n}(\aes) + p_\al \log T_n.
\label{hm:GQBIC_al}
\end{equation}
In the next section, we will formulate a two-step selection procedure for both $\gqaic$ and $\gqbic$.
%The overall (scale and drift) selection consistency can be proved as in the case of diffusion \cite{EguMas18a}.

\begin{rem}\normalfont
\label{hm:rem_other.QLF}
%In relation to Remark \ref{hm:rem_diffusion+1}, 
The proposed two-stage methodology itself is simple enough, and we believe that, in principle, it can be applied to other types of quasi-likelihoods such as the non-Gaussian stable one (\cite{CleGlo20}, \cite{JasKamMas19}, \cite{Mas19spa}, and \cite{masuda2023optimal} for details). 
As long as considering the ergodic case, the derivation of the AIC-type statistics remains valid since what is essential therein is the convergence of moments of the asymptotically normally distributed estimator.
The case of BIC-type statistics is easier to handle, for it is only based on the rate of convergence of the estimator;
of more interest is that different from the AIC type, it is not essential for the quasi-BIC statistics that the model is ergodic (see \cite{EguMas18a} and also \cite[Appendix]{EguMas19}).
Whatever the case, careful consideration and calculation are needed in terms of models.
We would like to leave these issues to future tasks.
\end{rem}

%%%%%
\section{Model comparison and asymptotic probability of relative model selection}
\label{se:sec_modsele}
In this section, we consider relative (pairwise) model selection probabilities of the GQAIC and GQBIC.
Below, we assume that $0<\sharp\mathfrak{M}_{1}<M_{1}$ and $0<\sharp\mathfrak{M}_{2}<M_{2}$, where 
$\sharp\mathfrak{M}_{1}$ and $\sharp\mathfrak{M}_{2}$ denote the numbers of elements of $\mathfrak{M}_{1}$ and $\mathfrak{M}_{2}$, respectively, with
\begin{align}
\mathfrak{M}_{1}&:=\{m_{1}\in\left\{1,\ldots,M_{1}\}: \text{there exists a } \gam_{m_{1},0}\in\Theta_{\gam_{m_{1}}} \text{ such that } c_{m_{1}}(\cdot,\gam_{m_{1},0})=C(\cdot)\right\}, \nn\\ %\label{se:def.scale.spec.set}\\
\mathfrak{M}_{2}&:=\{m_{2}\in\left\{1,\ldots,M_{2}\}: \text{there exists a } \al_{m_{2},0}\in\Theta_{\al_{m_{2}}} \text{ such that } a_{m_{2}}(\cdot,\al_{m_{2},0})=A(\cdot)\right\}. \nn %\label{se:def.drift.spec.set}
\end{align}
This means that the candidate coefficients $c_{1},\ldots,c_{M_{1}}$ and $a_{1},\ldots,a_{M_{2}}$ contain both correctly specified coefficients and misspecified coefficients.
In cases of L\'{e}vy driven SDEs where either or both of the drift and scale coefficients are misspecified, the asymptotic properties of estimators are shown in \cite{Ueh19} under suitable conditions.
Also, we formally use the $\gqaic_{1,n}$ and $\gqaic_{2,n}$ even for the possibly misspecified coefficients, although the assumptions of Theorem \ref{hm:thm_GQAIC.stepwise.se} may not hold.
Using the GQAIC, the stepwise model comparison is performed as follows.
\begin{itemize}
\item[(i)] We compute $\gqaic_{1,n}$ for each candidate scale coefficient, say $\gqaic_{1,n}^{(1)},\ldots,\gqaic_{1,n}^{(M_{1})}$, and select the best scale coefficient $c_{\hat{m}_{1,n}}$ having the minimum $\gqaic_{1,n}$-value:
\begin{align*}
\{\hat{m}_{1,n}\}=\argmin_{1\leq m_{1}\leq M_{1}}\gqaic_{1,n}^{(m_{1})}.
\end{align*}

\item[(ii)] Under the result of (i), we choose the best drift coefficient with index $\hat{m}_{2,n}$ such that
\begin{align*}
\{\hat{m}_{2,n}\}=\argmin_{1\leq m_{2}\leq M_{2}}\gqaic_{2,n}^{(m_{2}|\hat{m}_{1,n})},
\end{align*}
where $\gqaic_{2,n}^{(m_{2}|m_{1,n})}$ corresponds to \eqref{se:stepwise.gqaic_dri} with $c_{m_{1,n}}$ and $\hat{\gam}_{m_{1,n},n}$.
\end{itemize}
The total number of comparisons
in this procedure is $M_{1}+M_{2}$, and we can obtain the model $\mcm_{\hat{m}_{1,n},\hat{m}_{2,n}}$ as the final best model among the candidates.
When we use GQBIC for model comparison, the best model is selected by a similar procedure.

Let the functions $\mbbh_{1,n}^{(m_{1})}$ and $\mbbh_{2,n}^{(m_{2}|m_{1})}$ denote $\mbbh_{1,n}$ and $\mbbh_{2,n}$ in each candidate model $\mcm_{m_{1},m_{2}}$, respectively.
Then, we have
\begin{align*}
\frac{1}{n}\mbbh_{1,n}^{(m_{1})}(\gam_{m_{1}})
&\cip-\frac{1}{2}\int_{\mbbr^{d}}\left\{\trace\left(S(x,\gam_{m_{1}})^{-1}S(x)\right)+\log\left|S(x,\gam_{m_{1}})\right|\right\}\pi(dx) \\
&=:\mbbh_{1,0}^{(m_{1})}(\gam_{m_{1}}),
\end{align*}
where $S(x)=C(x)^{\otimes2}$.
We assume that the optimal scale parameter $\gam_{m_{1}}^{\ast}$ and scale index set $\mathfrak{M}_{1}^{\ast}$ are defined as
\begin{align*}
\{\gam_{m_{1}}^{\ast}\}&=\argmax_{\gam_{m_{1}}}\mbbh_{1,0}^{(m_{1})}(\gam_{m_{1}}), \\
\mathfrak{M}_{1}^{\ast}&=\argmin_{m_{1}\in\mathfrak{M}_{1}}\dim(\Theta_{\gam_{m_{1}}}),
\end{align*}
respectively.
For any $m_{1}\in\mathfrak{M}_{1}$, $\gam_{m_{1}}^{\ast}=\gam_{m_{1},0}$.

Next, for any fixed $m_{1}\in\{1,\ldots,M_{1}\}$,
\begin{align*}
\frac{1}{T_{n}}\mbbh_{2,n}^{(m_{2}|m_{1})}(\al_{m_{2}})&\cip-\frac{1}{2}\int_{\mbbr^{d}}S^{-1}(x,\gam_{m_{1}}^{\ast})\left[\left(a_{m_{2}}(x,\al_{m_{2}})-A(x) \right)^{\otimes 2}\right]\pi(dx) \\
&=:\mbbh_{2,0}^{(m_{2}|m_{1})}(\al_{m_{2}}),
\end{align*}
and assume that the optimal drift parameter $\al_{m_{2}}^{\ast}$ is given by maximizing $\mbbh_{2,0}^{(m_{2}|m_{1})}$:
\begin{align*}
\{\al_{m_{2}}^{\ast}\}&=\argmax_{\al_{m_{2}}}\mbbh_{2,0}^{(m_{2}|m_{1})}(\al_{m_{2}}).
\end{align*}
When $m_{2}$ is included in $\mathfrak{M}_{2}$, $\al_{m_{2}}^{\ast}=\al_{m_{2},0}$.
We also suppose that the drift index set $\mathfrak{M}_{2}^{\ast}$ is defined as
\begin{align*} 
\mathfrak{M}_{2}^{\ast}&=\argmin_{m_{2}\in\mathfrak{M}_{2}}\dim(\Theta_{\al_{m_{2}}}).
\end{align*}
From the assumptions and definitions of $\mbbh_{1,0}^{(m_{1})}$ and $\mbbh_{2,0}^{(m_{2}|m_{1})}$, $\mathfrak{M}_{1}=\argmax_{m_{1}}\mbbh_{1,0}^{(m_{1})}(\gam_{m_{1}}^{\ast})$ and $\mathfrak{M}_{2}=\argmax_{m_{2}}\mbbh_{2,0}^{(m_{2}|m_{1})}(\al_{m_{2}}^{\ast})$ hold.

Let $\Theta_{\gamma_{i_{1}}}\times\Theta_{\alpha_{i_{2}}}\subset\mbbr^{p_{\gamma_{i_{1}}}}\times\mbbr^{p_{\alpha_{i_{2}}}}$ and $\Theta_{\gamma_{j_{1}}}\times\Theta_{\alpha_{j_{2}}}\subset\mbbr^{p_{\gamma_{j_{1}}}}\times\mbbr^{p_{\alpha_{j_{2}}}}$ be the parameter space associated with model $\mcm_{i_{1},i_{2}}$ and $\mcm_{j_{1},j_{2}}$, respectively.
%If $p_{\gamma_{i_{1}}}<p_{\gamma_{j_{1}}}$ and there exista a matrix $F_{1}\in\mbbr^{p_{\gamma_{j_{1}}}\times p_{\gamma_{i_{1}}}}$ with $F_{1}^{\top}F_{1}=I_{p_{\gamma_{i_{1}}}\times p_{\gamma_{j_{1}}}}$ as well as a $c\in\mbbr^{p_{\gamma_{j_{1}}}}$ such that $c_{i_{1}}(\gam_{i_{1}})=c_{j_{1}}(F_{1}\gam_{i_{1}}+c_{1})$ for all $\gamma_{i_{1}}\in\Theta_{\gamma_{i_{1}}}$, we say $\Theta_{\gamma_{i_{1}}}$ is nested in $\Theta_{\gamma_{j_{1}}}$.
If $p_{\gamma_{i_{1}}}<p_{\gamma_{j_{1}}}$ and there exists a matrix $F_{1}\in\mbbr^{p_{\gamma_{j_{1}}}\times p_{\gamma_{i_{1}}}}$ with $F_{1}^{\top}F_{1}=I_{p_{\gamma_{i_{1}}}\times p_{\gamma_{i_{1}}}}$ as well as a $\mbbh_{1,n}^{(i_{1})}(\gam_{i_{1}})=\mbbh_{1,n}^{(j_{1})}(F_{1}\gamma_{i_{1}}+c_{1})$ for all $\gamma_{i_{1}}\in\Theta_{\gamma_{i_{1}}}$, we say $\Theta_{\gamma_{i_{1}}}$ is nested in $\Theta_{\gamma_{j_{1}}}$.
%When $\Theta_{\gamma_{i_{1}}}$ is nested in $\Theta_{\gamma_{j_{1}}}$, the equation $\mbbh_{1,n}^{(i_{1})}(\gam_{i_{1}})=\mbbh_{1,n}^{(j_{1})}(F_{1}\gamma_{i_{1}}+c_{1})$ holds.
It is defined in a similar manner that $\Theta_{\alpha_{i_{2}}}$ is nested in $\Theta_{\alpha_{j_{2}}}$.

\medskip

Now we are in a position to state the results. The theoretical properties of the GQAIC are given in Theorem \ref{se:thm_sele.prob}, and those of the GQBIC in Theorems \ref{se:thm_sele.prob2} and \ref{se:thm_sele.prob3}.
For convenience, we give the summaries before the statements:
\begin{itemize}
\item Theorem \ref{se:thm_sele.prob} 1(i) and 2(i) reveal that the probability of relative selection is asymptotically characterized by the non-central chi-squared distribution; in general, this happens when an estimator under consideration is asymptotically normally distributed with the asymptotic covariance matrix being of the sandwich form (see \cite{Ken82-2}).
Further, Theorem \ref{se:thm_sele.prob} 1(ii) and 2(ii) indicate that the probability that GQAIC chooses the misspecified coefficients tends to 0 as $n\to\infty$.
\item Theorem \ref{se:thm_sele.prob2} 1(i) shows that, when comparing correctly specified models, the probability that $\mathrm{GQBIC}_{1,n}^{\sharp}$ selects a larger model tends to 1.
Moreover, Theorem \ref{se:thm_sele.prob3} means that the GQBIC proposed by \eqref{hm:GQBIC_gam} and \eqref{hm:GQBIC_al} has the model selection consistency.
\end{itemize}
Let $\mathrm{GQBIC}_{1,n}^{(m_{1})}$ and $\mathrm{GQBIC}_{1,n}^{\sharp(m_{1})}$ denote the $\mathrm{GQBIC}_{1,n}$ and $\mathrm{GQBIC}_{1,n}^{\sharp}$ of the $m_{1}$-th candidate scale coefficient, respectively. Also, let $\mathrm{GQBIC}_{2,n}^{(m_{2}|m_{1,n})}$ correspond to \eqref{hm:GQBIC_al} associated with $c_{m_{1,n}}$ and $\hat{\gam}_{m_{1,n},n}$.

\begin{thm}
\label{se:thm_sele.prob}
Suppose that the assumptions of Theorem \ref{hm:thm_GQAIC.stepwise.se} hold for all candidate coefficients which are included in $\mathfrak{M}_{1}$ and $\mathfrak{M}_{2}$.
We also assume that indexes $m_{1}^{\ast}$ and $m_{2}^{\ast}$ satisfy $m_{1}^{\ast}\in\mathfrak{M}_{1}^{\ast}$ and $m_{2}^{\ast}\in\mathfrak{M}_{2}^{\ast}$, respectively.

\begin{itemize}
\item[1.] 
\begin{itemize}
\item[(i)] Let $m_{1}\in\mathfrak{M}_{1}\backslash\{m_{1}^{\ast}\}$.
If $\Theta_{\gam_{m_{1}^{\ast}}}$ is nested in $\Theta_{\gam_{m_{1}}}$ with map $F_{1}$, then
\begin{align*}
&\lim _{n\to\infty}\pr\left(\gqaic_{1,n}^{(m_{1}^{\ast})}-\gqaic_{1,n}^{(m_{1})}>0\right) \\
&=\pr\bigg[\sum_{j=1}^{p_{\gamma_{m_{1}}}}\lambda_{j}\chi_{j}^{2}>2\trace\left\{\Gam_{\gam_{m_{1}}}(\gam_{m_{1},0})^{-1}W_{\gam_{m_{1}}}(\gam_{m_{1},0})\right\} \\
&\hspace{40mm}-2\trace\left\{\Gam_{\gam_{m_{1}^{\ast}}}(\gam_{m_{1}^{\ast},0})^{-1}W_{\gam_{m_{1}^{\ast}}}(\gam_{m_{1}^{\ast},0})\right\}\bigg] \\
&>0,
\end{align*}
where 
\begin{align*}
G_{\gamma_{m_{1}}}(\gam_{m_{1},0})=\Gamma_{\gam_{m_{1}}}(\gam_{m_{1},0})^{-1}-F_{1}\left(F_{1}^{\top}\Gamma_{\gam_{m_{1}}}(\gam_{m_{1},0})F_{1}\right)^{-1}F_{1}^{\top},
\end{align*}
$(\chi_{j}^{2})$ is a sequence of independent $\chi^{2}$ random variables with one degree of freedom, and $\lambda_{1},\lambda_{2},\ldots,\lambda_{p_{\gamma_{m_{1}}}}$ are the the eigenvalues of $W_{\gamma_{m_{1}}}(\gam_{m_{1},0})^{1/2}G_{\gamma_{m_{1}}}(\gam_{m_{1},0})W_{\gamma_{m_{1}}}(\gam_{m_{1},0})^{1/2}$.

\item[(ii)] If $m_{1}\in\{1,\ldots,M_{1}\}\backslash\mathfrak{M}_{1}$, then
\begin{align*}
\lim _{n\to\infty}P\left(\gqaic_{1,n}^{(m_{1}^{\ast})} < \gqaic_{1,n}^{(m_{1})}\right)=1.
\end{align*}
\end{itemize}

\item[2.]
\begin{itemize}
\item[(i)] Let $m_{2}\in\mathfrak{M}_{2}\backslash\{m_{2}^{\ast}\}$. 
If $\Theta_{\al_{m_{2}^{\ast}}}$ is nested in $\Theta_{\al_{m_{2}}}$ with map $F_{2}$, then
\begin{align*}
&\lim _{n\to\infty}\pr\left(\gqaic_{2,n}^{(m_{2}^{\ast}|\hat{m}_{1,n})}-\gqaic_{2,n}^{(m_{2}|\hat{m}_{1,n})}>0\right)\\
&=\pr\bigg[\sum_{j=1}^{p_{\al_{m_{2}}}}\lambda_{j}^{\prime}\chi_{j}^{2}>2(p_{\al_{m_{2}}}-p_{\al_{m_{2}^{\ast}}})\bigg]>0,
\end{align*}
where 
\begin{align*}
G_{\al_{m_{2}}}(\al_{m_{2},0},\gam_{\hat{m}_{1,n},0})=\Gam_{\al_{m_{2}}}(\al_{m_{2},0},\gam_{\hat{m}_{1,n},0})^{-1}-F_{2}\left(F_{2}^{\top}\Gam_{\al_{m_{2}}}(\al_{m_{2},0},\gam_{\hat{m}_{1,n},0})F_{2}\right)^{-1}F_{2}^{\top}
\end{align*}
and $\lambda_{1}^{\prime},\lambda_{2}^{\prime},\ldots,\lambda_{p_{\al_{m_{2}}}}^{\prime}$ are the the eigenvalues of $\Gam_{\al_{m_{2}}}(\al_{m_{2},0},\gam_{\hat{m}_{1,n},0})^{1/2}G_{\al_{m_{2}}}(\al_{m_{2},0},\gam_{\hat{m}_{1,n},0})\allowbreak\Gam_{\al_{m_{2}}}(\al_{m_{2},0},\gam_{\hat{m}_{1,n},0})^{1/2}$.

\item[(ii)]If $m_{2}\in\{1,\ldots,M_{2}\}\backslash\mathfrak{M}_{2}$, then
\begin{align*}
\lim _{n\to\infty}\pr\left(\gqaic_{2,n}^{(m_{2}^{\ast}|\hat{m}_{1,n})} < \gqaic_{2,n}^{(m_{2}|\hat{m}_{1,n})}\right)=1.
\end{align*}
\end{itemize}
\end{itemize}
\end{thm}

\begin{thm}
\label{se:thm_sele.prob2}
Suppose that the assumptions of Theorem \ref{hm:thm_GQAIC.stepwise.se} hold for all candidate coefficients which are included in $\mathfrak{M}_{1}$.
We also assume that index $m_{1}^{\ast}$ satisfies $m_{1}^{\ast}\in\mathfrak{M}_{1}^{\ast}$.

\begin{itemize}
\item[(i)] Let $m_{1}\in\mathfrak{M}_{1}\backslash\{m_{1}^{\ast}\}$.
If $\Theta_{\gam_{m_{1}^{\ast}}}$ is nested in $\Theta_{\gam_{m_{1}}}$, then
\begin{align*}
&\lim _{n\to\infty}\pr\left(\mathrm{GQBIC}_{1,n}^{\sharp(m_{1}^{\ast})}>\mathrm{GQBIC}_{1,n}^{\sharp(m_{1})}\right)=1.
\end{align*}

\item[(ii)] If $m_{1}\in\{1,\ldots,M_{1}\}\backslash\mathfrak{M}_{1}$, then
\begin{align*}
&\lim _{n\to\infty}\pr\left(\mathrm{GQBIC}_{1,n}^{\sharp(m_{1}^{\ast})} < \mathrm{GQBIC}_{1,n}^{\sharp(m_{1})}\right)=1.
\end{align*}
\end{itemize}
\end{thm}

\begin{thm}
\label{se:thm_sele.prob3}
Suppose that the assumptions of Theorem \ref{hm:thm_GQAIC.stepwise.se} hold for all candidate coefficients which are included in $\mathfrak{M}_{1}$ and $\mathfrak{M}_{2}$.
We also assume that indexes $m_{1}^{\ast}$ and $m_{2}^{\ast}$ satisfy $m_{1}^{\ast}\in\mathfrak{M}_{1}^{\ast}$ and $m_{2}^{\ast}\in\mathfrak{M}_{2}^{\ast}$, respectively.

\begin{itemize}
\item[1.] 
\begin{itemize}
\item[(i)] Let $m_{1}\in\mathfrak{M}_{1}\backslash\{m_{1}^{\ast}\}$.
If $\Theta_{\gam_{m_{1}^{\ast}}}$ is nested in $\Theta_{\gam_{m_{1}}}$, then
\begin{align*}
&\lim _{n\to\infty}\pr\left(\mathrm{GQBIC}_{1,n}^{(m_{1}^{\ast})}<\mathrm{GQBIC}_{1,n}^{(m_{1})}\right)=1.
\end{align*}
\item[(ii)] If $m_{1}\in\{1,\ldots,M_{1}\}\backslash\mathfrak{M}_{1}$, then
\begin{align*}
&\lim _{n\to\infty}\pr\left(\mathrm{GQBIC}_{1,n}^{(m_{1}^{\ast})} < \mathrm{GQBIC}_{1,n}^{(m_{1})}\right)=1.
\end{align*}
\end{itemize}
\item[2.]
\begin{itemize}
\item[(i)] Let $m_{2}\in\mathfrak{M}_{2}\backslash\{m_{2}^{\ast}\}$. 
If $\Theta_{\al_{m_{2}^{\ast}}}$ is nested in $\Theta_{\al_{m_{2}}}$, then
\begin{align*}
&\lim _{n\to\infty}\pr\left(\mathrm{GQBIC}_{2,n}^{(m_{2}^{\ast}|\hat{m}_{1,n})} < \mathrm{GQBIC}_{2,n}^{(m_{2}|\hat{m}_{1,n})}\right)=1.
\end{align*}
\item[(ii)]If $m_{2}\in\{1,\ldots,M_{2}\}\backslash\mathfrak{M}_{2}$, then
\begin{align*}
&\lim _{n\to\infty}\pr\left(\mathrm{GQBIC}_{2,n}^{(m_{2}^{\ast}|\hat{m}_{1,n})}<\mathrm{GQBIC}_{2,n}^{(m_{2}|\hat{m}_{1,n})}\right)=1.
\end{align*}
\end{itemize}
\end{itemize}
\end{thm}

%%%%%
%%%%%
\section{Numerical experiments}
\label{hm:sec_sim}

In this section, we present simulation results to observe finite-sample performance of the proposed GQBIC and GQAIC.
We use the {\tt yuima} package on R (see \cite{YUIMA14}) for generating data.
All the Monte Carlo trials are based on 1000 independent sample paths, and the simulations are done for $(h,T_{n})=(0.01,10),$ $(0.005,10), (0.01,50)$, and $(0.005,50)$ (hence in each case, $n=1000, 2000, 5000$, and $10000$).

The sample data $\X_{n}=(X_{t_{j}})_{j=0}^{n}$ with $t_{j}=jh$ is obtained from
\begin{align*}
dX_{t}=-\frac{1}{2}X_{t}dt+\frac{3}{1+X_{t-}^{2}}dZ_{t},\quad t\in[0,T_{n}],\quad X_{0}=0,
\end{align*}
where $T_{n}=nh_{n}$.
The numerical experiments are conducted in three situations:
\begin{itemize}
\item[(i)] $\mathcal{L}(Z_{t})=NIG(10,0,10t,0)$,
\item[(ii)] $\mathcal{L}(Z_{t})=bGamma(t,\sqrt{2},t,\sqrt{2})$,
\item[(iii)] $\mathcal{L}(Z_{t})=NIG\left(\frac{25}{3},\frac{20}{3},\frac{9}{5}t,-\frac{12}{5}t\right)$.
% (iv) $\mathcal{L}(Z_{t})=NIG\left(5,4,\frac{27}{25}t,-\frac{36}{25}t\right)$.
%The density functions of (i) -- (iii) are drawn in Figure \ref{density}.
\end{itemize}
Here $NIG$ and $bGamma$ refer to the normal inverse-Gaussian and bilateral gamma distributions, respectively (see \cite{IacYos18} for the definitions).
In this example, we consider the following candidate scale (Scale) and drift (Drift) coefficients:
\begin{align*}
&\;{\bf Scale}\;{\bf 1:} c_{1}(x,\gamma_{1})=\gamma_{1};
\;{\bf Scale}\;{\bf 2:} c_{2}(x,\gamma_{2})=\frac{\gamma_{2}}{1+x^{2}}; \\
&\;{\bf Scale}\;{\bf 3:} c_{3}(x,\gamma_{3})=\frac{\gamma_{3,1}+\gamma_{3,2}x^{2}}{1+x^{2}};
\;{\bf Scale}\;{\bf 4:} c_{4}(x,\gamma_{4})=\frac{\gamma_{4,1}+\gamma_{4,2}x+\gamma_{4,3}x^{2}}{1+x^{2}},
\end{align*}
and
\begin{align*}
{\bf Drift}\;{\bf 1:}\; a_{1}(x,\alpha_{1})=-\alpha_{1};
\;{\bf Drift}\;{\bf 2:}\; a_{2}(x,\alpha_{2})=-\alpha_{2}x;
\;{\bf Drift}\;{\bf 3:}\; a_{3}(x,\alpha_{3})=-\alpha_{3,1}x-\alpha_{3,2}.
\end{align*}
Each candidate model is given by a combination of the scale and drift coefficients, 
and the stochastic differential equation models based on local Gaussian approximate models of L\'evy stochastic differential equations are assumed as misspecified candidate models.
For example, in the case of Scale 1 and Drift 1, the statistical model is
%not L\'evy driven SDE
%\begin{align*}
%dX_{t}=-\al_{1}dt+\gam_{1}dZ_{t},
%\end{align*}
%but 
a stochastic differential equation model given by
\begin{align*}
\mathcal{L}(X_{t_{j}}|X_{t_{j-1}}=x) = N\left(x+\al_{1}h,h\gam_{1}^{2}\right).
\end{align*}
\color{black}
Then, the Scale 2 and Drift 2 with $(\gamma_{2},\alpha_{2})=(3,\frac{1}{2})$ are the true coefficients, and the coefficients Scale 3, 4, and Drift 3 include the true coefficient.

We compare model selection frequency through GQAIC, GQBIC, and $\mathrm{GQBIC}^{\sharp}$.
Also, we formally use classical AIC, say formal AIC(fAIC), and compare the model selection results with those of the proposed criteria.
The fAIC for scale and drift are given by
\begin{align*}
\mathrm{fAIC}_{1,n}&=-2\mbbh_{1,n}(\ges)+2p_{\gam}, \\
\mathrm{fAIC}_{2,n}&=-2\mbbh_{2,n}(\hat{\alpha}_{n})+2p_{\al},
\end{align*}
respectively.
\color{black}
Tables \ref{simu:tab1}, \ref{simu:tab2}, and \ref{simu:tab3} summarize the comparison results of model selection frequency.
The GQAIC and GQBIC select the true coefficients Scale 2 and Drift 2 with high frequency in all cases, while the fAIC selects the true coefficients with less frequency in (ii) and (iii) cases.
Moreover, we can observe the frequencies of selecting the true coefficient by GQBIC become larger as sample size $n$ increases.
In the (ii) and (iii) cases, also observed is that the frequencies that the misspecified coefficients Scale 1 and Drift 1 are chosen by GQAIC become lower as $n$ increases.

%%%%%
%%%%%

\begin{table}[t]
\begin{center}
\caption{Computation results of (i) case. 
Model selection frequencies for various situations are shown. The true model consists of Scale 2 and Drift 2.}
\scalebox{0.9}[0.9]{
\begin{tabular}{l r r l l r r r r} \hline
fAIC & $T_{n}$ & $h$ & & & Scale 1 & Scale $2^{\ast}$ & Scale 3 & Scale 4 \\ \hline
& 10 & 0.01 & Drift 1 & & 0 & 1 & 0 & 0 \\[1mm]
& \multicolumn{2}{r}{$(n=1000)$} & Drift $2^{\ast}$ & & 0 & {\bf 536} & 124 & 199 \\[1mm]
& & & Drift 3 & & 0 & 74 & 27 & 39 \\[1mm]
& 10 & 0.005 & Drift 1 & & 0 & 0 & 1 & 0 \\[1mm]
& \multicolumn{2}{r}{$(n=2000)$} & Drift $2^{\ast}$ & & 0 & {\bf 458} & 110 & 285 \\[1mm]
& & & Drift 3 & & 0 & 67 & 25 & 54 \\[1mm]
& 50 & 0.01 & Drift 1 & & 0 & 0 & 0 & 0 \\[1mm]
& \multicolumn{2}{r}{$(n=5000)$} & Drift $2^{\ast}$ & & 0 & {\bf 466} & 202 & 231 \\[1mm]
& & & Drift 3 & & 0 & 51 & 28 & 32 \\[1mm]
& 50 & 0.005 & Drift 1 & & 0 & 0 & 0 & 0 \\[1mm]
& \multicolumn{2}{r}{$(n=10000)$} & Drift $2^{\ast}$ & & 0 & {\bf 402} & 150 & 352 \\[1mm]
& & & Drift 3 & & 0 & 51 & 13 & 32 \\[1mm] \hline\hline
%%%%%
GQAIC & $T_{n}$ & $h_{n}$ & & & Scale 1 & Scale $2^{\ast}$ & Scale 3 & Scale 4 \\ \hline
& 10 & 0.01 & Drift 1 & & 0 & 1 & 0 & 0 \\[1mm]
& \multicolumn{2}{r}{$(n=1000)$} & Drift $2^{\ast}$ & & 0 & {\bf 714} & 77 & 64 \\[1mm]
& & & Drift 3 & & 0 & 110 & 20 & 14 \\[1mm]
& 10 & 0.005 & Drift 1 & & 0 & 1 & 0 & 0 \\[1mm]
& \multicolumn{2}{r}{$(n=2000)$} & Drift $2^{\ast}$ & & 0 & {\bf 733} & 62 & 56 \\[1mm]
& & & Drift 3 & & 0 & 117 & 17 & 14 \\[1mm]
& 50 & 0.01 & Drift 1 & & 0 & 0 & 0 & 0 \\[1mm]
& \multicolumn{2}{r}{$(n=5000)$} & Drift $2^{\ast}$ & & 0 & {\bf 713} & 122 & 64 \\[1mm]
& & & Drift 3 & & 0 & 83 & 11 & 7 \\[1mm]
& 50 & 0.005 & Drift 1 & & 0 & 0 & 0 & 0 \\[1mm]
& \multicolumn{2}{r}{$(n=10000)$} & Drift $2^{\ast}$ & & 0 & {\bf 765} & 79 & 59 \\[1mm]
& & & Drift 3 & & 0 & 88 & 5 & 4 \\[1mm] \hline\hline
%%%%%
GQBIC & $T_{n}$ & $h_{n}$ & & & Scale 1 & Scale $2^{\ast}$ & Scale 3 & Scale 4 \\ \hline
& 10 & 0.01 & Drift 1 & & 0 & 1 & 0 & 0 \\[1mm]
& \multicolumn{2}{r}{$(n=1000)$} & Drift $2^{\ast}$ & & 0 & {\bf 861} & 0 & 0 \\[1mm]
& & & Drift 3 & & 0 & 138 & 0 & 0 \\[1mm]
& 10 & 0.005 & Drift 1 & & 0 & 1 & 0 & 0 \\[1mm]
& \multicolumn{2}{r}{$(n=2000)$} & Drift $2^{\ast}$ & & 0 & {\bf 866} & 0 & 0 \\[1mm]
& & & Drift 3 & & 0 & 133 & 0 & 0 \\[1mm]
& 50 & 0.01 & Drift 1 & & 0 & 0 & 0 & 0 \\[1mm]
& \multicolumn{2}{r}{$(n=5000)$} & Drift $2^{\ast}$ & & 0 & {\bf 965} & 0 & 0 \\[1mm]
& & & Drift 3 & & 0 & 35 & 0 & 0 \\[1mm]
& 50 & 0.005 & Drift 1 & & 0 & 0 & 0 & 0 \\[1mm]
& \multicolumn{2}{r}{$(n=10000)$} & Drift $2^{\ast}$ & & 0 & {\bf 964} & 0 & 0 \\[1mm]
& & & Drift 3 & & 0 & 36 & 0 & 0 \\[1mm] \hline\hline
%%%%%
$\mathrm{GQBIC}^{\sharp}$ & $T_{n}$ & $h_{n}$ & & & Scale 1 & Scale $2^{\ast}$ & Scale 3 & Scale 4 \\ \hline
& 10 & 0.01 & Drift 1 & & 0 & 1 & 0 & 0 \\[1mm]
& \multicolumn{2}{r}{$(n=1000)$} & Drift $2^{\ast}$ & & 0 & {\bf 788} & 49 & 30 \\[1mm]
& & & Drift 3 & & 0 & 119 & 7 & 6 \\[1mm]
& 10 & 0.005 & Drift 1 & & 0 & 1 & 0 & 0 \\[1mm]
& \multicolumn{2}{r}{$(n=2000)$} & Drift $2^{\ast}$ & & 0 & {\bf 759} & 69 & 45 \\[1mm]
& & & Drift 3 & & 0 & 107 & 9 & 10 \\[1mm]
& 50 & 0.01 & Drift 1 & & 0 & 0 & 0 & 0 \\[1mm]
& \multicolumn{2}{r}{$(n=5000)$} & Drift $2^{\ast}$ & & 0 & {\bf 882} & 64 & 19 \\[1mm]
& & & Drift 3 & & 0 & 32 & 2 & 1 \\[1mm]
& 50 & 0.005 & Drift 1 & & 0 & 0 & 0 & 0 \\[1mm]
& \multicolumn{2}{r}{$(n=10000)$} & Drift $2^{\ast}$ & & 0 & {\bf 862} & 68 & 34 \\[1mm]
& & & Drift 3 & & 0 & 32 & 2 & 2 \\[1mm] \hline
\end{tabular}
}
\label{simu:tab1}
\end{center}
\end{table}

\begin{table}[t]
\begin{center}
\caption{Computation results of (ii) case. 
Model selection frequencies for various situations are shown. The true model consists of Scale 2 and Drift 2.}
\scalebox{0.9}[0.9]{
\begin{tabular}{l r r l l r r r r} \hline
fAIC & $T_{n}$ & $h$ & & & Scale 1 & Scale $2^{\ast}$ & Scale 3 & Scale 4 \\ \hline
& 10 & 0.01 & Drift 1 & & 0 & 0 & 0 & 9 \\[1mm]
& \multicolumn{2}{r}{$(n=1000)$} & Drift $2^{\ast}$ & & 5 & {\bf 106} & 26 & 767 \\[1mm]
& & & Drift 3 & & 0 & 10 & 1 & 76 \\[1mm]
& 10 & 0.005 & Drift 1 & & 0 & 0 & 0 & 9 \\[1mm]
& \multicolumn{2}{r}{$(n=2000)$} & Drift $2^{\ast}$ & & 2 & {\bf 83} & 24 & 793 \\[1mm]
& & & Drift 3 & & 0 & 6 & 2 & 81 \\[1mm]
& 50 & 0.01 & Drift 1 & & 0 & 0 & 0 & 0 \\[1mm]
& \multicolumn{2}{r}{$(n=5000)$} & Drift $2^{\ast}$ & & 2 & {\bf 85} & 47 & 782 \\[1mm]
& & & Drift 3 & & 0 & 6 & 4 & 74 \\[1mm]
& 50 & 0.005 & Drift 1 & & 0 & 0 & 0 & 0 \\[1mm]
& \multicolumn{2}{r}{$(n=10000)$} & Drift $2^{\ast}$ & & 0 & {\bf 68} & 25 & 826 \\[1mm]
& & & Drift 3 & & 0 & 4 & 3 & 74 \\[1mm] \hline\hline
%%%%%
GQAIC & $T_{n}$ & $h_{n}$ & & & Scale 1 & Scale $2^{\ast}$ & Scale 3 & Scale 4 \\ \hline
& 10 & 0.01 & Drift 1 & & 1 & 0 & 0 & 4 \\[1mm]
& \multicolumn{2}{r}{$(n=1000)$} & Drift $2^{\ast}$ & & 97 & {\bf 591} & 10 & 186 \\[1mm]
& & & Drift 3 & & 15 & 86 & 0 & 10 \\[1mm]
& 10 & 0.005 & Drift 1 & & 1 & 0 & 0 & 5 \\[1mm]
& \multicolumn{2}{r}{$(n=2000)$} & Drift $2^{\ast}$ & & 99 & {\bf 584} & 7 & 191 \\[1mm]
& & & Drift 3 & & 15 & 88 & 0 & 10 \\[1mm]
& 50 & 0.01 & Drift 1 & & 0 & 0 & 0 & 0 \\[1mm]
& \multicolumn{2}{r}{$(n=5000)$} & Drift $2^{\ast}$ & & 29 & {\bf 741} & 36 & 100 \\[1mm]
& & & Drift 3 & & 4 & 73 & 5 & 12 \\[1mm]
& 50 & 0.005 & Drift 1 & & 0 & 0 & 0 & 0 \\[1mm]
& \multicolumn{2}{r}{$(n=10000)$} & Drift $2^{\ast}$ & & 27 & {\bf 747} & 27 & 104 \\[1mm]
& & & Drift 3 & & 4 & 74 & 5 & 12 \\[1mm] \hline\hline
%%%%%
GQBIC & $T_{n}$ & $h_{n}$ & & & Scale 1 & Scale $2^{\ast}$ & Scale 3 & Scale 4 \\ \hline
& 10 & 0.01 & Drift 1 & & 3 & 0 & 0 & 1 \\[1mm]
& \multicolumn{2}{r}{$(n=1000)$} & Drift $2^{\ast}$ & & 142 & {\bf 700} & 2 & 36 \\[1mm]
& & & Drift 3 & & 14 & 101 & 0 & 1 \\[1mm]
& 10 & 0.005 & Drift 1 & & 3 & 0 & 0 & 2 \\[1mm]
& \multicolumn{2}{r}{$(n=2000)$} & Drift $2^{\ast}$ & & 142 & {\bf 700} & 1 & 37 \\[1mm]
& & & Drift 3 & & 14 & 100 & 0 & 1 \\[1mm]
& 50 & 0.01 & Drift 1 & & 0 & 0 & 0 & 0 \\[1mm]
& \multicolumn{2}{r}{$(n=5000)$} & Drift $2^{\ast}$ & & 42 & {\bf 890} & 9 & 22 \\[1mm]
& & & Drift 3 & & 3 & 33 & 0 & 1 \\[1mm]
& 50 & 0.005 & Drift 1 & & 0 & 0 & 0 & 0 \\[1mm]
& \multicolumn{2}{r}{$(n=10000)$} & Drift $2^{\ast}$ & & 38 & {\bf 894} & 6 & 23 \\[1mm]
& & & Drift 3 & & 3 & 35 & 0 & 1 \\[1mm] \hline\hline
%%%%%
$\mathrm{GQBIC}^{\sharp}$ & $T_{n}$ & $h_{n}$ & & & Scale 1 & Scale $2^{\ast}$ & Scale 3 & Scale 4 \\ \hline
& 10 & 0.01 & Drift 1 & & 0 & 0 & 0 & 10 \\[1mm]
& \multicolumn{2}{r}{$(n=1000)$} & Drift $2^{\ast}$ & & 15 & {\bf 189} & 32 & 685 \\[1mm]
& & & Drift 3 & & 2 & 14 & 2 & 51 \\[1mm]
& 10 & 0.005 & Drift 1 & & 0 & 0 & 0 & 10 \\[1mm]
& \multicolumn{2}{r}{$(n=2000)$} & Drift $2^{\ast}$ & & 10 & {\bf 130} & 35 & 743 \\[1mm]
& & & Drift 3 & & 1 & 12 & 2 & 57 \\[1mm]
& 50 & 0.01 & Drift 1 & & 0 & 0 & 0 & 0 \\[1mm]
& \multicolumn{2}{r}{$(n=5000)$} & Drift $2^{\ast}$ & & 3 & {\bf 200} & 84 & 686 \\[1mm]
& & & Drift 3 & & 0 & 4 & 2 & 21 \\[1mm]
& 50 & 0.005 & Drift 1 & & 0 & 0 & 0 & 0 \\[1mm]
& \multicolumn{2}{r}{$(n=10000)$} & Drift $2^{\ast}$ & & 1 & {\bf 159} & 49 & 763 \\[1mm]
& & & Drift 3 & & 0 & 3 & 2 & 23 \\[1mm] \hline
\end{tabular}
}
\label{simu:tab2}
\end{center}
\end{table}

\begin{table}[t]
\begin{center}
\caption{Computation results of (iii) case. 
Model selection frequencies for various situations are shown. The true model consists of Scale 2 and Drift 2.}
\scalebox{0.9}[0.9]{
\begin{tabular}{l r r l l r r r r} \hline
fAIC & $T_{n}$ & $h$ & & & Scale 1 & Scale $2^{\ast}$ & Scale 3 & Scale 4 \\ \hline
& 10 & 0.01 & Drift 1 & & 0 & 0 & 0 & 0 \\[1mm]
& \multicolumn{2}{r}{$(n=1000)$} & Drift $2^{\ast}$ & & 0 & {\bf 128} & 52 & 600 \\[1mm]
& & & Drift 3 & & 0 & 30 & 30 & 160 \\[1mm]
& 10 & 0.005 & Drift 1 & & 0 & 0 & 0 & 1 \\[1mm]
& \multicolumn{2}{r}{$(n=2000)$} & Drift $2^{\ast}$ & & 0 & {\bf 88} & 37 & 653 \\[1mm]
& & & Drift 3 & & 0 & 25 & 16 & 180 \\[1mm]
& 50 & 0.01 & Drift 1 & & 0 & 0 & 0 & 0 \\[1mm]
& \multicolumn{2}{r}{$(n=5000)$} & Drift $2^{\ast}$ & & 0 & {\bf 88} & 75 & 731 \\[1mm]
& & & Drift 3 & & 0 & 12 & 4 & 9 \\[1mm]
& 50 & 0.005 & Drift 1 & & 0 & 0 & 0 & 0 \\[1mm]
& \multicolumn{2}{r}{$(n=10000)$} & Drift $2^{\ast}$ & & 0 & {\bf 70} & 39 & 785 \\[1mm]
& & & Drift 3 & & 0 & 6 & 1 & 99 \\[1mm] \hline\hline
%%%%%
GQAIC & $T_{n}$ & $h_{n}$ & & & Scale 1 & Scale $2^{\ast}$ & Scale 3 & Scale 4 \\ \hline
& 10 & 0.01 & Drift 1 & & 0 & 0 & 0 & 0 \\[1mm]
& \multicolumn{2}{r}{$(n=1000)$} & Drift $2^{\ast}$ & & 18 & {\bf 589} & 33 & 155 \\[1mm]
& & & Drift 3 & & 32 & 125 & 33 & 15 \\[1mm]
& 10 & 0.005 & Drift 1 & & 0 & 0 & 0 & 1 \\[1mm]
& \multicolumn{2}{r}{$(n=2000)$} & Drift $2^{\ast}$ & & 15 & {\bf 602} & 30 & 149 \\[1mm]
& & & Drift 3 & & 32 & 122 & 32 & 17 \\[1mm]
& 50 & 0.01 & Drift 1 & & 0 & 0 & 0 & 0 \\[1mm]
& \multicolumn{2}{r}{$(n=5000)$} & Drift $2^{\ast}$ & & 0 & {\bf 672} & 102 & 122 \\[1mm]
& & & Drift 3 & & 0 & 80 & 8 & 16 \\[1mm]
& 50 & 0.005 & Drift 1 & & 0 & 0 & 0 & 0 \\[1mm]
& \multicolumn{2}{r}{$(n=10000)$} & Drift $2^{\ast}$ & & 0 & {\bf 710} & 77 & 110 \\[1mm]
& & & Drift 3 & & 0 & 78 & 4 & 21 \\[1mm] \hline\hline
%%%%%
GQBIC & $T_{n}$ & $h_{n}$ & & & Scale 1 & Scale $2^{\ast}$ & Scale 3 & Scale 4 \\ \hline
& 10 & 0.01 & Drift 1 & & 0 & 0 & 0 & 0 \\[1mm]
& \multicolumn{2}{r}{$(n=1000)$} & Drift $2^{\ast}$ & & 24 & {\bf 798} & 5 & 5 \\[1mm]
& & & Drift 3 & & 35 & 133 & 0 & 0 \\[1mm]
& 10 & 0.005 & Drift 1 & & 0 & 0 & 0 & 0 \\[1mm]
& \multicolumn{2}{r}{$(n=2000)$} & Drift $2^{\ast}$ & & 25 & {\bf 795} & 4 & 4 \\[1mm]
& & & Drift 3 & & 37 & 135 & 0 & 0 \\[1mm]
& 50 & 0.01 & Drift 1 & & 0 & 0 & 0 & 0 \\[1mm]
& \multicolumn{2}{r}{$(n=5000)$} & Drift $2^{\ast}$ & & 0 & {\bf 943} & 28 & 0 \\[1mm]
& & & Drift 3 & & 0 & 28 & 1 & 0 \\[1mm]
& 50 & 0.005 & Drift 1 & & 0 & 0 & 0 & 0 \\[1mm]
& \multicolumn{2}{r}{$(n=10000)$} & Drift $2^{\ast}$ & & 0 & {\bf 957} & 16 & 0 \\[1mm]
& & & Drift 3 & & 0 & 26 & 1 & 0 \\[1mm] \hline\hline
%%%%%
$\mathrm{GQBIC}^{\sharp}$ & $T_{n}$ & $h_{n}$ & & & Scale 1 & Scale $2^{\ast}$ & Scale 3 & Scale 4 \\ \hline
& 10 & 0.01 & Drift 1 & & 0 & 0 & 0 & 0 \\[1mm]
& \multicolumn{2}{r}{$(n=1000)$} & Drift $2^{\ast}$ & & 0 & {\bf 264} & 64 & 477 \\[1mm]
& & & Drift 3 & & 5 & 37 & 45 & 108 \\[1mm]
& 10 & 0.005 & Drift 1 & & 0 & 0 & 0 & 1 \\[1mm]
& \multicolumn{2}{r}{$(n=2000)$} & Drift $2^{\ast}$ & & 0 & {\bf 181} & 48 & 565 \\[1mm]
& & & Drift 3 & & 2 & 35 & 29 & 139 \\[1mm]
& 50 & 0.01 & Drift 1 & & 0 & 0 & 0 & 0 \\[1mm]
& \multicolumn{2}{r}{$(n=5000)$} & Drift $2^{\ast}$ & & 0 & {\bf 210} & 136 & 619 \\[1mm]
& & & Drift 3 & & 0 & 7 & 1 & 27 \\[1mm]
& 50 & 0.005 & Drift 1 & & 0 & 0 & 0 & 0 \\[1mm]
& \multicolumn{2}{r}{$(n=10000)$} & Drift $2^{\ast}$ & & 0 & {\bf 166} & 96 & 708 \\[1mm]
& & & Drift 3 & & 0 & 4 & 1 & 25 \\[1mm] \hline
\end{tabular}
}
\label{simu:tab3}
\end{center}
\end{table}

%%%%%
%%%%%

%%%%%%%

%\appendix

%%%%%
%%%%%
\section{Proofs}
\label{hm:sec_proofs}

%%%%%
\subsection{Proof of Theorem \ref{hm:thm_GQMLE}}
\label{hm:sec_proof_stepwise.GQMLE}

The proofs are essentially the same as in those of \cite[Theorem 2.7]{Mas13}, \cite[Theorem 3.4]{MasUeh17-2}, making use of the general machinery \cite{Yos11}. Hence we only mention the formal difference, omitting the further details:
The only difference to be mentioned is that the proofs are based on the two-stage procedure for $M$-estimators as in \cite[Section 6]{Yos11}, where the first-stage random field is
\begin{equation}
u_\gam \mapsto \log \mbbz_{1,n}(u_\gam):= h\left( \mbbh_{1,n}(\gam_0 + T_n^{-1/2}u_\gam) - \mbbh_{1,n}(\gam_0)\right),
\nonumber
\end{equation}
and the second-stage one (depending on $\ges$) is the rescaled
\begin{equation}
u_\al \mapsto \log \mbbz_{2,n}(u_\al):=\mbbh_{2,n}(\al_0 + T_n^{-1/2}u_\al) - \mbbh_{2,n}(\al_0).
\nonumber
\end{equation}
Note the resolution in handling the scale coefficient in the first stage is corrected by the multiplicative factor ``$h$''.
The Studentization \eqref{hm:GQMLE.asn} can be verified exactly as in \cite[Corollary 2.8]{Mas13}.

%%%%%
\subsection{Proof of Proposition \ref{hm:prop+1}}
%Concerned with \eqref{hm:gqaic_thm-1}, we recall the expression \eqref{hm:W4-hat}.
Let us recall the expression \eqref{hm:W4-hat} for $\hat{W}_{\gam,n}=(\hat{W}_{\gam,n}^{(qr)})$.
Under the integrability conditions, the sequence $(\hat{W}_{\gam,n})_n$ is $L^q(\pr)$-bounded for any $q>0$.
To see this, %we may and do suppose that $p_\gam=1$ without loss of generality.
let $\chi_j := \D_j X - ha_{j-1}(\al_0)$, and write $O^\ast_p(1)$ for a random sequence $(\zeta_n)_n$ such that $\sup_n\E(|\zeta_n|^q)<\infty$ for any $q>0$. Write $\E^{j-1}$ for the expectation conditional on $\mcf_{t_{j-1}}$, where $(\mcf_t)$ denotes the underlying filtration to which all the stochastic processes are adapted.
Then, by compensation and Burkholder's inequality, we have
\begin{align}
|\hat{W}_{\gam,n}| &\lesssim \frac{1}{T_n}\sumj (1+|X_{t_{j-1}}|)^C |\hat{\chi}_j|^4 \nn\\
&\lesssim \frac{1}{T_n}\sumj (1+|X_{t_{j-1}}|)^C |\chi_j|^4 + O^\ast_p(1) \nn\\
&\lesssim \frac{1}{T_n}\sumj (1+|X_{t_{j-1}}|)^C E^{j-1}\left[|\chi_j|^4\right] + O^\ast_p(T_n^{-1/2})+O^\ast_p(1) \nn\\
&\lesssim \frac{1}{n}\sumj (1+|X_{t_{j-1}}|)^C + O^\ast_p(T_n^{-1/2})+O^\ast_p(1) = O^\ast_p(1).
\nonumber
\end{align}
%Therefore, the uniform integrability of $(|\hat{\Gam}_{\gam,n}^{-1}|)_n$ is enough to conclude \eqref{hm:gqaic_thm-1}. 
Since $|\hat{\Gam}_{\gam,n}^{-1}|$ is bounded by a universal-constant multiple of $\lam_{\min}^{-1}(\hat{\Gam}_{\gam,n})$,
we can apply H\"{o}lder' inequality to ensure that \eqref{hm:v0.6_add1} is sufficient for \eqref{hm:gqaic_thm-1}.

%%%%%
\subsection{Proof of Lemma \ref{hm:lem_imb}}
To begin with, let $k\ll n$ be a positive integer not depending on $n$, and let $m:=\lfloor n/k\rfloor$; without loss of generality, we set $k\le n/2$.
Then, for $u\in\mbbr^{p_{\gamma}}$,
\begin{align}
\hat{\Gam}_{\gam,n} [u^{\otimes 2}]
&= \frac{1}{2n} \sumj \trace\left\{\left(\hat{S}^{-1}_{j-1}(\p_{\gam}\hat{S}_{j-1})[u]\right)^2\right\}
\nn\\
&\ge \frac{1}{2n} \sumj \frac1d \left\{\trace\left((\hat{S}^{-1}_{j-1}\p_{\gam}\hat{S}_{j-1})[u]\right)\right\}^{2} \nn\\
&= \frac{1}{n} \sumj \frac{1}{2d}\left\{ 
%\big( (\p_\gam \log|S_{j-1}|)(\ges) \big)
\zeta_{j-1}(\ges)
[u]\right\}^{2} \nn\\
&\gtrsim \frac{1}{mk} \sum_{j=1}^{mk} \left\{ 
%\big( (\p_\gam \log|S_{j-1}|)(\ges) \big)
\zeta_{j-1}(\ges)[u] \right\}^{2} 
=: \frac{1}{k} \cdot \frac{1}{m} \sum_{i=1}^{m} V_{i}(\ges)[u^{\otimes 2}],
\label{hm:imb-1}
\end{align}
where $V_{i}(\gam) := \sum_{j=(i-1)k+1}^{ik} \zeta_{j-1}(\gam)^{\otimes 2}$;
% with $\xi_{j-1}(\gam):=\big( (\p_\gam \log|S_{j-1}|)(\gam) \big)$;
the first inequality is due to the Cauchy-Schwarz inequality: $\trace(A^2)\ge d^{-1}\trace(A)^2$ %$\mathrm{rank}(A)\trace(A^2)\ge\trace(A)^2$ 
for any real square matrix $A$ with real eigenvalues.
Observe that by \eqref{hm:imb-1} and the Jensen inequality,
\begin{align}
\lam_{\min}^{-q}(\hat{\Gam}_{\gam,n}) 
&\le \left( \inf_{u\in\mbbs}\hat{\Gam}_{\gam,n}[u^{\otimes 2}] \right)^{-q} \nn\\
%&\le C_{d,k,q} \left( \frac1m \sum_{i=1}^{m} \inf_{|u|=1} Q_{i}(\ges)[u^{\otimes 2}] \right)^{-q} \nn\\
&\lesssim \left( \frac1m \sum_{i=1}^{m} \inf_{u\in\mbbs} \inf_\gam V_{i}(\gam)[u^{\otimes 2}] \right)^{-q} \nn\\
&\lesssim \frac1m \sum_{i=1}^{m} \left( \inf_{u\in\mbbs} \inf_\gam V_{i}(\gam)[u^{\otimes 2}] \right)^{-q}.
\nonumber
\end{align}
%where and in what follows, $C_{d,k,q}$ denotes a universal constant depending only on $(d,k,q)$.
It suffices for \eqref{hm:v0.6_add1} to have
\begin{equation}
\sup_{n}\sup_{i\in\mbbn}\E\left[ \left( \inf_{u\in\mbbs} \inf_\gam V_{i}(\gam)[u^{\otimes 2}] \right)^{-q} \right] < \infty.
\label{hm:imb-2}
\end{equation}

Fix a constant $r>0$ in the sequel. The expectation in \eqref{hm:imb-2} equals
\begin{align}
& \int_0^\infty \pr\left[ \left( \inf_{u\in\mbbs} \inf_\gam V_{i}(\gam)[u^{\otimes 2}] \right)^{-q} > s\right] ds \nn\\
&\le 1 + \int_1^\infty \pr\left[ \inf_{u\in\mbbs} \inf_\gam \sum_{j=(i-1)k+1}^{ik} \big(\zeta_{j-1}(\gam)[u]\big)^2
< s^{-1/q} \right] ds \nn\\
&\le 1 + \int_1^\infty \left(
\pr\left[ \sum_{j=(i-1)k+1}^{ik} \sup_\gam |\zeta_{j-1}(\gam)|^2 \ge s^{r/q} \right] 
\right.
\nn\\
& \qquad \left. +\pr\left[ \inf_{u\in\mbbs} \inf_\gam \sum_{j=(i-1)k+1}^{ik} \big(\zeta_{j-1}(\gam)[u]\big)^2< s^{-1/q},~ 
\sum_{j=(i-1)k+1}^{ik} \sup_\gam |\zeta_{j-1}(\gam)|^2 \le s^{r/q}\right] 
\right) ds 
\nn\\
&=: 1 + \int_1^\infty \left( I'_{i,k,q}(s) + I''_{i,k,q}(s)\right)ds.
\nn
%\\
%&\le ... \nn\\
%&\le 1 + \int_1^\infty \pr\left[ \bigcap_{j=(i-1)k+1}^{ik} \left\{ \inf_{u\in\mbbs} \inf_\gam \big(\zeta_{j-1}(\gam)[u]\big)^2 < s^{-1/q} \right\} \right] ds \nn\\
%&\le 1 + \int_1^\infty \bigg( 
%\pr\left[ \bigcap_{j=(i-1)k+1}^{ik} \left\{ \inf_{u\in\mbbs} \inf_\gam \big(\zeta_{j-1}(\gam)[u]\big)^2 < s^{-1/q} \right\} \right] 
%\bigg)ds 
%\nn
\end{align}
Since $k$ is fixed and
\begin{equation}
I'_{i,k,q}(s) \lesssim s^{-2} \sup_{t}\E\left[ \sup_\gam 
\big| \zeta(X_t,\gam) \big|^{4q/r}
\right] \lesssim s^{-2} \left(1 + \sup_t \E[|X_t|^C]\right) \lesssim s^{-2}
\nonumber
\end{equation}
whatever $r>0$ is under the present assumptions, it remains to be shown that
\begin{equation}
\int_1^\infty I''_{i,k,q}(s) ds \lesssim 1.
\label{hm:imb-5}
\end{equation}

First, to handle the infimum for $u$, we will apply Lemma \ref{hm:lem_FinWei.A3}.
Let $r':=(r+1)/2$.
Taking $\del = s^{-r'/q}$ in Lemma \ref{hm:lem_FinWei.A3}, we can pick some elements $u'_1,\dots,u'_{D(s)}$ with the integer
\begin{equation}
D(s) = O\big(s^{r' (p_\gam -1)/q}\big) \qquad s\uparrow \infty,
\nonumber
\end{equation}
which is constant for $p_\gam=1$.
Write $\mbbs_l := \{u\in\mbbs:\,|u-u'_l|<s^{-r'/q}\}$.
Since $\inf_{u\in\mbbs_l} \inf_\gam \big| \zeta_{j-1}(\gam)[u]\big| \ge \inf_\gam \big| \zeta_{j-1}(\gam)[u'_l]\big| - s^{-r'/q} \sup_\gam |\zeta_{j-1}(\gam)|$ for each $u\in\mbbs_l$, we have
\begin{align}
I''_{i,k,q}(s) 
%&\le \sum_{l=1}^{D(s)}
%\pr\left[ \inf_{u\in\mbbs_l} \inf_\gam \sum_{j=(i-1)k+1}^{ik} \big(\zeta_{j-1}(\gam)[u]\big)^2 < s^{-1/q},~ 
%\sum_{j=(i-1)k+1}^{ik} \sup_\gam |\zeta_{j-1}(\gam)|^2 \le s^{r/q}\right] \nn\\ 
&\le \sum_{l=1}^{D(s)}
\pr\left[ \bigcap_{j=(i-1)k+1}^{ik} 
\left\{ \inf_{u\in\mbbs_l} \inf_\gam \big| \zeta_{j-1}(\gam)[u]\big| < s^{-1/(2q)},~ 
\sup_\gam |\zeta_{j-1}(\gam)| \le s^{r/(2q)}\right\}\right] \nn\\
&\le \sum_{l=1}^{D(s)}
\pr\left[ \bigcap_{j=(i-1)k+1}^{ik} 
\left\{ \inf_\gam \big| \zeta_{j-1}(\gam)[u'_l]\big| < 2s^{-1/(2q)},~ 
\sup_\gam |\zeta_{j-1}(\gam)| \le s^{r/(2q)}\right\}\right]
\nn\\
&\le \sum_{l=1}^{D(s)}
\pr\left[ \bigcap_{j=(i-1)k+1}^{ik} \left\{ \inf_\gam \big| \zeta_{j-1}(\gam)[u'_l]\big| < 2s^{-1/(2q)}\right\}\right].
\label{hm:imb-3}
\end{align}

Next, we get rid of the infimum with respect to $\gam$.
Since $\overline{\Theta}_\gam\subset\mbbr^{p_\gam}$ is compact, we can cover it by finitely many hypercubes $\mcu_1,\dots,\mcu_{H(s)}$, each with side length $s^{-1/(2q)}$ and the number
\begin{equation}
H(s) = O\big(s^{p_\gam/(2q)}\big) \qquad s\uparrow\infty.
\nonumber
\end{equation}
Pick elements $\gam'_b \in \mcu_b$ arbitrarily ($b=1,\dots,H(s)$).
We also have $\sup_\gam|\p_\gam \zeta_{j-1}(\gam)| \lesssim 1 + |X_{t_{j-1}}|^C$.
%\begin{equation}
%\forall K>0,\quad \sup_{n} \sup_{j\le n}\E\left[\sup_\gam \big|\p_\gam \zeta_{j-1}(\gam)\big|^K\right] < \infty.
%\nonumber
%\end{equation}
With these observations, for some nonnegative function $F$ such that $F(x) \lesssim 1+|x|^C$ we can continue \eqref{hm:imb-3} as follows:
\begin{align}
I''_{i,k,q}(s) 
&\le \sum_{l=1}^{D(s)} \sum_{b=1}^{H(s)}
\pr\left[ \bigcap_{j=(i-1)k+1}^{ik} \left\{ \inf_{\gam\in\mcu_{b}} \big| \zeta_{j-1}(\gam)[u'_l]\big| < 2s^{-1/(2q)}\right\}\right]
\nn\\
&\le \sum_{l=1}^{D(s)} \sum_{b=1}^{H(s)} \pr\left[ \bigcap_{j=(i-1)k+1}^{ik} \left\{ \big| \zeta_{j-1}(\gam'_b)[u'_l]\big| \le s^{-1/(2q)} F_{j-1}\right\}\right] \nn\\
&\le 
\sum_{l=1}^{D(s)} \sum_{b=1}^{H(s)} \pr\left[ \bigcap_{j=(i-1)k+1}^{ik} \left\{ \big| \zeta_{j-1}(\gam'_b)[u'_l]\big| \le s^{-1/(4q)}\right\}\right]
+\sum_{l=1}^{D(s)} \sum_{b=1}^{H(s)} \pr\left[ F_{j-1} \ge s^{1/(4q)}\right].
\label{hm:imb-4}
\end{align}
Whatever $r>0$ is, the second term in \eqref{hm:imb-4} can be bounded by $C s^{-2}$ through the Markov inequality with sufficiently high-order moments.
Finally, by iterative conditioning, the first term in \eqref{hm:imb-4} is a.s. bounded by $C \sum_{l=1}^{D(s)} \sum_{b=1}^{H(s)} (s^{-1/(4q)})^{k \rho}$ with $\rho>0$ of \eqref{hm:lem_imb_A*}; here is the only place where we used the condition \eqref{hm:lem_imb_A*}.
Given $p_\gam$, $q=1+\del$, $r>0$, and $\rho>0$, we can pick a sufficiently large $k\in\mbbn$ to ensure that $C \sum_{l=1}^{D(s)} \sum_{b=1}^{H(s)} (s^{-1/(4q)})^{k \rho} \lesssim s^{-2}$, thus concluding \eqref{hm:imb-5}. The proof is complete.

%%%%%
\subsection{Proof of Corollary \ref{hm:lem_imb-cor1}}
By the iterative conditioning with taking $k$ large enough, we can bound the first term in \eqref{hm:imb-4} from above by $\sum_{l=1}^{D(s)} \sum_{b=1}^{H(s)} s^{-\rho k/(2q)} \lesssim s^{-2}$. Hence \eqref{hm:imb-5}.

%%%%%
\subsection{Proof of Corollary \ref{hm:lem_imb-cor2}}
\eqref{hm:imb-2} readily follows from \eqref{hm:imb-6}; in this case we may set $k=1$ in the proof of Lemma \ref{hm:lem_imb}.

%%%%%
\subsection{Proof of Theorem \ref{hm:thm_gqbic1}}
\label{hm:sec_thm_gqbic1_proof}

Roughly, \eqref{hm:gqbic_expa1_EU} will be proved by expanding $\mbbh_{1,n}(\gam)$ around $\ges$ with vicinity size of order $n^{-1/2}$, while \eqref{hm:gqbic_expa1} around $\gam_0$ with vicinity size of order $T_n^{-1/2}$.
%The former is the well-known classic route (see e.g. \cite{CavNea99}), while not appropriate for consistent model selection.

%%%
\subsubsection{Proof of \eqref{hm:gqbic_expa1_EU}}

Our proof is achieved in an analogous way to the derivation of the Bernstein-von Mises theorem given in \cite[Theorem 1]{JasKamMas19}, which dealt with a more complicated two-step non-ergodic setting.
\footnote{
Therefore, it should be remarked that without any essential change the proof below could be easily extended to the non-ergodic framework, where the matrices $\Sig(\tz)$ and $\Gam(\tz)$ are random.
}
%We will give a sketch of the proof for the sake of completeness.

By the change of variable $\gam=\ges+n^{-1/2}u$, we have
\begin{align}
\mfF_{1,n}(1)
%&= -\frac{1}{T_n}\log\left(\int_{\Theta_\gam}\exp\{h \mbbh_{1,n}(\gam)\} \pi_1(\gam)d\gam\right) \nn\\
%&= -\frac{1}{n} \mbbh_{1,n}(\ges) + \frac{p_\gam}{2n}\log n - \frac{1}{n}\log\left(\int_{\hat{U}_{1,n}} \hat{\mbbz}_{1,n}(u) \pi_1(\ges+n^{-1/2}u)du\right) \nn\\
&= -\frac{1}{n}\mbbh_{1,n}(\ges) + \frac{p_\gam}{2n}\log n - \frac{1}{n}\log \hat{\mbbz}^\ast_{1,n},
\nonumber
\end{align}
where $\hat{\mbbz}^\ast_{1,n} = \int_{\hat{U}_{1,n}} \hat{\mbbz}_{1,n}(u) \pi_1(\ges+n^{-1/2}u)du$ with $\hat{U}_{1,n}:=\{v\in\mbbr^{p_{\gam}}:\, \ges+n^{-1/2}v \in \Theta_\gam\}$ and $\hat{\mbbz}_{1,n}(u) := \exp\{\mbbh_{1,n}(\ges+n^{-1/2}u)-\mbbh_{1,n}(\ges)\}$.
It suffices to show that $\log \hat{\mbbz}^\ast_{1,n}=O_p(1)$.

We need some notation and preliminary remarks.
Let $\mbby_{1,n}(\gam):=n^{-1}\{ \mbbh_{1,n}(\gam)-\mbbh_{1,n}(\gam_0)\}$ and 
$\mbby_{1,0}(\gam):=-(1/2)\int\{\trace\left(S(x,\gam)^{-1}S(x,\gam_0) - I_d\right) + \log(|S(x,\gam)|/|S(x,\gam_0)|)\}\pi(dx)$. Then, Assumption \ref{hm:ass_iden} ensures that we can find a constant $\chi_\gam>0$ such that 
%$\mbby_{1,0}(\gam) \le -\chi_\gam|\gam-\gam_0|^2$ for every $\gam$, namely 
\begin{equation}
\sup_{\gam:\,|\gam-\gam_0|\ge \del}\mbby_{1,0}(\gam) \le -\chi_\gam \del^2
\label{hm:gqbic_proof-2}
\end{equation}
for every $\del>0$.
By \cite[Lemma 4.3]{Mas13} we know that $\sup_{\theta}\sqrt{T_n}|\mbby_{1,n}(\gam)-\mbby_{1,0}(\gam)|=O_p(1)$.
Moreover, note that $\D_{1,n}(\ges)=h\,T_n^{-1/2}\p_\gam \mbbh_{1,n}(\ges)=0$ if $\ges\in\Theta_\gam$, that $\Gam_{\gam,n}(\tz) \cip \Gam_\gam(\gam_0)$ (see \eqref{hm:Gam1->}), and that $n^{-1}\sup_\gam |\p_\theta^3\mbbh_{1,n}(\gam)|=O_p(1)$.
We pick a constant (recall \eqref{hm:h-rate-add})
\begin{equation}
0 < c_0 < \frac{c_1}{4}.
%0 < c_0 < \left(\frac12 \wedge \frac{c_1}{4}\right) = \frac{c_1}{4}.
\label{hm:gqbic_proof-3}
\end{equation}
%Then we have $c_0<1/2$.
Put $\ep_n=n^{-c_0}$ in the sequel.

We introduce the following auxiliary event for constants $M,\lam>0$ (recall $\hat{u}_{\gam,n} := \sqrt{T_n}(\ges-\gam_0)$):
\begin{align}
G_{1,n}(M,\lam) 
&:= \bigg\{
\ges\in\Theta_\gam, \quad \left| \Gam_{\gam,n}(\ges) - \Gam_\gam(\ges) \right| \le \lam,
\quad \lam_{\min}(\Gam_{\gam}(\ges)) \ge 4\lam,
\nn\\
&{}\qquad 
\big|\sqrt{T_n}\sup_{\theta}(\mbby_{1,n}(\gamma)-\mbby_{1}(\gam))\big| \vee |\hat{u}_{\gam,n}| < M,
%\nn\\
%&{}
\quad \frac{1}{n} \sup_\gam \left|\p_\theta^3\mbbh_{1,n}(\gam)\right| \le \frac{3\lam}{\ep_n}
\bigg\}.
\nonumber
\end{align}
Fix any $\ep>0$.
Then, we can find a pair $(M_{1},\lam_{1})$ and an $N_{1}\in\mbbn$ such that
\begin{equation}
\sup_{n\ge N_{1}}\pr\left\{G_{1,n}(M,\lam)^{c}\right\}<\ep
\nonumber
\end{equation}
holds for every $M\ge M_{1}$ and $\lam\in(0,\lam_{1}]$.
Therefore, to deduce $\log \hat{\mbbz}^\ast_{1,n}=O_p(1)$, we may and do focus on the event $G_{1,n}(M,\lam)$ with $M=M(\ep)$ and $\lam=\lam(\ep)$ being sufficiently large and small, respectively.

Let $A_{1,n}:=\{u\in\mbbr^{p_\gam}:\, |u| \le \ep_n\sqrt{n}\}$.
For $u\in A_{1,n}$ and on $G_{1,n}(M,\lam)$, we apply the third-order Taylor expansion to conclude that
\begin{equation}
\log \hat{\mbbz}_{1,n}(u) \le -\lam |u|^2,
\nonumber
\end{equation}
by noting that, for some random point $\tilde{\gam}_n$ on the segment joining $\ges$ and $\gam_0$, we have
$\log \hat{\mbbz}_{1,n}(u) = n^{-1/2} \p_\gam\mbbh_{1,n}(\ges)[u] -(1/2)\{
\Gam_{\gam}(\ges) + (\Gam_{\gam,n}(\ges)-\Gam_{\gam}(\ges)) - (3n)^{-1}\p_\gam^3\mbbh_{1,n}(\tilde{\gam}_n)[n^{-1/2}u]\}[u^{\otimes 2}]$.
%\begin{align}
%\log \hat{\mbbz}_{1,n}(u) &= \frac{1}{\sqrt{n}} \p_\gam\mbbh_{1,n}(\ges)[u] \nn\\
%&{}\qquad -\frac12\left\{
%\Gam_{\gam}(\ges) + (\Gam_{\gam,n}(\ges)-\Gam_{\gam,n}(\ges)) - \frac{1}{3n}\p_\gam^3\mbbh_{1,n}(\tilde{\gam}_n)\left[\frac{u}{\sqrt{n}}\right]
%\right\}[u^{\otimes 2}].
%\nonumber
%\end{align}
This entails that $\sup_n \sup_{u\in A_{1,n}\cap \hat{U}_{1,n}} \hat{\mbbz}_{1,n}(u) \pi_1(\ges+n^{-1/2}u) \lesssim \exp(-|u|^2)$, 
%\begin{equation}
%\sup_n \sup_{u\in A_{1,n}\cap \hat{U}_{1,n}} \hat{\mbbz}_{1,n}(u) \pi_1(\ges+n^{-1/2}u) \lesssim \exp(-|u|^2),
%\nonumber
%\end{equation}
followed by
\begin{align}
& \left|\int_{A_{1,n}\cap\hat{U}_{1,n}} \hat{\mbbz}_{1,n}(u) \{\pi_1(\ges+n^{-1/2}u) - \pi_1(\ges)\}du\right| \nn\\
&\le \sup_{u\in A_{1,n}\cap\hat{U}_{1,n}}\left|\pi_1(\ges+n^{-1/2}u) - \pi_1(\ges)\right|
\int_{A_{1,n}\cap\hat{U}_{1,n}}\hat{\mbbz}_{1,n}(u)du \nn\\
&\lesssim \sup_{|v_n|\le \ep_n}|\pi_1(\ges+v_n) - \pi_1(\ges)| \cip 0.
\nonumber
\end{align}
For $\mbbz^0_{1}(u) := \exp\{-(1/2)\Gam_{\gam}(\gam_0)[u^{\otimes 2}]\}$, we can deduce that $\int_{A_{1,n}\cap\hat{U}_{1,n}} \hat{\mbbz}_{1,n}(u)du = \int_{A_{1,n}\cap\hat{U}_{1,n}} \mbbz^0_{1}(u)du + o_p(1)$ (on $G_{1,n}(M,\lam)$) by using the subsequence argument in much the same way as in \cite{JasKamMas19}.
Since $\int_{\mbbr^{p_{\gamma}}} \mbbz^0_{1}(u)du=(2\pi)^{p_{\gam}/2}|\Gam_\gam(\gam_0)|^{-1/2}+o(1)$, we conclude that 
$\log\{\int_{A_{1,n}\cap\hat{U}_{1,n}} \hat{\mbbz}_{1,n}(u) \pi_1(\ges+n^{-1/2}u)du\} = \log\pi(\gam_0) + (p_\gam/2)\log(2\pi) - (1/2)\log|\Gam_\gam(\gam_0)|+ o_p(1) = O_p(1)$.

We are left to proving $\int_{A_{1,n}^c\cap\hat{U}_{1,n}} \hat{\mbbz}_{1,n}(u) \pi_1(\ges+n^{-1/2}u)du = o_p(1)$; since $\pi_1$ is bounded, it suffices to show that $\int_{A_{1,n}^c\cap\hat{U}_{1,n}} \hat{\mbbz}_{1,n}(u)du = o_p(1)$.
We have on $G_{1,n}(M,\lam)$,
\begin{align}
\sup_{u\in A_{1,n}^c \cap\hat{U}_{1,n}}\log\hat{\mbbz}_{1,n}(u)
&\le n\sup_{u\in A_{1,n}^c\cap\hat{U}_{1,n}}\left\{
\mbby_{1,n}\left(\ges+\frac{u}{\sqrt{n}}\right) - \mbby_{1,0}\left(\ges+\frac{u}{\sqrt{n}}\right) \right\} 
\nn\\
&{}\qquad 
+ \mbby_{1,0}\left(\ges+\frac{u}{\sqrt{n}}\right) \nn\\
&\le n\left(
\sup_\gam |\mbby_{1,n}(\gam)-\mbby_{1,0}(\gam)| 
+ \sup_{u\in A_{1,n}^c\cap\hat{U}_{1,n}}\mbby_{1,0}\left(\ges+\frac{u}{\sqrt{n}}\right)
\right).
%\nn\\
%&= n\left( o_p(1)
%+ \sup_{u\in A_{1,n}^c\cap\hat{U}_{1,n}}\mbby_{1,0}\left(\ges+\frac{u}{\sqrt{n}}\right)
%\right).
\nn\\
&\le n\left( M n^{-c_1/2} + \sup_{u\in A_{1,n}^c\cap\hat{U}_{1,n}}\mbby_{1,0}\left(\ges+\frac{u}{\sqrt{n}}\right)
\right)
\label{hm:gqbic_proof-1}
\end{align}
Observe that
\begin{align}
\inf_{u\in A_{1,n}^c}\left|\left(\ges+\frac{u}{\sqrt{n}}\right) - \gam_0\right|
&\ge \inf_{u\in A_{1,n}^c}\frac{|u|}{\sqrt{n}} - \frac{|\hat{u}_{\gam,n}|}{\sqrt{T_n}}
\ge n^{-c_0}\left(1-M n^{c_0-c_1/2}\right) \ge \frac12 n^{-c_0}
\nonumber
\end{align}
for every $n$ large enough.
Recalling \eqref{hm:gqbic_proof-2} and \eqref{hm:gqbic_proof-3}, we can continue the estimate \eqref{hm:gqbic_proof-1} as follows:
\begin{align}
\sup_{u\in A_{1,n}^c \cap\hat{U}_{1,n}}\log\hat{\mbbz}_{1,n}(u)
&\lesssim -n^{1-2c_0} \left(1- n^{2c_0-c_1/2}\right) \lesssim -n^{1-2c_0} \downarrow -\infty.
\nonumber
\end{align}
Thus $\int_{A_{1,n}^c\cap\hat{U}_{1,n}} \hat{\mbbz}_{1,n}(u)du \lesssim \exp(-C n^{1-2c_0})\int_{|u|\le C\sqrt{n}} du \lesssim n^{p_\gam/2}\exp(-C n^{1-2c_0}) \to 0$, concluding that $\int_{A_{1,n}^c\cap\hat{U}_{1,n}} \hat{\mbbz}_{1,n}(u)du = o_p(1)$.
The proof is complete.

%%%
\subsubsection{Proof of \eqref{hm:gqbic_expa1}}

The proof is similar to the proof of \eqref{hm:gqbic_expa1_EU} and much closer to that of \cite[Theorem 1]{JasKamMas19}.

Let $\mbbz_{1,n}(u):=\exp\{h(\mbbh_{1,n}(\gam_0+T_n^{-1/2}u) - \mbbh_{1,n}(\gam_0))\}$. 
The change of variable $\gam=\gam_0+T_n^{-1/2}u$ yields
\begin{align}
\mfF_{1,n}(h)
%&= -\frac{1}{T_n}\log\left(\int_{\Theta_\gam}\exp\{h \mbbh_{1,n}(\gam)\} \pi_1(\gam)d\gam\right) \nn\\
&= -\frac{1}{T_n} h\mbbh_{1,n}(\gam_0) + \frac{p_\gam}{2T_n}\log T_n - \frac{1}{T_n}\log\left(\int_{U_{1,n}} \mbbz_{1,n}(u) \pi_1(\gam_0+T_N^{-1/2}u)du\right) \nn\\
&=: -\frac{1}{T_n}h\mbbh_{1,n}(\gam_0) + \frac{p_\gam}{2T_n}\log T_n - \frac{1}{T_n}\log \overline{\mbbz}_{1,n},
\nonumber
\end{align}
where $U_{1,n}:=\{v\in\mbbr^{p_{\gam}}:\, \gam_0+T_n^{-1/2}v \in \Theta_\gam\}$; a relevant form already appeared in Remark \ref{hm:rem_mixed.rates.M}.
By Theorem \ref{hm:thm_GQMLE} (and its proof), it is easily seen that $h\mbbh_{1,n}(\gam_0) = h\mbbh_{1,n}(\ges) + O_p(1)$, so that it suffices to show $\log \overline{\mbbz}_{1,n}=O_p(1)$.

Recall that $\D_{1,n}(\gam_0):=h\,T_n^{-1/2}\p_\gam \mbbh_{1,n}(\gam_0)=O_p(1)$.
In the present case, the auxiliary event is given as follows (we keep using $\ep_n=n^{-c_0}$ such that \eqref{hm:gqbic_proof-3} holds): for constants $\lam\in(0,\lam_{\min}(\Gam_{\gam}(\gam_0))/4)$ and $M>0$,
\begin{align}
G_{1,n}(M,\lam) 
&:= \bigg\{
|\Delta_{1,n}(\gam_0)| \le M,
\quad \left| \Gam_{\gam,n}(\tz) - \Gam_\gam(\gam_0) \right| < \lam,
\nn\\
&{}\qquad 
\sqrt{T_n}\sup_{\theta}|\mbby_{1,n}(\gamma)-\mbby_{1}(\gam)|< M, \quad 
\frac{1}{n} \sup_\gam \left|\p_\theta^3\mbbh_{1,n}(\gam)\right| \le \frac{3\lam}{\ep_n}
 \bigg\}.
\nonumber
\end{align}
With this $G_{1,n}(M,\lam)$, the remaining arguments are almost identical to those of \cite{JasKamMas19}, hence omitted.

%%%%%
\subsection{Proof of Theorem \ref{se:thm_sele.prob}}

Under the assumptions, the coefficients $c_{m_{1}^{\ast}}$ and $a_{m_{2}^{\ast}}$ are correctly specified, $\gam_{m_{1}^{\ast}}^{\ast}=\gam_{m_{1}^{\ast},0}$, and $\al_{m_{2}^{\ast}}^{\ast}=\al_{m_{2}^{\ast},0}$.

\subsubsection{Proof of 1} \label{se:prf_thm_sele.prob1}
(i) In this case, $c_{m_{1}}$ is correctly specified, and $\gam_{m_{1}}^{\ast}=\gam_{m_{1},0}$.
Furthermore, the equation $\mbbh_{1,0}^{(m_{1}^{\ast})}(\gam_{m_{1}^{\ast},0})=\mbbh_{1,0}^{(m_{1})}(\gam_{m_{1},0})$ holds.
Define the map $f_{1}:\Theta_{\gamma_{m_{1}^{\ast}}}\to\Theta_{\gamma_{m_{1}}}$ by $f_{1}(\gamma_{m_{1}^{\ast}})=F_{1}\gamma_{m_{1}^{\ast}}+c_{1}$, where $F_{1}$ and $c_{1}$ satisfy the equation $\mbbh_{1,n}^{(m_{1}^{\ast})}(\gamma_{m_{1}^{\ast}})=\mbbh_{1,n}^{(m_{1})}\big(f_{1}(\gamma_{m_{1}^{\ast}})\big)$ for any $\gamma_{m_{1}^{\ast}}\in\Theta_{\gamma_{m_{1}^{\ast}}}$.
When $f_{1}(\gamma_{m_{1}^{\ast},0})\neq\gamma_{m_{1},0}$, $\mbbh_{1,0}^{(m_{1}^{\ast})}(\gamma_{m_{1}^{\ast},0})=\mbbh_{1,0}^{(m_{1})}\big(f_{1}(\gamma_{m_{1}^{\ast},0})\big)<\mbbh_{1,0}^{(m_{1})}(\gamma_{m_{1},0})$.
Hence, we have $f_{1}(\gamma_{m_{1}^{\ast},0})=\gamma_{m_{1},0}$.

Using the fact that $P\{\p_{\gam}\mbbh_{1,n}^{(m_{1})}(\hat{\gam}_{m_{1},n})=0\}\to1$ and Taylor expansion of $\mbbh_{1,n}^{(m_{1})}$,
\begin{align*}
\mbbh_{1,n}^{(m_{1}^{\ast})}(\hat{\gamma}_{m_{1}^{\ast},n})
&=\mbbh_{1,n}^{(m_{1})}\left(f_{1}(\hat{\gamma}_{m_{1}^{\ast},n})\right)\\
&=\mbbh_{1,n}^{(m_{1})}(\hat{\gamma}_{m_{1},n})-\frac{1}{2}\left(-\p_{\gamma_{m_{1}}}^{2}\mbbh_{1,n}^{(m_{1})}(\tilde{\gamma}_{m_{1},n})\right)\left[\left\{\hat{\gamma}_{m_{1},n}-f_{1}(\hat{\gamma}_{m_{1}^{\ast},n})\right\}^{\otimes2}\right]
\end{align*}
where $\tilde{\gamma}_{m_{1},n}=\hat{\gamma}_{m_{1},n}-\xi_{1}\big\{f_{1}(\hat{\gamma}_{m_{1}^{\ast},n})-\hat{\gamma}_{m_{1},n}\big\}$ for some $0<\xi_{1}<1$ and $\tilde{\gamma}_{m_{1},n}\cip\gamma_{m_{1},0}$ as $n\to\infty$.
Therefore, the difference between $\gqaic_{1,n}^{(m_{1}^{\ast})}$ and $\gqaic_{1,n}^{(m_{1})}$ is given by
\begin{align*}
\gqaic_{1,n}^{(m_{1}^{\ast})}-\gqaic_{1,n}^{(m_{1})}
&=\left(-\p_{\gamma_{m_{1}}}^{2}\mbbh_{1,n}^{(m_{1})}(\tilde{\gamma}_{m_{1},n})\right)\left[\left\{\hat{\gamma}_{m_{1},n}-f_{1}(\hat{\gamma}_{m_{1}^{\ast},n})\right\}^{\otimes2}\right] \\
&\qquad+\frac{2}{h}\trace\left(\hat{\Gam}_{\gam_{m_{1}^{\ast}},n}^{-1}\hat{W}_{\gam_{m_{1}^{\ast}},n}\right)-\frac{2}{h}\trace\left(\hat{\Gam}_{\gam_{m_{1}},n}^{-1}\hat{W}_{\gam_{m_{1}},n}\right) \\
&=\frac{1}{h}\left(-\frac{1}{n}\p_{\gamma_{m_{1}}}^{2}\mbbh_{1,n}^{(m_{1})}(\tilde{\gamma}_{m_{1},n})\right)\left[\left\{\sqrt{T_{n}}\big(\hat{\gamma}_{m_{1},n}-f_{1}(\hat{\gamma}_{m_{1}^{\ast},n})\big)\right\}^{\otimes2}\right] \\
&\qquad+\frac{2}{h}\trace\left(\hat{\Gam}_{\gam_{m_{1}^{\ast}},n}^{-1}\hat{W}_{\gam_{m_{1}^{\ast}},n}\right)-\frac{2}{h}\trace\left(\hat{\Gam}_{\gam_{m_{1}},n}^{-1}\hat{W}_{\gam_{m_{1}},n}\right).
\end{align*}
Since the chain rule gives $\p_{\gamma_{m_{1}^{\ast}}}\mbbh_{1,n}^{(m_{1}^{\ast})}(\gamma_{m_{1}^{\ast},0})=F_{1}^{\top}\p_{\gamma_{m_{1}}}\mbbh_{1,n}^{(m_{1})}(\gamma_{m_{1},0})$ and $\p_{\gamma_{m_{1}^{\ast}}}^{2}\mbbh_{1,n}^{(m_{1}^{\ast})}(\gamma_{m_{1}^{\ast},0})=F_{1}^{\top}\p_{\gamma_{m_{1}}}^{2}\mbbh_{1,n}^{(m_{1})}(\gamma_{m_{1},0})F_{1}$, 
\begin{align*}
&\sqrt{T_{n}}\left\{f_{1}(\hat{\gamma}_{m_{1}^{\ast},n})-\gamma_{m_{1},0}\right\} \\
&=\sqrt{T_{n}}\left\{f_{1}(\hat{\gamma}_{m_{1}^{\ast},n})-f_{1}(\gamma_{m_{1}^{\ast},0})\right\} \\
&=F_{1}\sqrt{T_{n}}(\hat{\gamma}_{m_{1}^{\ast},n}-\gamma_{m_{1}^{\ast},0}) \\
&=F_{1}\left(-\frac{1}{n}\p_{\gamma_{m_{1}^{\ast}}}^{2}\mbbh_{1,n}^{(m_{1}^{\ast})}(\gamma_{m_{1}^{\ast},0})\right)^{-1}\left(\sqrt{\frac{h}{n}}\p_{\gamma_{m_{1}^{\ast}}}\mbbh_{1,n}^{(m_{1}^{\ast})}(\gamma_{m_{1}^{\ast},0})\right)+O_{p}(T_{n}^{-1/2}) \\
&=F_{1}\left[F_{1}^{\top}\left\{-\frac{1}{n}\p_{\gamma_{m_{1}}}^{2}\mbbh_{1,n}^{(m_{1})}\left(\gamma_{m_{1},0}\right)\right\}F_{1}\right]^{-1}F_{1}^{\top}\left(\sqrt{\frac{h}{n}}\p_{\gamma_{m_{1}}}\mbbh_{1,n}^{(m_{1})}(\gamma_{m_{1},0})\right)+O_{p}(T_{n}^{-1/2}) \\
&\cil F_{1}\left(F_{1}^{\top}\Gamma_{\gam_{m_{1}}}(\gam_{m_{1},0})F_{1}\right)^{-1}F_{1}^{\top}W_{\gamma_{m_{1}}}(\gam_{m_{1},0})^{1/2}{\bf N}_{m_{1}},
\end{align*}
where ${\bf N}_{m_{1}}\sim N_{p_{\gamma_{m_{1}}}}(0,I_{p_{\gamma_{m_{1}}}})$ (under $P$ without loss of generality).
Moreover, we have
\begin{align*}
\sqrt{T_{n}}\big(\hat{\gamma}_{m_{1},n}-f_{1}(\hat{\gamma}_{m_{1}^{\ast},n})\big)
&=\sqrt{T_{n}}\left\{(\hat{\gamma}_{m_{1},n}-\gamma_{m_{1},0})-\big(f_{1}(\hat{\gamma}_{m_{1}^{\ast},n})-\gamma_{m_{1},0}\big)\right\} \\
&\cil \left\{\Gamma_{\gam_{m_{1}}}(\gam_{m_{1},0})^{-1}-F_{1}\left(F_{1}^{\top}\Gamma_{\gam_{m_{1}}}(\gam_{m_{1},0})F_{1}\right)^{-1}F_{1}^{\top}\right\}W_{\gamma_{m_{1}}}(\gam_{m_{1},0})^{1/2}{\bf N}_{m_{1}} \\
&=G_{\gamma_{m_{1}}}(\gam_{m_{1},0})W_{\gamma_{m_{1}}}(\gam_{m_{1},0})^{1/2}{\bf N}_{m_{1}}.
\end{align*}
Hence, 
\begin{align*}
&\pr\left(\gqaic_{1,n}^{(m_{1}^{\ast})}-\gqaic_{1,n}^{(m_{1})}>0\right) \\
&=\pr\bigg[\frac{1}{h}\left(-\frac{1}{n}\p_{\gamma_{m_{1}}}^{2}\mbbh_{1,n}^{(m_{1})}(\tilde{\gamma}_{m_{1},n})\right)\left[\left(G_{\gamma_{m_{1}}}(\gam_{m_{1},0})W_{\gamma_{m_{1}}}(\gam_{m_{1},0})^{1/2}{\bf N}_{m_{1}}\right)^{\otimes2}\right] \\
&\quad\qquad+\frac{2}{h}\trace\left(\hat{\Gam}_{\gam_{m_{1}^{\ast}},n}^{-1}\hat{W}_{\gam_{m_{1}^{\ast}},n}\right)-\frac{2}{h}\trace\left(\hat{\Gam}_{\gam_{m_{1}},n}^{-1}\hat{W}_{\gam_{m_{1}},n}\right)>0\bigg] \\
&=\pr\bigg[\left(-\frac{1}{n}\p_{\gamma_{m_{1}}}^{2}\mbbh_{1,n}^{(m_{1})}(\tilde{\gamma}_{m_{1},n})\right)\left[\left(G_{\gamma_{m_{1}}}(\gam_{m_{1},0})W_{\gamma_{m_{1}}}(\gam_{m_{1},0})^{1/2}{\bf N}_{m_{1}}\right)^{\otimes2}\right]> \\
&\quad\qquad 2\trace\left(\hat{\Gam}_{\gam_{m_{1}},n}^{-1}\hat{W}_{\gam_{m_{1}},n}\right)-2\trace\left(\hat{\Gam}_{\gam_{m_{1}^{\ast}},n}^{-1}\hat{W}_{\gam_{m_{1}^{\ast}},n}\right)\bigg] \\
&\to\pr\bigg[\Gamma_{\gamma_{m_{1}}}(\gam_{m_{1},0})\left[\left(G_{\gamma_{m_{1}}}(\gam_{m_{1},0})W_{\gamma_{m_{1}}}(\gam_{m_{1},0})^{1/2}{\bf N}_{m_{1}}\right)^{\otimes2}\right] \\
&\quad\qquad >2\trace\left\{\Gam_{\gam_{m_{1}}}(\gam_{m_{1},0})^{-1}W_{\gam_{m_{1}}}(\gam_{m_{1},0})\right\}-2\trace\left\{\Gam_{\gam_{m_{1}^{\ast}}}(\gam_{m_{1}^{\ast},0})^{-1}W_{\gam_{m_{1}^{\ast}}}(\gam_{m_{1}^{\ast},0})\right\}\bigg] \\
&\to\pr\bigg[{\bf N}_{m_{1}}^{\top}W_{\gamma_{m_{1}}}(\gam_{m_{1},0})^{1/2}G_{\gamma_{m_{1}}}(\gam_{m_{1},0})W_{\gamma_{m_{1}}}(\gam_{m_{1},0})^{1/2}{\bf N}_{m_{1}} \\
&\quad\qquad >2\trace\left\{\Gam_{\gam_{m_{1}}}(\gam_{m_{1},0})^{-1}W_{\gam_{m_{1}}}(\gam_{m_{1},0})\right\}-2\trace\left\{\Gam_{\gam_{m_{1}^{\ast}}}(\gam_{m_{1}^{\ast},0})^{-1}W_{\gam_{m_{1}^{\ast}}}(\gam_{m_{1}^{\ast},0})\right\}\bigg]
\end{align*}
as $n\to\infty$.
%Since $W_{\gamma_{m_{1}}}(\gam_{m_{1},0})^{1/2}G_{\gamma_{m_{1}}}(\gam_{m_{1},0})\Gamma_{\gamma_{m_{1}}}(\gam_{m_{1},0})G_{\gamma_{m_{1}}}(\gam_{m_{1},0})W_{\gamma_{m_{1}}}(\gam_{m_{1},0})^{1/2}$ is the symmetric matrix,  there exists an orthogonal matrix $M$ such that 
%\begin{align*}
%&M^{\top}\left(W_{\gamma_{m_{1}}}(\gam_{m_{1},0})^{1/2}G_{\gamma_{m_{1}}}(\gam_{m_{1},0})\Gamma_{\gamma_{m_{1}}}(\gam_{m_{1},0})G_{\gamma_{m_{1}}}(\gam_{m_{1},0})W_{\gamma_{m_{1}}}(\gam_{m_{1},0})^{1/2}\right)M \\
%&=\diag(\lambda_{1},\lambda_{2},\ldots,\lambda_{p_{\gamma_{m_{1}}}}).
%\end{align*}
Since $W_{\gamma_{m_{1}}}(\gam_{m_{1},0})^{1/2}G_{\gamma_{m_{1}}}(\gam_{m_{1},0})W_{\gamma_{m_{1}}}(\gam_{m_{1},0})^{1/2}$ is the symmetric matrix,  there exists an orthogonal matrix $M$ such that 
\begin{align*}
&M^{\top}\left(W_{\gamma_{m_{1}}}(\gam_{m_{1},0})^{1/2}G_{\gamma_{m_{1}}}(\gam_{m_{1},0})W_{\gamma_{m_{1}}}(\gam_{m_{1},0})^{1/2}\right)M \\
&=\diag(\lambda_{1},\lambda_{2},\ldots,\lambda_{p_{\gamma_{m_{1}}}}).
\end{align*}
Therefore,
\begin{align*}
%&\Gamma_{\gamma_{m_{1}}}(\gam_{m_{1},0})\left[\left(G_{\gamma_{m_{1}}}(\gam_{m_{1},0})W_{\gamma_{m_{1}}}(\gam_{m_{1},0})^{1/2}{\bf N}_{m_{1}}\right)^{\otimes2}\right] \\
%&={\bf N}_{m_{1}}^{\top}W_{\gamma_{m_{1}}}(\gam_{m_{1},0})^{1/2}G_{\gamma_{m_{1}}}(\gam_{m_{1},0})\Gamma_{\gamma_{m_{1}}}(\gam_{m_{1},0})G_{\gamma_{m_{1}}}(\gam_{m_{1},0})W_{\gamma_{m_{1}}}(\gam_{m_{1},0})^{1/2}{\bf N}_{m_{1}} \\
&{\bf N}_{m_{1}}^{\top}W_{\gamma_{m_{1}}}(\gam_{m_{1},0})^{1/2}G_{\gamma_{m_{1}}}(\gam_{m_{1},0})W_{\gamma_{m_{1}}}(\gam_{m_{1},0})^{1/2}{\bf N}_{m_{1}} \\
&={\bf N}_{m_{1}}^{\top}M\diag(\lambda_{1},\lambda_{2},\ldots,\lambda_{p_{\gamma_{m_{1}}}})M^{\top}{\bf N}_{m_{1}} \\
&=\sum_{j=1}^{p_{\gamma_{m_{1}}}}\lambda_{j}\chi_{j}^{2}
\end{align*}
in distribution.

(ii) Since $\mathfrak{M}_{1}=\argmax_{m_{1}}\mbbh_{1,0}^{(m_{1})}(\gam_{m_{1}}^{\ast})$ holds, the inequality $\mbbh_{1,0}^{(m_{1})}(\gamma_{m_{1}}^{\ast})<\mbbh_{1,0}^{(m_{1}^{\ast})}(\gamma_{m_{1}^{\ast}}^{\ast})\big(=\mbbh_{1,0}^{(m_{1}^{\ast})}(\gamma_{m_{1}^{\ast},0})\big)$ is satisfied.
Further, we have
\begin{align*}
\frac{1}{n}\mbbh_{1,n}^{(m_{1})}(\hat{\gamma}_{m_{1},n})
&=\mbbh_{1,0}^{(m_{1})}(\gamma_{m_{1}}^{\ast})+o_{p}(1), \\
\frac{1}{n}\mbbh_{1,n}^{(m_{1}^{\ast})}(\hat{\gamma}_{m_{1}^{\ast},n})
&=\mbbh_{1,0}^{(m_{1}^{\ast})}(\gamma_{m_{1}^{\ast},0})+o_{p}(1).
\end{align*}
Hence,
\begin{align*}
&\pr\left(\mathrm{GQAIC}^{(m_{1}^{\ast})}-\mathrm{GQAIC}^{(m_{1})}>0\right) \\
&=\pr\left[-\frac{2}{n}\left(\mbbh_{1,n}^{(m_{1}^{\ast})}(\hat{\gamma}_{m_{1}^{\ast},n})-\mbbh_{1,n}^{(m_{1})}(\hat{\gamma}_{m_{1},n})\right)>\frac{2}{n}\left\{\trace\left(\hat{\Gam}_{\gam_{m_{1}},n}^{-1}\hat{W}_{\gam_{m_{1}},n}\right)-\trace\left(\hat{\Gam}_{\gam_{m_{1}^{\ast}},n}^{-1}\hat{W}_{\gam_{m_{1}^{\ast}},n}\right)\right\}\right] \\
&=\pr\left\{-2\left(\mbbh_{1,0}^{(m_{1}^{\ast})}(\gamma_{m_{1}^{\ast},0})-\mbbh_{1,0}^{(m_{1})}(\gamma_{m_{1}}^{\ast})\right)>o_{p}(1)\right\} \\
&\to0
\end{align*}
as $n\to\infty$.

\subsubsection{Proof of 2} \label{se:prf_thm_sele.prob1}
(i) The claims of 1 of this theorem mean that $\pr(\hat{m}_{1,n}\in\mathfrak{M}_{1})\to1$.
For any $m_{2}\in\mathfrak{M}_{2}\backslash\{m_{2}^{\ast}\}$, we have
\begin{align*}
&\pr\left(\gqaic_{2,n}^{(m_{2}^{\ast}|\hat{m}_{1,n})}-\gqaic_{2,n}^{(m_{2}|\hat{m}_{1,n})}>0\right) \\
&=\pr\left(\gqaic_{2,n}^{(m_{2}^{\ast}|\hat{m}_{1,n})}-\gqaic_{2,n}^{(m_{2}|\hat{m}_{1,n})}>0 \; | \; \hat{m}_{1,n}\in\mathfrak{M}_{1}\right)P\left(\hat{m}_{1,n}\in\mathfrak{M}_{1}\right) \\
&\qquad+\pr\left(\gqaic_{2,n}^{(m_{2}^{\ast}|\hat{m}_{1,n})}-\gqaic_{2,n}^{(m_{2}|\hat{m}_{1,n})}>0 \;| \; \hat{m}_{1,n}\notin\mathfrak{M}_{1}\right)P\left(\hat{m}_{1,n}\notin\mathfrak{M}_{1}\right) \\
&=\pr\left(\gqaic_{2,n}^{(m_{2}^{\ast}|\hat{m}_{1,n})}-\gqaic_{2,n}^{(m_{2}|\hat{m}_{1,n})}>0 \; | \; \hat{m}_{1,n}\in\mathfrak{M}_{1}\right)+o(1).
\end{align*}
Below, we focus on the cases where $\hat{m}_{1,n}\in\mathfrak{M}_{1}$.
Then, $\hat{\gam}_{\hat{m}_{1,n},n}$ converges to $\gam_{\hat{m}_{1,n},0}$ in probability, and $c_{\hat{m}_{1,n}}(\cdot,\gam_{\hat{m}_{1,n},0})=C(\cdot)$.

Because of assumptions, $a_{m_{2}}$ is correctly specified, and $\alpha_{m_{2}}^{\ast}=\alpha_{m_{2},0}$.
Define the map $f_{2}:\Theta_{\al_{m_{2}^{\ast}}}\to\Theta_{\al_{m_{2}}}$ by $f_{2}(\al_{m_{2}^{\ast}})=F_{2}\al_{m_{2}^{\ast}}+c_{2}$, where $F_{2}$ and $c_{2}$ satisfy the equation $\mbbh_{2,n}^{(m_{2}^{\ast}|m_{1})}(\al_{m_{2}^{\ast}})=\mbbh_{2,n}^{(m_{2}|m_{1})}\big(f_{2}(\al_{m_{2}^{\ast}})\big)$ for any $\al_{m_{2}^{\ast}}\in\Theta_{\al_{m_{2}^{\ast}}}$.
If $f_{2}(\al_{m_{2}^{\ast},0})=\al_{m_{2},0}$, then the inequality $\mbbh_{2,0}^{(m_{2}^{\ast}|m_{1})}(\al_{m_{2}^{\ast},0})=\mbbh_{2,0}^{(m_{2}|m_{1})}\big(f_{2}(\al_{m_{2}^{\ast},0})\big)<\mbbh_{2,0}^{(m_{2}|m_{1})}(\al_{m_{2},0})$ holds, and the assumption of $\mathfrak{M}_{2}$ is not satisfied. 
Hence, we have $f_{2}(\al_{m_{2}^{\ast},0})=\al_{m_{2},0}$.

Considering the fact that $P\{\p_{\al}\mbbh_{2,n}^{(m_{2}|\hat{m}_{1,n})}(\hat{\al}_{m_{2},n})=0\}\to1$ and Taylor expansion of $\mbbh_{2,n}^{(m_{2}|\hat{m}_{1,n})}$,
\begin{align*}
\mbbh_{2,n}^{(m_{2}^{\ast}|\hat{m}_{1,n})}(\hat{\al}_{m_{2}^{\ast},n})
&=\mbbh_{2,n}^{(m_{2}|\hat{m}_{1,n})}\big(f_{2}(\hat{\al}_{m_{2}^{\ast},n})\big) \\
&=\mbbh_{2,n}^{(m_{2}|\hat{m}_{1,n})}(\hat{\al}_{m_{2},n}) \\
&\qquad -\frac{1}{2}\left(-\frac{1}{T_{n}}\p_{\al_{m_{2}}}^{2}\mbbh_{2,n}^{(m_{2}|\hat{m}_{1,n})}(\tilde{\al}_{m_{2},n})\right)\left[\left\{\sqrt{T_{n}}\big(\hat{\al}_{m_{2},n}-f_{2}(\hat{\al}_{m_{2}^{\ast},n})\big)\right\}^{\otimes2}\right], 
%&=\mbbh_{2,n}^{(m_{2}|\hat{m}_{1,n})}(\hat{\al}_{m_{2},n}) \\
%&\qquad -\frac{1}{2}\left(-\frac{1}{T_{n}}\p_{\al_{m_{2}}}^{2}\mbbh_{2,n}^{(m_{2}|\hat{m}_{1,n})}(\tilde{\al}_{m_{2},n},\gam_{\hat{m}_{1,n},0})\right)\left[\left\{\sqrt{T_{n}}\big(\hat{\al}_{m_{2},n}-f_{2}(\hat{\al}_{m_{2}^{\ast},n})\big)\right\}^{\otimes2}\right]+O_{p}(T_{n}^{-1/2}),
\end{align*}
where $\tilde{\al}_{m_{2},n}=\hat{\al}_{m_{2},n}-\xi_{2}\left\{f_{2}(\hat{\al}_{m_{2}^{\ast},n})-\hat{\al}_{m_{2},n}\right\}$ for some $0<\xi_{2}<1$ and $\tilde{\al}_{m_{2},n}\cip\al_{m_{2},0}$ as $n\to\infty$. 
Moreover,
\begin{align*}
&\sqrt{T_{n}}(\hat{\al}_{m_{2},n}-\al_{m_{2},0}) \\
&=\left(-\frac{1}{T_{n}}\p_{\al_{m_{2}}}^{2}\mbbh_{2,n}^{(m_{2}|\hat{m}_{1,n})}(\al_{m_{2},0},\hat{\gam}_{\hat{m}_{1,n},n})\right)^{-1}\left(\frac{1}{\sqrt{T_{n}}}\p_{\al_{m_{2}}}\mbbh_{2,n}^{(m_{2}|\hat{m}_{1,n})}(\al_{m_{2},0},\hat{\gam}_{\hat{m}_{1,n},n})\right)+O_{p}(T_{n}^{-1/2}) \\
&=\left(-\frac{1}{T_{n}}\p_{\al_{m_{2}}}^{2}\mbbh_{2,n}^{(m_{2}|\hat{m}_{1,n})}(\al_{m_{2},0},\gam_{\hat{m}_{1,n},0})\right)^{-1}\left(\frac{1}{\sqrt{T_{n}}}\p_{\al_{m_{2}}}\mbbh_{2,n}^{(m_{2}|\hat{m}_{1,n})}(\al_{m_{2},0},\gam_{\hat{m}_{1,n},0}) \right. \\
&\qquad\qquad \left. -\frac{1}{T_{n}}\p_{\al_{m_{2}}}\p_{\gam_{\hat{m}_{1,n}}}\mbbh_{2,n}^{(m_{2}|\hat{m}_{1,n})}(\al_{m_{2},0},\gam_{\hat{m}_{1,n},0})\left[\sqrt{T_{n}}(\hat{\gam}_{\hat{m}_{1,n},n}-\gam_{\hat{m}_{1,n},0})\right]\right)+O_{p}(T_{n}^{-1/2}) \\
&\cil \Gam_{\al_{m_{2}}}(\al_{m_{2},0},\gam_{\hat{m}_{1,n},0})^{-1/2}{\bf N}_{m_{2}},
\end{align*}
where ${\bf N}_{m_{2}}\sim N_{p_{\al_{m_{2}}}}\left(0,I_{p_{\al_{m_{2}}}}\right)$. 
By the above and chain rule,
%gives $\p_{\al_{m_{2}^{\ast}}}\mbbh_{2,n}^{(m_{2}^{\ast}|\hat{m}_{1,n})}(\al_{m_{2}^{\ast},0})=F_{2}^{\top}\p_{\al_{m_{2}}}\mbbh_{2,n}^{(m_{2}|\hat{m}_{1,n})}(\al_{m_{2},0})$ and $\p_{\al_{m_{2}^{\ast}}}\mbbh_{2,n}^{(m_{2}^{\ast}|\hat{m}_{1,n})}(\al_{m_{2}^{\ast}})=F_{2}^{\top}\p_{\al_{m_{2}}}\mbbh_{2,n}^{(m_{2}|\hat{m}_{1,n})}(\al_{m_{2}})F_{2}$,
\begin{align*}
&\sqrt{T_{n}}\left\{\hat{\al}_{m_{2},n}-f_{2}(\hat{\al}_{m_{2}^{\ast},n})\right\} \\
&=\sqrt{T_{n}}\left\{(\hat{\al}_{m_{2},n}-\al_{m_{2},0})-\big(f_{2}(\hat{\al}_{m_{2}^{\ast},n})-f_{2}(\al_{m_{2}^{\ast},0})\big)\right\} \\
&=\sqrt{T_{n}}(\hat{\al}_{m_{2},n}-\al_{m_{2},0})-F_{2}\sqrt{T_{n}}(\hat{\al}_{m_{2}^{\ast},n}-\al_{m_{2}^{\ast},0}) \\
&=\sqrt{T_{n}}(\hat{\al}_{m_{2},n}-\al_{m_{2},0})-F_{2}\left(-\frac{1}{T_{n}}\p_{\al_{m_{2}^{\ast}}}^{2}\mbbh_{2,n}^{(m_{2}^{\ast}|\hat{m}_{1,n})}(\al_{m_{2}^{\ast},0},\gam_{\hat{m}_{1,n},0})\right)^{-1} \\
&\qquad\times\left(\frac{1}{\sqrt{T_{n}}}\p_{\al_{m_{2}^{\ast}}}\mbbh_{2,n}^{(m_{2}^{\ast}|\hat{m}_{1,n})}(\al_{m_{2}^{\ast},0},\gam_{\hat{m}_{1,n},0})\right)+O_{p}(T_{n}^{-1/2}) \\
&=\sqrt{T_{n}}(\hat{\al}_{m_{2},n}-\al_{m_{2},0})-F_{2}\left[F_{2}^{\top}\left\{-\frac{1}{T_{n}}\p_{\al_{m_{2}}}^{2}\mbbh_{2,n}^{(m_{2}|\hat{m}_{1,n})}\big(\al_{m_{2},0},\gam_{\hat{m}_{1,n},0}\big)\right\}F_{2}\right]^{-1} \\
&\qquad\times F_{2}^{\top}\left(\frac{1}{\sqrt{T_{n}}}\p_{\al_{m_{2}}}\mbbh_{2,n}^{(m_{2}|\hat{m}_{1,n})}(\al_{m_{2},0},\gam_{\hat{m}_{1,n},0})\right)+O_{p}(T_{n}^{-1/2}) \\
&\cil \left\{\Gam_{\al_{m_{2}}}(\al_{m_{2},0},\gam_{\hat{m}_{1,n},0})^{-1}-F_{2}\left(F_{2}^{\top}\Gam_{\al_{m_{2}}}(\al_{m_{2},0},\gam_{\hat{m}_{1,n},0})F_{2}\right)^{-1}F_{2}^{\top}\right\}\Gam_{\al_{m_{2}}}(\al_{m_{2},0},\gam_{\hat{m}_{1,n},0})^{1/2}{\bf N}_{m_{2}} \\
&=G_{\al_{m_{2}}}(\al_{m_{2},0},\gam_{\hat{m}_{1,n},0})\Gam_{\al_{m_{2}}}(\al_{m_{2},0},\gam_{\hat{m}_{1,n},0})^{1/2}{\bf N}_{m_{2}}.
\end{align*}
%where $G_{\al_{m_{2}}}(\al_{m_{2},0},\gam_{\hat{m}_{1,n},0})=\Gam_{\al_{m_{2}}}(\al_{m_{2},0},\gam_{\hat{m}_{1,n},0})^{-1}-F_{2}\left(F_{2}^{\top}\Gam_{\al_{m_{2}}}(\al_{m_{2},0},\gam_{\hat{m}_{1,n},0})F_{2}\right)^{-1}F_{2}^{\top}$.
Therefore,
\begin{align*}
&\pr\left(\gqaic_{2,n}^{(m_{2}^{\ast}|\hat{m}_{1,n})}-\gqaic_{2,n}^{(m_{2}|\hat{m}_{1,n})}>0\right) \\
&=\pr\left(\gqaic_{2,n}^{(m_{2}^{\ast}|\hat{m}_{1,n})}-\gqaic_{2,n}^{(m_{2}|\hat{m}_{1,n})}>0 \; | \; \hat{m}_{1,n}\in\mathfrak{M}_{1}\right)+o(1) \\
&=\pr\bigg[\left(-\frac{1}{T_{n}}\p_{\al_{m_{2}}}^{2}\mbbh_{2,n}^{(m_{2}|\hat{m}_{1,n})}(\tilde{\al}_{m_{2},n})\right)\left[\left\{\sqrt{T_{n}}\big(\hat{\al}_{m_{2},n}-f_{2}(\hat{\al}_{m_{2}^{\ast},n})\big)\right\}^{\otimes2}\right] \\
&\qquad+2p_{\al_{m_{2}^{\ast}}}-2p_{\al_{m_{2}}}>0 \; | \; \hat{m}_{1,n}\in\mathfrak{M}_{1}\bigg]+o(1) \\
&\to \pr\left[\Gam_{\al_{m_{2}}}(\al_{m_{2},0},\gam_{\hat{m}_{1,n},0})\left[\left(G_{\al_{m_{2}}}(\al_{m_{2},0},\gam_{\hat{m}_{1,n},0})\Gam_{\al_{m_{2}}}(\al_{m_{2},0},\gam_{\hat{m}_{1,n},0})^{1/2}{\bf N}_{m_{2}}\right)^{\otimes2}\right]>2(p_{\al_{m_{2}}}-p_{\al_{m_{2}^{\ast}}})\right] \\
&=\pr\bigg[{\bf N}_{m_{2}}^{\top}\Gam_{\al_{m_{2}}}(\al_{m_{2},0},\gam_{\hat{m}_{1,n},0})^{1/2}G_{\al_{m_{2}}}(\al_{m_{2},0},\gam_{\hat{m}_{1,n},0})\Gam_{\al_{m_{2}}}(\al_{m_{2},0},\gam_{\hat{m}_{1,n},0})^{1/2}{\bf N}_{m_{2}}>2(p_{\al_{m_{2}}}-p_{\al_{m_{2}^{\ast}}})\bigg]
\end{align*}
as $n\to\infty$.
In a similar way as (i) of Section \ref{se:prf_thm_sele.prob1}, we can show that
\begin{align*}
&{\bf N}_{m_{2}}^{\top}\Gam_{\al_{m_{2}}}(\al_{m_{2},0},\gam_{\hat{m}_{1,n},0})^{1/2}G_{\al_{m_{2}}}(\al_{m_{2},0},\gam_{\hat{m}_{1,n},0})\Gam_{\al_{m_{2}}}(\al_{m_{2},0},\gam_{\hat{m}_{1,n},0})^{1/2}{\bf N}_{m_{2}}=\sum_{j=1}^{p_{\al_{m_{2}}}}\lambda_{j}^{\prime}\chi_{j}^{2}
\end{align*}
in distribution.

(ii) The set $\mathfrak{M}_{2}$ satisfies that $\mathfrak{M}_{2}=\argmax_{m_{2}}\mbbh_{2,0}^{(m_{2}|\hat{m}_{1,n})}(\al_{m_{2}}^{\ast})$, so that the inequality $\mbbh_{2,0}^{(m_{2}|\hat{m}_{1,n})}(\al_{m_{2}}^{\ast})<\mbbh_{2,0}^{(m_{2}^{\ast}|\hat{m}_{1,n})}(\al_{m_{2}^{\ast}}^{\ast})\big(=\mbbh_{2,0}^{(m_{2}^{\ast}|\hat{m}_{1,n})}(\al_{m_{2}^{\ast},0})\big)$ holds.
Since
\begin{align*}
\frac{1}{T_{n}}\mbbh_{2,n}^{(m_{2}|\hat{m}_{1,n})}(\hat{\al}_{m_{2},n})
%&=\mbbh_{2,n}^{(m_{2}|\hat{m}_{1,n})}(\al_{m_{2}}^{\ast})+o_{p}(1)
&=\mbbh_{2,0}^{(m_{2}|\hat{m}_{1,n})}(\al_{m_{2}}^{\ast})+o_{p}(1), \\
\frac{1}{T_{n}}\mbbh_{2,n}^{(m_{2}^{\ast}|\hat{m}_{1,n})}(\hat{\al}_{m_{2}^{\ast},n})
%&=\mbbh_{2,n}^{(m_{2}^{\ast}|\hat{m}_{1,n})}(\al_{m_{2}^{\ast},0})+o_{p}(1)
&=\mbbh_{2,0}^{(m_{2}^{\ast}|\hat{m}_{1,n})}(\al_{m_{2}^{\ast},0})+o_{p}(1),
\end{align*}
we have
\begin{align*}
&\pr\left(\gqaic_{2,n}^{(m_{2}^{\ast}|\hat{m}_{1,n})}-\gqaic_{2,n}^{(m_{2}|\hat{m}_{1,n})}>0\right) \\
&=\pr\left\{-\frac{2}{T_{n}}\left(\mbbh_{2,n}^{(m_{2}^{\ast}|\hat{m}_{1,n})}(\hat{\al}_{m_{2}^{\ast},n})-\frac{1}{T_{n}}\mbbh_{2,n}^{(m_{2}|\hat{m}_{1,n})}(\hat{\al}_{m_{2},n})\right)>\frac{2}{T_{n}}(p_{\al_{m_{2}}}-p_{\al_{m_{2}^{\ast}}})\right\} \\
&=\pr\left\{-2\left(\mbbh_{2,0}^{(m_{2}^{\ast}|\hat{m}_{1,n})}(\al_{m_{2}^{\ast},0})-\mbbh_{2,0}^{(m_{2}|\hat{m}_{1,n})}(\al_{m_{2}}^{\ast})>o_{p}(1)\right)\right\} \\
&\to0
\end{align*}
as $n\to\infty$.

%%%%%
\subsection{Proof of Theorem \ref{se:thm_sele.prob2}}

(i) Since $h\log n \to 0$, in a similar way as the proof of Theorem \ref{se:thm_sele.prob} 1(i), we obtain
\begin{align*}
&\pr\left(\mathrm{GQBIC}_{1,n}^{\sharp(m_{1}^{\ast})}-\mathrm{GQBIC}_{1,n}^{\sharp(m_{1})}>0\right) \\
&=\pr\left[\left(-\frac{1}{n}\p_{\gamma_{m_{1}}}^{2}\mbbh_{1,n}^{(m_{1})}(\tilde{\gamma}_{m_{1},n})\right)\left[\left(G_{\gamma_{m_{1}}}(\gam_{m_{1},0})W_{\gamma_{m_{1}}}(\gam_{m_{1},0})^{1/2}{\bf N}_{m_{1}}\right)^{\otimes2}\right]>\left(\gam_{m_{1}}-\gam_{m_{1}^{\ast}}\right)h\log n\right] \\
&\to\pr\left[\Gamma_{\gamma_{m_{1}}}(\gam_{m_{1},0})\left[\left(G_{\gamma_{m_{1}}}(\gam_{m_{1},0})W_{\gamma_{m_{1}}}(\gam_{m_{1},0})^{1/2}{\bf N}_{m_{1}}\right)^{\otimes2}\right]>0\right] \\
&=1
\end{align*}
as $n\to\infty$.

(ii) As with Theorem \ref{se:thm_sele.prob} 1(ii), we can get
\begin{align*}
\pr\left(\mathrm{GQBIC}_{1,n}^{\sharp(m_{1}^{\ast})}-\mathrm{GQBIC}_{1,n}^{\sharp(m_{1})}>0\right)
&=\pr\left\{-2\left(\mbbh_{1,0}^{(m_{1}^{\ast})}(\gamma_{m_{1}^{\ast},0})-\mbbh_{1,0}^{(m_{1})}(\gamma_{m_{1}}^{\ast})\right)>\left(\gam_{m_{1}}-\gam_{m_{1}^{\ast}}\right)\frac{\log n}{n}\right\} \\
&\to0
\end{align*}
as $n\to\infty$.

%%%%%
\subsection{Proof of Theorem \ref{se:thm_sele.prob3}}
Since the proof of claim 2 can be handled analogously as claim 1 and Theorem \ref{se:thm_sele.prob} 2, we only prove claim 1.

(i) In a similar way as the proof of Theorem \ref{se:thm_sele.prob} 1(i), we have
\begin{align*}
&\pr\left(\mathrm{GQBIC}_{1,n}^{(m_{1}^{\ast})}-\mathrm{GQBIC}_{1,n}^{(m_{1})}>0\right) \\
&=\pr\left[\left(-\frac{1}{n}\p_{\gamma_{m_{1}}}^{2}\mbbh_{1,n}^{(m_{1})}(\tilde{\gamma}_{m_{1},n})\right)\left[\left(G_{\gamma_{m_{1}}}(\gam_{m_{1},0})W_{\gamma_{m_{1}}}(\gam_{m_{1},0})^{1/2}{\bf N}_{m_{1}}\right)^{\otimes2}\right]>\left(\gam_{m_{1}}-\gam_{m_{1}^{\ast}}\right)\log T_{n}\right] \\
&\to 0,
\end{align*}
as $n\to\infty$.

(ii) As with Theorem \ref{se:thm_sele.prob} 1(ii), we can deduce that
\begin{align*}
\pr\left(\mathrm{GQBIC}_{1,n}^{(m_{1}^{\ast})}-\mathrm{GQBIC}_{1,n}^{(m_{1})}>0\right)
&=\pr\left\{-2\left(\mbbh_{1,0}^{(m_{1}^{\ast})}(\gamma_{m_{1}^{\ast},0})-\mbbh_{1,0}^{(m_{1})}(\gamma_{m_{1}}^{\ast})\right)>\left(\gam_{m_{1}}-\gam_{m_{1}^{\ast}}\right)\frac{\log T_{n}}{T_{n}}\right\} \\
&\to0
\end{align*}
as $n\to\infty$.

%%%%%%
%\subsection{Proof of Theorem \ref{hm:thm_gqbic2}}

%%%%%
%%%%%
\bigskip

\noindent
\textbf{Acknowledgements.}
The authors are grateful to the anonymous reviewers for their valuable comments, which led to significant improvements in the first version.
This work was partially supported by JST CREST Grant Numbers JPMJCR14D7 and JPMJCR2115, and by JSPS KAKENHI Grant Numbers JP19K14593 and 22H01139.

%%\clearpage
\bigskip 
\bibliographystyle{abbrv} % plain,jplain,abbrv,unsrt,alpha,apalike
%\bibliography{SE_bibs}

\begin{thebibliography}{10}

\bibitem{Aka73}
H.~Akaike.
\newblock Information theory and an extension of the maximum likelihood
  principle.
\newblock In {\em Second {I}nternational {S}ymposium on {I}nformation {T}heory
  ({T}sahkadsor, 1971)}, pages 267--281. Akad\'emiai Kiad\'o, Budapest, 1973.

\bibitem{AsmRos01}
S.~Asmussen and J.~Rosi{\'n}ski.
\newblock Approximations of small jumps of {L}\'evy processes with a view
  towards simulation.
\newblock {\em J. Appl. Probab.}, 38(2):482--493, 2001.

\bibitem{BhaPap91}
R.~J. Bhansali and F.~Papangelou.
\newblock Convergence of moments of least squares estimators for the
  coefficients of an autoregressive process of unknown order.
\newblock {\em Ann. Statist.}, 19(3):1155--1162, 1991.

\bibitem{YUIMA14}
A.~Brouste, M.~Fukasawa, H.~Hino, S.~M. Iacus, K.~Kamatani, Y.~Koike,
  H.~Masuda, R.~Nomura, T.~Ogihara, Y.~Shimizu, M.~Uchida, and N.~Yoshida.
\newblock The yuima project: A computational framework for simulation and
  inference of stochastic differential equations.
\newblock {\em Journal of Statistical Software}, 57(4):1--51, 2014.

\bibitem{ChaIng11}
N.~H. Chan and C.-K. Ing.
\newblock Uniform moment bounds of {F}isher's information with applications to
  time series.
\newblock {\em Ann. Statist.}, 39(3):1526--1550, 2011.

\bibitem{CleGlo20}
E.~Cl\'{e}ment and A.~Gloter.
\newblock Joint estimation for {SDE} driven by locally stable {L}\'{e}vy
  processes.
\newblock {\em Electron. J. Stat.}, 14(2):2922--2956, 2020.

\bibitem{EguMas18a}
S.~Eguchi and H.~Masuda.
\newblock Schwarz type model comparison for {LAQ} models.
\newblock {\em Bernoulli}, 24(3):2278--2327, 2018.

\bibitem{EguMas19}
S.~Eguchi and H.~Masuda.
\newblock Data driven time scale in {G}aussian quasi-likelihood inference.
\newblock {\em Stat. Inference Stoch. Process.}, 22(3):383--430, 2019.

\bibitem{EguUeh21}
S.~Eguchi and Y.~Uehara.
\newblock Schwartz-type model selection for ergodic stochastic differential
  equation models.
\newblock {\em Scand. J. Stat.}, 48(3):950--968, 2021.

\bibitem{FinWei02}
D.~F. Findley and C.-Z. Wei.
\newblock A{IC}, overfitting principles, and the boundedness of moments of
  inverse matrices for vector autoregressions and related models.
\newblock {\em J. Multivariate Anal.}, 83(2):415--450, 2002.

\bibitem{IacYos18}
S.~M. Iacus and N.~Yoshida.
\newblock {\em Simulation and inference for stochastic processes with {YUIMA}}.
\newblock Use R! Springer, Cham, 2018.
\newblock A comprehensive R framework for SDEs and other stochastic processes.

\bibitem{JasKamMas19}
A.~Jasra, K.~Kamatani, and H.~Masuda.
\newblock Bayesian inference for stable {L}\'{e}vy-driven stochastic
  differential equations with high-frequency data.
\newblock {\em Scand. J. Stat.}, 46(2):545--574, 2019.

\bibitem{Ken82-2}
J.~T. Kent.
\newblock Robust properties of likelihood ratio tests.
\newblock {\em Biometrika}, 69(1):19--27, 1982.

\bibitem{KonKit96}
S.~Konishi and G.~Kitagawa.
\newblock Generalised information criteria in model selection.
\newblock {\em Biometrika}, 83(4):875--890, 1996.

\bibitem{LipShi01-1}
R.~S. Liptser and A.~N. Shiryaev.
\newblock {\em Statistics of random processes. {I}}, volume~5 of {\em
  Applications of Mathematics (New York)}.
\newblock Springer-Verlag, Berlin, expanded edition, 2001.
\newblock General theory, Translated from the 1974 Russian original by A. B.
  Aries, Stochastic Modelling and Applied Probability.

\bibitem{MagNeu79}
J.~R. Magnus and H.~Neudecker.
\newblock The commutation matrix: some properties and applications.
\newblock {\em Ann. Statist.}, 7(2):381--394, 1979.

\bibitem{Mas10ejs}
H.~Masuda.
\newblock Approximate self-weighted {LAD} estimation of discretely observed
  ergodic {O}rnstein-{U}hlenbeck processes.
\newblock {\em Electron. J. Stat.}, 4:525--565, 2010.

\bibitem{Mas13}
H.~Masuda.
\newblock Convergence of {G}aussian quasi-likelihood random fields for ergodic
  {L}\'evy driven {SDE} observed at high frequency.
\newblock {\em Ann. Statist.}, 41(3):1593--1641, 2013.

\bibitem{Mas19spa}
H.~Masuda.
\newblock Non-{G}aussian quasi-likelihood estimation of {SDE} driven by locally
  stable {L}\'{e}vy process.
\newblock {\em Stochastic Process. Appl.}, 129(3):1013--1059, 2019.

\bibitem{masuda2023optimal}
H.~Masuda.
\newblock Optimal stable {Ornstein}-{Uhlenbeck} regression.
\newblock {\em Japanese Journal of Statistics and Data Science}, accepted,
  2023.

\bibitem{MasUeh17-2}
H.~Masuda and Y.~Uehara.
\newblock On stepwise estimation of {L}\'evy driven stochastic differential
  equation ({Japanese}).
\newblock {\em Proc. Inst. Statist. Math.}, 65(1):21--38, 2017.

\bibitem{OgiYos11}
T.~Ogihara and N.~Yoshida.
\newblock Quasi-likelihood analysis for the stochastic differential equation
  with jumps.
\newblock {\em Stat. Inference Stoch. Process.}, 14(3):189--229, 2011.

\bibitem{Rad08}
P.~Radchenko.
\newblock Mixed-rates asymptotics.
\newblock {\em Ann. Statist.}, 36(1):287--309, 2008.

\bibitem{Sch78}
G.~Schwarz.
\newblock Estimating the dimension of a model.
\newblock {\em Ann. Statist.}, 6(2):461--464, 1978.

\bibitem{Uch10}
M.~Uchida.
\newblock Contrast-based information criterion for ergodic diffusion processes
  from discrete observations.
\newblock {\em Ann. Inst. Statist. Math.}, 62(1):161--187, 2010.

\bibitem{UchYos12}
M.~Uchida and N.~Yoshida.
\newblock Adaptive estimation of an ergodic diffusion process based on sampled
  data.
\newblock {\em Stochastic Process. Appl.}, 122(8):2885--2924, 2012.

\bibitem{Ueh19}
Y.~Uehara.
\newblock Statistical inference for misspecified ergodic {L}\'{e}vy driven
  stochastic differential equation models.
\newblock {\em Stochastic Process. Appl.}, 129(10):4051--4081, 2019.

\bibitem{Wat13}
S.~Watanabe.
\newblock A widely applicable {B}ayesian information criterion.
\newblock {\em J. Mach. Learn. Res.}, 14:867--897, 2013.

\bibitem{ShiYos06}
S.~Yasutaka and Y.~Nakahiro.
\newblock Estimation of parameters for diffusion processes with jumps from
  discrete observations.
\newblock {\em Statistical Inference and Stochastic Processes}, 9(3):227--277,
  2006.

\bibitem{Yos11}
N.~Yoshida.
\newblock Polynomial type large deviation inequalities and quasi-likelihood
  analysis for stochastic differential equations.
\newblock {\em Ann. Inst. Statist. Math.}, 63(3):431--479, 2011.

\end{thebibliography}

\end{document}